%% file: II_genfeedforward.tex
\documentclass[11pt, a4paper, final]{article} 			
\usepackage[english]{babel}
\usepackage[utf8]{inputenc}  							
\usepackage{amsmath}
\usepackage{amssymb}
\usepackage{amsthm}
\usepackage{mathtools}
\usepackage{txfonts}
\usepackage{color}
\usepackage{enumerate}									
\usepackage[inline]{enumitem}
\usepackage{xspace}										
\usepackage{caption}
\usepackage{multirow}
\usepackage{authblk}									
\usepackage{float}
\usepackage{bbm}										
\usepackage{etoolbox}									
\usepackage{subcaption}
\usepackage{ifdraft}

\usepackage{standalone}									
\usepackage{xcolor}										
\usepackage{tikz}										
\usetikzlibrary{%
	arrows,
	shapes.geometric,
	arrows.meta,
	fit,
	positioning,
	calc
}
\usepackage{pgfplots}
\usepgfplotslibrary{colormaps}

\usepackage[
	colorinlistoftodos,
	textsize = small,
	color=blue!20!white,
	bordercolor=red
	]{todonotes}										

\usepackage{ragged2e}									
\usepackage[
	style = numeric,
	sorting = nyt,
	sortcites,
	firstinits = true,
	isbn = false,
	doi = false, 
	clearlang = true,
	backend = biber
	]{biblatex}
\AtEveryBibitem{\clearfield{pagetotal}}					
\AtEveryBibitem{\clearfield{note}}						
\DeclareSourcemap{
	\maps[datatype=bibtex]{
		\map{
			\pernottype{misc} 
			\step[fieldset=url, null] 
			\step[fieldset=urldate, null] 
		}
	}
}

\usepackage{csquotes}	
\DeclareRedundantLanguages{English,english,german,french,eng}{english,german,ngerman,french} 

\renewbibmacro{in:}{%
	\ifentrytype{article}{}{\printtext{\bibstring{in}\intitlepunct}}}		

\renewcommand{\textcite}{\cite}							

\usepackage{geometry}
\usepackage{aliascnt}
\usepackage[pdfpagelabels=true,plainpages=false, hidelinks]{hyperref}
\usepackage{cleveref}									

\theoremstyle{plain}
\newtheorem{theorem}{Theorem}[section]

\newtheorem*{theorem*}{Theorem}

\newaliascnt{lemma}{theorem}
\newtheorem{lemma}[lemma]{Lemma}
\aliascntresetthe{lemma}
\crefname{lemma}{lemma}{lemmas}

\newaliascnt{cor}{theorem}
\newtheorem{cor}[cor]{Corollary}
\aliascntresetthe{cor}

\newaliascnt{prop}{theorem}
\newtheorem{prop}[prop]{Proposition}
\aliascntresetthe{prop}
\crefname{prop}{proposition}{propositions}

\theoremstyle{definition}
\newtheorem*{defi}{Definition}

\theoremstyle{remark}
\newaliascnt{remk}{theorem}
\newtheorem{remk}[remk]{Remark}
\aliascntresetthe{remk}
\crefname{remk}{remark}{remarks}

\newtheorem{ex}{Example}[section]

\newtheorem{step}{Step}
\crefname{step}{step}{steps}
\AtBeginEnvironment{proof}{\setcounter{step}{0}}
\crefformat{step}{\textit{step}~#2#1#3}
\Crefformat{step}{\textit{Step}~#2#1#3}
\crefrangeformat{step}{\textit{steps}~#3#1#4 to~#5#2#6}
\Crefrangeformat{step}{\textit{Steps}~#3#1#4 to~#5#2#6}
\crefmultiformat{step}{\textit{steps}~#2#1#3}{ and~#2#1#3}{, #2#1#3}{ and~#2#1#3}
\Crefmultiformat{step}{\textit{Steps}~#2#1#3}{ and~#2#1#3}{, #2#1#3}{ and~#2#1#3}

\newcommand{\RR}{\mathbb{R}}

\newcommand{\OO}{\mathcal{O}}

\newcommand{\LL}{\mathcal L}

\newcommand{\gl}[1]{\mathfrak{gl}(#1)}
\newcommand{\spec}[1]{\operatorname{spec(#1)}}

\newcommand{\nole}{\trianglelefteq}
\newcommand{\noge}{\trianglerighteq}
\newcommand{\nol}{\vartriangleleft}
\newcommand{\nog}{\vartriangleright}

\newcommand{\loopeq}{\multimapboth}

\newcommand{\idv}{\mathbbm{1}_V}

\newcommand{\dsup}{d}

\newcommand{\Rsup}{R}

\newcommand{\xisup}{\Xi}
\newcommand{\xisub}{\Xi}

\newcommand{\musup}{\mu}
\newcommand{\musub}{\mu}
\newcommand{\Mu}{\hat{\mu}}
\newcommand{\Qsup}{Q}

\newcommand{\Qn}{\mathbf{Q}}
\newcommand{\Dsup}{D}
\newcommand{\Dsub}{D}
\newcommand{\Bsup}{B}
\newcommand{\Bsub}{B}

\newcommand{\SSigma}{\Sigma^\star}

\DeclareMathOperator{\Id}{Id}

\definecolor{wred}{rgb}{0.7,0.18,0.12}
\definecolor{wgreen}{rgb}{0.1,0.53,0.37}

\numberwithin{equation}{section} 
\renewcommand{\theequation}{\arabic{section}.\arabic{equation}} 

\definecolor{red}{RGB}{226,0,26}									
\definecolor{blue}{RGB}{0,156,209}									
\definecolor{grey}{RGB}{59,81,91}									

\begin{document}

\title{Amplified steady state bifurcations in feedforward networks}
\author{Sören von der Gracht\thanks{Department of Mathematics, Universität Hamburg, Germany, \href{mailto:soeren.von.der.gracht@uni-hamburg.de}{soeren.von.der.gracht@uni-hamburg.de}}, Eddie Nijholt\thanks{\mbox{ICMC, University of São Paulo, São Carlos, Brazil, \href{mailto:eddie.nijholt@gmail.com}{eddie.nijholt@gmail.com}}}, Bob Rink\thanks{\mbox{Department of Mathematics, Vrije Universiteit Amsterdam, The Netherlands, \href{mailto:b.w.rink@vu.nl}{b.w.rink@vu.nl}}}}
\date{}
\maketitle

\input{II_genfeedforward_in.tex}

\input{II_genfeedforward_pr.tex}

\input{II_genfeedforward_td.tex}

\input{II_genfeedforward_bi.tex}

\input{II_genfeedforward_ex.tex}

\section*{Acknowledgement}
Parts of this work originated in Sören von der Gracht's doctoral project and are contained in his thesis (``Genericity in Network Dynamics'', 2019 \cite{Schwenker.2019}), written under the primary supervision of Reiner Lauterbach (Universität Hamburg) and co-examined by Bob Rink (Vrije Universiteit Amsterdam) and Ana Paula Dias (Universidade do Porto). The author wishes to express his gratitude to the examiners for helpful comments, discussions and support.\\

\vspace{.5cm}
\noindent
\begin{minipage}{0.85\textwidth}
	\begin{flushleft}
		This research is partly financed by the Dutch Research Council (NWO) via Eddie Nijholt's research program ``Designing Network Dynamical Systems through Algebra''.
	\end{flushleft}
\end{minipage}\hfill
\begin{minipage}{0.1\textwidth}
	\hfill \includegraphics[width=1cm]{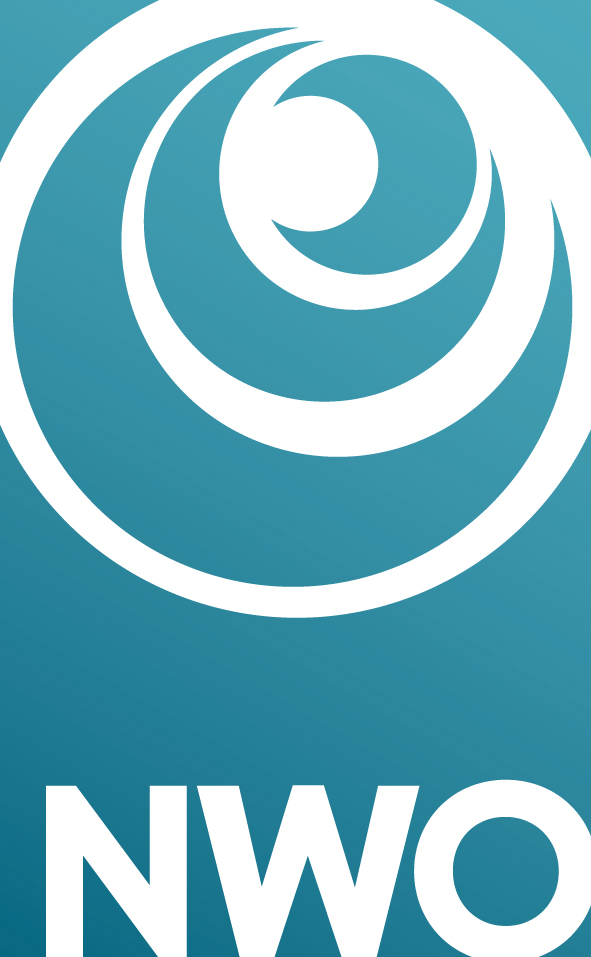}
\end{minipage}\\

\vspace{.5cm}
\noindent
Bob Rink is happy to acknowledge the hospitality and financial support of the Sydney Mathematical Research Institute.

\input{II_genfeedforward_ap.tex}

\begingroup
\RaggedRight
\printbibliography
\endgroup
\phantomsection

\end{document}

%% file: II_genfeedforward_in.tex
\begin{abstract}
	\noindent
	We investigate bifurcations in feedforward coupled cell networks. Feedforward structure (the absence of feedback) can be defined by a partial order on the cells. We use this property to study generic one-parameter steady state bifurcations for such networks. Branching solutions and their asymptotics are described in terms of Taylor coefficients of the internal dynamics. They can be determined via an algorithm that only exploits the network structure. Similar to previous results on feedforward chains, we observe amplifications of the growth rates of steady state branches induced by the feedforward structure. However, contrary to these earlier results, as the interaction scenarios can be more complicated in general feedforward networks, different branching patterns and different amplifications can occur for different regions in the space of Taylor coefficients.
\end{abstract}

\section*{Introduction}

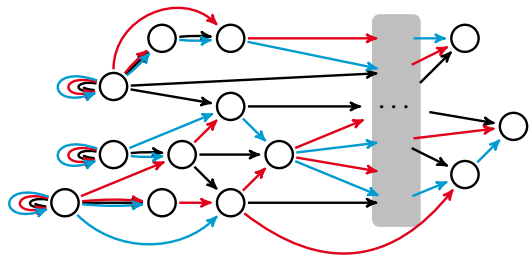
\begin{figure}[b]
	\begin{center}
		\resizebox{.475\linewidth}{!}{
			\begin{tikzpicture}[->,
				>=stealth',
				shorten >=1pt,
				auto,
				node distance=.5cm,
				main node/.style={line width=1.5pt, circle, scale = 3, draw, font=\sffamily\tiny, inner sep=2pt}
				]
				\node[main node] (1) {};
				\node[main node, above of=1] (2) {};
				\node[main node, below left of=1] (3) {};
				\node[main node, above right of=2] (4) {};
				\node[main node, right of=1] (5) {};
				\node[main node, above right of=5] (6) {};
				\node[main node, below right of=1] (7) {};
				\node[main node, below right of=6] (8) {};
				\node[main node, below right of=5] (9) {};
				\node[main node, right of=4] (10) {};
				
				\node[main node, above right of=8,draw=none] (c1) {};
				\node[main node, right of=c1,draw=none] (d1) {$\dots$};
				\node[main node, above of=d1,draw=none, node distance=.25cm] (d2) {};
				\node[main node, above of=d1,draw=none] (d3) {};
				\node[main node, below of=d1,draw=none, node distance=.25cm] (d4) {};
				\node[main node, below of=d1,draw=none] (d5) {};
				\coordinate[main node, below right of=8,draw=none] (c2) {};
				\node[main node, right of=c2,draw=none] (d6) {};
				
				\node[main node, right of=d3] (11) {};
				\node[main node, right of=d5] (12) {};
				\node[main node, above right of=12] (13) {};
				
				\node[line width=1.5pt, color=lightgray, draw, fill, rectangle, rounded corners, fit=(d3)(d6)](ld) {};
				\node[main node, right of=c1,draw=none] (d) {$\dots$};
				
				\path[every node/.style={font=\sffamily\small}, line width =1.5pt]
				(1) edge [loop left, in=190, out=170, looseness=10] node {} (1)
				(2) edge [loop left, in=190, out=170, looseness=10] node {} (2)
				(3) edge [loop left, in=190, out=170, looseness=10] node {} (3)
				(2) edge [] node {} (4)
				(1) edge [bend left=5] node {} (5)
				(2) edge [] node {} (6)
				(3) edge [] node {} (7)
				(5) edge [] node {} (8)
				(5) edge [] node {} (9)
				(4) edge [bend left=5] node {} (10)
				(6) edge [] node {} (d1)
				(2) edge [] node {} (d2)
				(9) edge [] node {} (d6)
				(d1) edge [] node {} (11)
				(d4) edge [] node {} (12)
				(d1) edge [] node {} (13)
				
				(1) edge [loop left, color=red, in=195, out=165, looseness=11] node {} (1)
				(2) edge [loop left, color=red, in=195, out=165, looseness=11] node {} (2)
				(3) edge [loop left, color=red, in=195, out=165, looseness=11] node {} (3)
				(2) edge [color=red, bend left=5] node {} (4)
				(3) edge [color=red] node {} (5)
				(5) edge [color=red] node {} (6)
				(3) edge [color=red, bend left=5] node {} (7)
				(9) edge [color=red] node {} (8)
				(7) edge [color=red] node {} (9)
				(2) edge [color=red, out=90, in=135, looseness=1.2] node {} (10)
				(8) edge [color=red] node {} (d1)
				(10) edge [color=red] node {} (d3)
				(8) edge [color=red] node {} (d5)
				(d2) edge [color=red] node {} (11)
				(9) edge [color=red, bend left=5, out=315, in=225, looseness=1.2] node {} (12)
				(d4) edge [color=red] node {} (13)
				
				(1) edge [loop left, color=blue, in=200, out=160, looseness=12] node {} (1)
				(2) edge [loop left, color=blue, in=200, out=160, looseness=12] node {} (2)
				(3) edge [loop left, color=blue, in=200, out=160, looseness=12] node {} (3)
				(2) edge [color=blue, bend right=5] node {} (4)
				(1) edge [color=blue, bend right=5] node {} (5)
				(1) edge [color=blue] node {} (6)
				(3) edge [color=blue, bend right=5] node {} (7)
				(6) edge [color=blue] node {} (8)
				(3) edge [color=blue, out=-45, in=-135] node {} (9)
				(4) edge [color=blue, bend right=5] node {} (10)
				(10) edge [color=blue] node {} (d2)
				(8) edge [color=blue] node {} (d4)
				(8) edge [color=blue] node {} (d6)
				(d3) edge [color=blue] node {} (11)
				(d6) edge [color=blue] node {} (12)
				(12) edge [color=blue] node {} (13)
				;
			\end{tikzpicture}
		}
		\caption{A general homogeneous feedforward network with asymmetric inputs.}
		\label{fig:genff}
	\end{center}
\end{figure}

\paragraph{Summary of main results.}
We investigate \emph{feedforward} network dynamical systems for their bifurcation behavior. Most generally, feedforward structure is defined by the absence of feedback except for self-loops, i.e., there are no directed cycles consisting of $2$ or more cells. Under this structural assumption, we classify generic steady state bifurcations of \emph{homogeneous} networks with \emph{asymmetric inputs} (compare to \Cref{fig:genff}) by proving
\begin{theorem*}[\Cref{thm:maxcrit,thm:abff_maxnoncritsuper,thm:abff_maxnoncritsub,thm:abff_maxnoncritgen}]
	For each branch of steady state solutions emerging in a generic one-parameter steady state bifurcation in a feedforward network, there exists a unique \emph{root subnetwork}, i.e., a subnetwork that is surrounded by \emph{critical} cells, in which the cells remain synchronous. The state of a cell that is not in the root subnetwork grows asymptotically in the bifurcation parameter as $\sim|\lambda|^{\frac{1}{n}}$ where $n$ is given by the maximal number of critical cells along paths from the root subnetwork to this cell. This effect is called \emph{amplification}.
\end{theorem*}
\noindent
The result generalizes earlier classifications of generic steady state branches for feedforward chains and layered feedforward networks. We elaborate on this relation in the remainder of this introduction after illustrating the main result in three numerically investigated examples.

\paragraph{Some examples.}
\begin{figure}[t]
	\centering
	\resizebox{.5\linewidth}{!}{	
		\input{5cellfeedforwardwithoddbifurcationsextnet}%
	}%
	\\[\smallskipamount]
	\includegraphics[width=.491\textwidth]{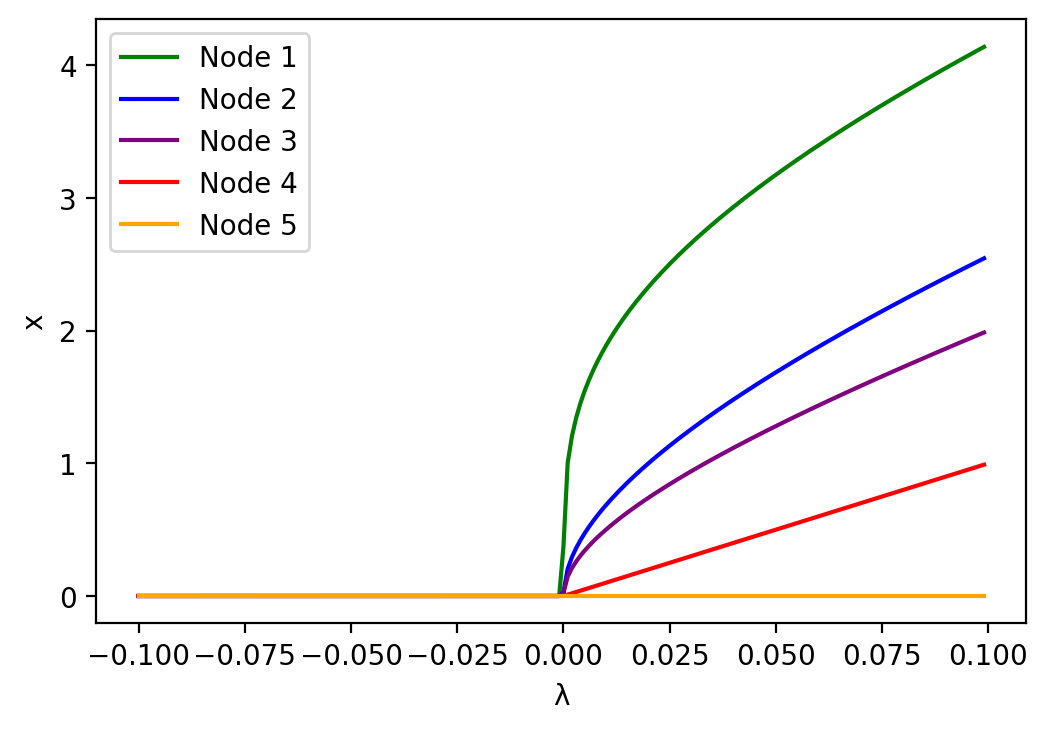}\hfill
	\includegraphics[width=.501\textwidth]{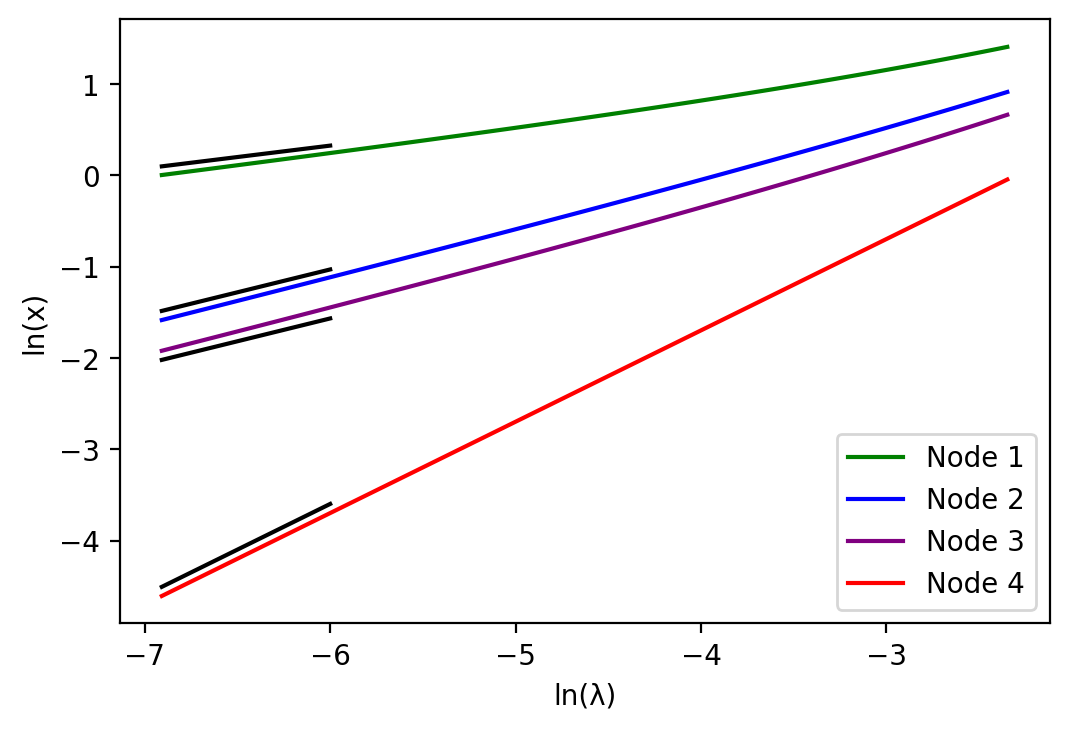}
	\caption{An amplifying bifurcation in a feedforward network. Top: a graphical depiction of the network. Bottom left: numerical evidence of a bifurcation in which a synchrony-breaking branch of steady state points emerges. Bottom right: a log–log plot of the synchrony-breaking branch on the left for positive values of $\lambda$. Black line-segments have fixed slopes $1/4, 1/2, 1/2$ and $1$, from top to bottom, indicating asymptotics of the different nodes.}\label{fig:numerics1}
\end{figure}
Consider the feedforward network shown at the top of \Cref{fig:numerics1}. A general family of admissible vector fields for this network is given by 
\begin{equation}
	\label{eq:numerics1} 
	\gamma_{f} (x) =
	\begin{pmatrix}
		f(x_{1},\textcolor{red}{x_{2}},\textcolor{blue}{x_{3}},\textcolor{grey}{x_{4}},\textcolor{magenta}{x_{5}}, \lambda)\\
		f(x_{2},\textcolor{red}{x_{5}},\textcolor{blue}{x_{4}},\textcolor{grey}{x_{5}},\textcolor{magenta}{x_{5}}, \lambda)\\
		f(x_{3},\textcolor{red}{x_{4}},\textcolor{blue}{x_{5}},\textcolor{grey}{x_{5}},\textcolor{magenta}{x_{5}}, \lambda)\\
		f(x_{4},\textcolor{red}{x_{5}},\textcolor{blue}{x_{5}},\textcolor{grey}{x_{5}},\textcolor{magenta}{x_{5}}, \lambda)\\
		f(x_{5},\textcolor{red}{x_{5}},\textcolor{blue}{x_{5}},\textcolor{grey}{x_{5}},\textcolor{magenta}{x_{5}}, \lambda)\\
	\end{pmatrix}.
\end{equation}
where $\lambda \in \mathbb{R}$ is a bifurcation parameter and the colors of the variables correspond to the different arrow-types in the figure. In the graph we have left out self-loops corresponding to the first argument of $f$, as they are understood as internal dynamics.

The bottom left of Figure \ref{fig:numerics1} shows a numerically computed bifurcation branch in a system of the form \eqref{eq:numerics1}.  More precisely, we have used a specific choice of response function $f: \mathbb{R}^4 \times \mathbb{R} \rightarrow \mathbb{R}$, given by 
\begin{equation*}
	f(x,y,z,v,w, \lambda) = y+2z-4w+5\lambda x -0.5x^2
\end{equation*}
(strictly speaking the fact that $f$ does not depend on $v$ cancels the influence of the \textcolor{grey}{grey} arrows). For each of $200$ fixed values of $\lambda \in [-0.1, 0.1]$, we forward integrated the system of Equation  \eqref{eq:numerics1} up to $t = 10000$, using Euler's method with time steps of $0.1$. The values of $\lambda$ are evenly spaced and the integration started each time from $(0.01, 0.02, 0.03, 0.04, -0.05)$. 
The figure shows all five coordinates of the final point of integration. For negative values of $\lambda$, all components end up indistinguishably close to zero, indicating that the origin is stable. For positive values of $\lambda$, however, stability is passed on to a fully non-synchronous point.

Writing $(x_1(\lambda), \dots, x_5(\lambda))$ for this non-synchronous branch of steady state points, the figure suggests that $x_5(\lambda) = 0$ and $x_4(\lambda) \sim \lambda$. The components $x_3(\lambda)$ and $x_2(\lambda)$ seem to grow at a steeper rate, hinting at the amplifying behavior $x_2(\lambda), x_3(\lambda)\sim \lambda^{1/2}$ often observed in feedforward structures. Finally, $x_1(\lambda)$ grows at the steepest rate, which might indicate $x_1(\lambda) \sim \lambda^{1/4}$.

The bottom right of \cref{fig:numerics1} corroborates these suggested growth rates for $x_1(\lambda)$ to $x_4(\lambda)$. Shown here is a plot of $\ln(\lambda)$ against $\ln(x_i)$ for $i \in \{1, \dots, 4\}$, computed as before for $200$ evenly spaced values of $\ln(\lambda) \in [\ln(0.001), \ln(0.1)]$. The black line-segments have fixed slopes $1/4, 1/2, 1/2$ and $1$, from top to bottom. This suggests that we indeed have the asymptotics $x_4(\lambda) \sim \lambda$, $x_3(\lambda), x_2(\lambda) \sim \lambda^{1/2}$ and $x_1(\lambda) \sim \lambda^{1/4}$. Our results in this paper predict this bifurcation in a system of the form \eqref{eq:numerics1}  for an open set of Taylor coefficients of $f$. Additional (unstable) branches are furthermore typically present. We investigate the bifurcations in the network in \Cref{fig:numerics1} analytically in \Cref{sec:ex}.

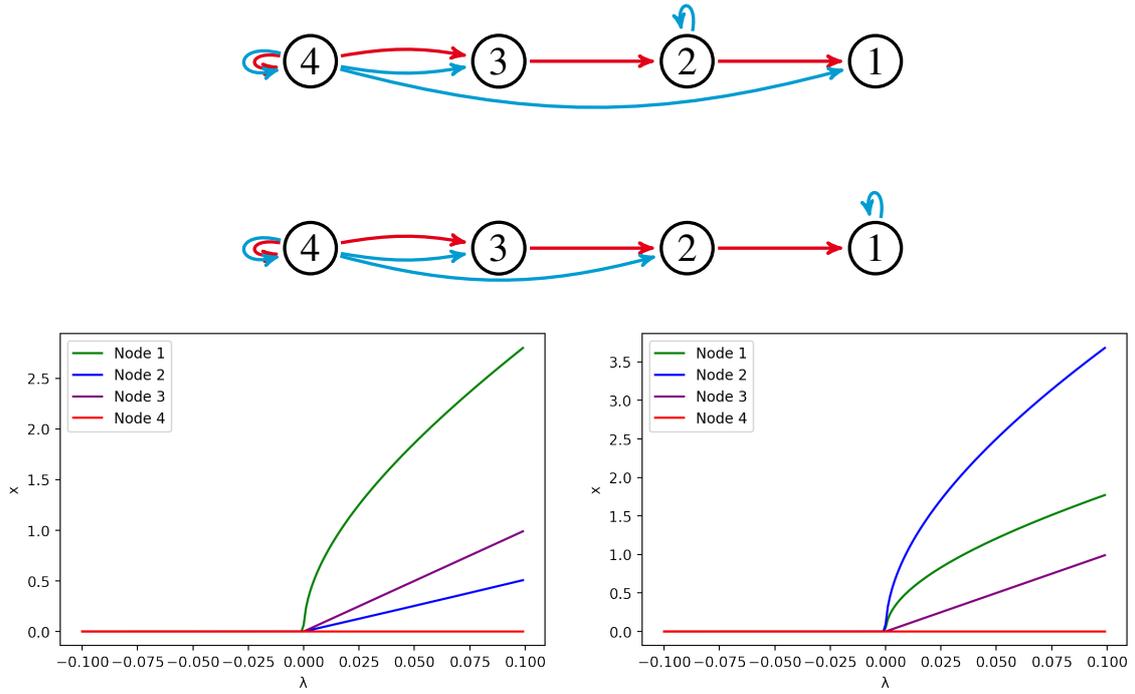
\begin{figure}[t]
	\centering
	\resizebox{.6\linewidth}{!}{
		\begin{tikzpicture}[->,
			>=stealth',
			shorten >=1pt,
			auto,
			node distance=1cm,
			main node/.style={line width=1.5pt, circle, scale = 3, draw, font=\sffamily\tiny, inner sep=1pt}]
			\node[main node] (1) {$1$};
			\node[main node, left of=1] (2) {$2$};
			\node[main node, left of=2] (3) {$3$};
			\node[main node, left of=3] (4) {$4$};
			\path[every node/.style={font=\sffamily\small}, line width =1.5pt]
			(2) edge [color = {red}] (1)
			(3) edge [color = {red}] (2)
			(4) edge [color = {red}, bend left=10] (3)
			(4) edge [in=190, out=170, looseness = 8, color = {red}] (4)
			(4) edge [color = {blue}, bend right=15] (1)
			(2) edge [in=100, out=80, looseness = 8, color = {blue}] (2)
			(4) edge [color = {blue}, bend right=10] (3)
			(4) edge [in=195, out=165, looseness = 8, color = {blue}] (4)
			;
			
			\node[main node, below of=1] (1s) {$1$};
			\node[main node, below of=2] (2s) {$2$};
			\node[main node, below of=3] (3s) {$3$};
			\node[main node, below of=4] (4s) {$4$};
	
			\path[every node/.style={font=\sffamily\small}, line width =1.5pt]
			(2s) edge [color = {red}] (1s)
			(3s) edge [color = {red}] (2s)
			(4s) edge [color = {red}, bend left=10] (3s)
			(4s) edge [in=190, out=170, looseness = 8, color = {red}] (4s)
			
			(1s) edge [in=100, out=80, looseness = 8, color = {blue}] (1s)
			(4s) edge [color = {blue}, bend right=15] (2s)
			(4s) edge [color = {blue}, bend right=10] (3s)
			(4s) edge [in=195, out=165, looseness = 8, color = {blue}] (4s)
			;
			
		\end{tikzpicture}
	}
	\\[\bigskipamount]
	\includegraphics[width=.49\textwidth]{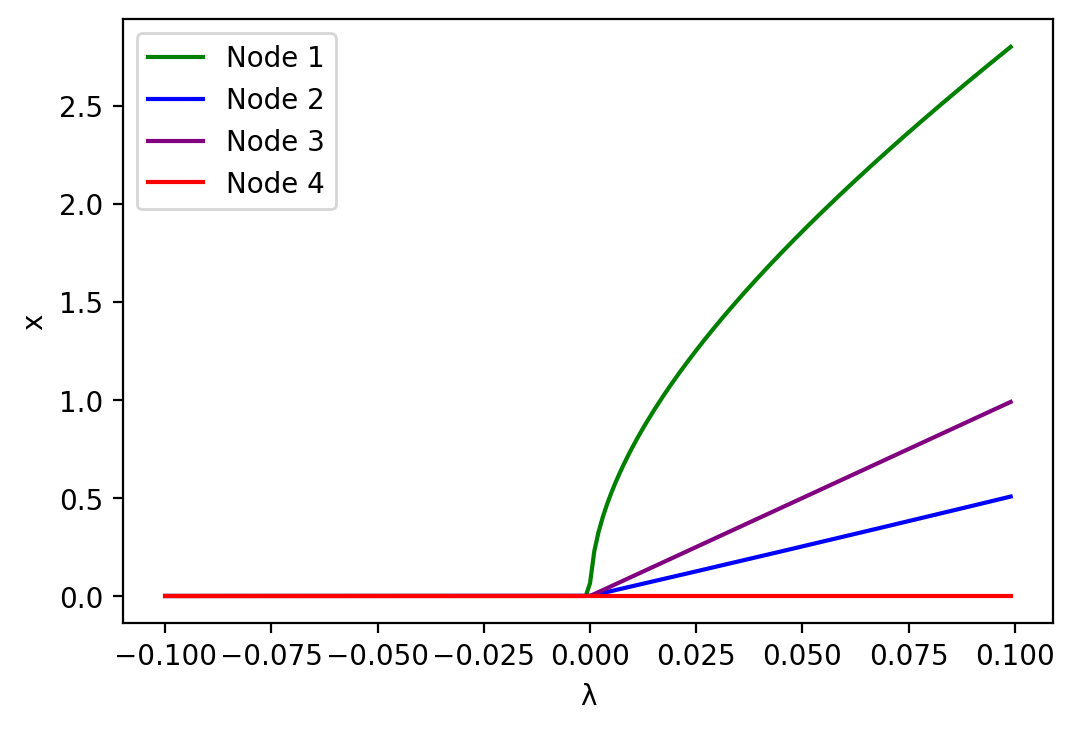}\hfill
	\includegraphics[width=.49\textwidth]{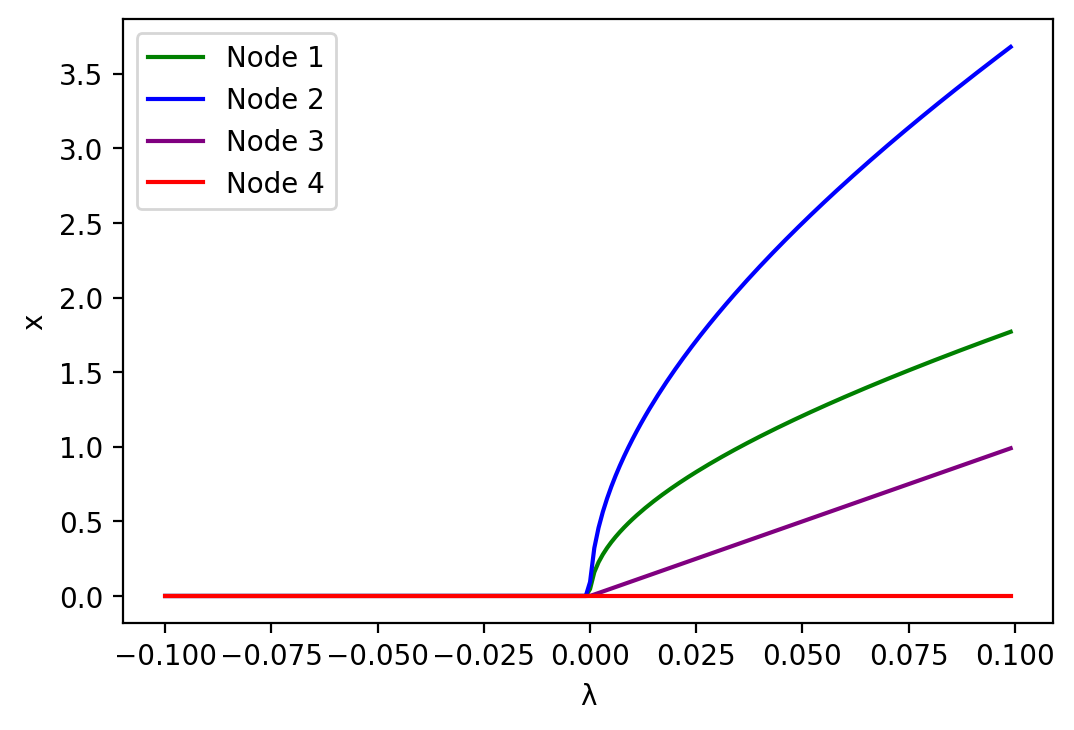}
	\caption{Top: two feedforward chains that differ only slightly. Bottom left: numerical evidence of a bifurcation in the top network,  where one node remains at $0$, two nodes grow linearly and one node appears to grow at a sub-linear rate. Bottom right: numerical evidence of a bifurcation in the bottom network, with one node staying at $0$, one growing linearly and two growing at a sub-linear rate. The same response function was used for both networks.}\label{fig:numerics2}
\end{figure}

Next, we consider the two networks shown at the top of Figure \ref{fig:numerics2}. Their admissible vector fields are of the form
\begin{equation*}\label{eq:numerics2} 
	\gamma^1_f(x, \lambda) =
		\begin{pmatrix}
			f(x_{1},\textcolor{red}{x_{2}},\textcolor{blue}{x_{4}}, \lambda)\\
			f(x_{2},\textcolor{red}{x_{3}},\textcolor{blue}{x_{2}}, \lambda)\\
			f(x_{3},\textcolor{red}{x_{4}},\textcolor{blue}{x_{4}}, \lambda)\\
			f(x_{4},\textcolor{red}{x_{4}},\textcolor{blue}{x_{4}}, \lambda)\\
		\end{pmatrix}  
	\quad\text{and}\quad 
	\gamma^2_f(x, \lambda) = 
		\begin{pmatrix}
			f(x_{1},\textcolor{red}{x_{2}},\textcolor{blue}{x_{1}}, \lambda)\\
			f(x_{2},\textcolor{red}{x_{3}},\textcolor{blue}{x_{4}}, \lambda)\\
			f(x_{3},\textcolor{red}{x_{4}},\textcolor{blue}{x_{4}}, \lambda)\\
			f(x_{4},\textcolor{red}{x_{4}},\textcolor{blue}{x_{4}}, \lambda)\\
		\end{pmatrix},
\end{equation*}
for the top and bottom network, respectively. The only difference between these two networks is the placement of a blue self-loop at node 1 or 2, as opposed to a connection to node 4. The bottom left and right of \Cref{fig:numerics2} show the branches of a bifurcation in the top and bottom network, respectively. They are computed numerically, in the same way as in the previous example, using the response function 
\begin{equation*}
	f(x,y,z, \lambda) = y-2z+ \lambda x -0.1x^2\, .
\end{equation*}
In each of the two examples and for each of 200 values of $\lambda$, the forward Euler-method started at $(0.001,0.002,0.003,-0.004) \in \mathbb{R}^4$ and was performed up to $t=10000$ with time steps of $0.1$.

The bifurcation plots imply a progression of node-asymptotics, from $x_i(\lambda) = 0$ to $x_i(\lambda) \sim \lambda$  to  $x_i(\lambda) \sim \lambda^{1/2}$, starting at node 4 and going through the network nodes from left to right. When going through the nodes in this order,  the asymptotics of a node changes with respect to the previous one, only if it has no self-loops attached. For instance, in both networks node 3 has no self-loops. As a result, we have $x_4(\lambda) = 0$ but $x_3(\lambda) \sim \lambda$. In the top network node 2 does have a self-loop, and so we have $x_2(\lambda) \sim \lambda$. That is, we have the same asymptotics as node 3. In the bottom network node 2 has no self-loops, and as a result the asymptotics of node 2 differs from that of node 3. We therefore find  $x_2(\lambda) \sim \lambda^{1/2}$ instead. In the top network node 1 has no self-loops, and so the asymptotics differs from that of node 2. In the bottom network we do have a self-loop on node 1, and so the asymptotics of nodes 1 and 2 are the same. As a result, we see $x_1(\lambda) \sim \lambda^{1/2}$ in both networks. We will show that such rules of amplification -- progressing through the network and depending on the presence or absence of self-loops -- are a typical occurrence in feedforward networks.

\paragraph{Background.}
A structural feature that arises frequently in the sciences is that of a network of clearly distinguishable units that are connected in a specific configuration to influence each other. Examples arise in engineering (e.g. power grids), biology (e.g. food webs or neural networks), computer science (e.g. deep learning), and many others. In mathematics the study of networks takes place in numerous disciplines -- as well as through interdisciplinary approaches -- such as graph theory, algebra, stochastics or dynamical systems. See for example the survey of the field of network science with historical context and applications in \textcite{Newman.2006} and the references therein. Here, we are interested in dynamical systems with the underlying structure of a network, so called \emph{coupled cell systems}. In particular, we investigate systems of ordinary differential equations where the state variable of one (or multiple) equations is also an argument of the equation for another. We say that the corresponding cell receives an input from other cells. The configuration by which the cells influence each other is often encoded by a (directed) graph.

The field of network dynamical systems has seen enormous activity in recent years and multiple formalisms have been put forward. Most prominently we mention the \emph{groupoid formalism} in \textcite{Golubitsky.2006,Golubitsky.2004} and its equivalent definition in \textcite{Field.2004}. More recently, so-called \emph{open systems} (see e.g. \textcite{Lerman.2018,Schultz.2020}) and \emph{asynchronous networks} (\textcite{Bick.2016,Bick.2017,Bick.2017b}) have been introduced to model more complex applications. Network dynamical systems exhibit interesting phenomena that are highly anomalous in general dynamical systems. Notable examples include synchronization effects and pattern formation, as well as unusual bifurcations and spectral degeneracies. These phenomena are only possible because of the network structure, but the precise mechanisms leading to them often remain unclear. They do, however, resemble observations made in symmetric or equivariant dynamics. It was the introduction of \emph{graph fibrations} from category theory and their implications for dynamical systems (see \textcite{DeVille.2010,DeVille.2015,DeVille.2015b}) that allowed for the development of a theory to interpret networks as algebraic structures that have a direct connection to symmetry, so-called \emph{hidden symmetry} (\textcite{Rink.2013,Rink.2014,Rink.2015,Nijholt.2016,Nijholt.2017,Nijholt.2017c}). This theory applies to the class of networks with \emph{asymmetric inputs} and most results that have been established focus on \emph{homogeneous networks}. We briefly recapitulate the basics in \Cref{sec:pr}.

Additional structure of the network itself can aid the dynamical investigation by providing additional analytical means. A prominent example is that of \emph{feedforward structure}. Broadly speaking a network exhibits feedforward structure if information can only flow in one direction. Information one cell emits cannot become an input into that same cell, not even indirectly, i.e. there are no feedback effects. This rather simple structure has the convenient effect that there is a natural partition of the cells such that the first part receives no inputs from anywhere else in the network, the second receives inputs only from the first, the third receives inputs from the first and the second and so on. This greatly simplifies the investigation of such networks, both in mathematical analysis due to technical simplifications, but also conceptually as it allows to study the network inductively. While feedforward structure and similar weaker notions are abundant in networks -- as a matter of fact, in any network that is not (indirectly) all-to-all coupled, we can find feedforward structure between parts of the network --, they are a prominent feature in deep learning via artificial neural networks, where one type of information is processed by cells within one part and then passed on to the next part until some output is generated. For more information and historical background on this see \textcite{Schmidhuber.2015} and the extensive list of references therein.

The network of feedforward type that was first considered in the network dynamical systems literature is the $3$-cell homogeneous network
	\begin{center}
		\resizebox{.4\linewidth}{!}{
			\begin{tikzpicture}[->,
			>=stealth',
			shorten >=1pt,
			auto,
			node distance=1cm,
			main node/.style={line width=1.5pt, circle, scale = 3, draw, font=\sffamily\tiny, inner sep=2pt}]
			\node[main node] (1) {};
			\node[main node, right of=1] (2) {};
			\node[main node, right of=2] (3) {};
			\path[every node/.style={font=\sffamily\small}, line width =1.5pt]
			(1) edge [loop left] node {} (1)
			(1) edge [] node {} (2)
			(2) edge [] node {} (3)
			;
			\end{tikzpicture}
		}
	\end{center}
(see \textcite{Elmhirst.2006,Golubitsky.2012,Golubitsky.2006}). The first cell is not influenced by any other cell, the second only by the first, and the last only by the second. We refer to this setting as a \emph{feedforward chain}. Note that the self-loop of the first cell, although it seems to contradict the `no-feedback' assumption, is due to a mere convention where we allow cells to influence themselves. It was observed that dynamical systems with the underlying structure of this $3$-cell feedforward network exhibit surprising generic Hopf bifurcations: if a fully synchronous steady state looses stability through a pair of imaginary eigenvalues, the corresponding system can exhibit a Hopf branch in which the first cell remains in the steady state, the amplitude of the second cell grows with rate $\sim|\lambda|^{\frac{1}{2}}$, and the amplitude of the last cell grows with rate $\sim|\lambda|^\frac{1}{6}$, where $\lambda$ is the bifurcation parameter. In particular the growth in the last cell is much faster than expected in `standard' Hopf bifurcations. The effect is also referred to as \emph{amplification} and is forced by the network structure. In \textcite{Rink.2013}, the anomalous Hopf bifurcation result is generalized to feedforward chains of arbitrary length, where the amplification is observed to increase the `further down in the chain' the cell is located. Furthermore, a similar result for steady state bifurcations is proved. Since then, more general classes of feedforward networks, not restricted to chains, were investigated. In \textcite{Nijholt.2017d}, as an example, the authors introduce so-called \emph{ring-feedforward networks} which are feedforward chains where the first cell is replaced by an oriented ring. Most recently, the steady state bifurcation result (as well as other investigations) has been generalized to certain \emph{layered feedforward networks}:
\begin{center}
	\resizebox{.55\linewidth}{!}{
		\begin{tikzpicture}[->,
		>=stealth',
		shorten >=1pt,
		auto,
		node distance=1cm,
		main node/.style={line width=1.5pt, circle, scale = 3, draw, font=\sffamily\tiny, inner sep=2pt}
	]
		\coordinate (c1);
		\node[main node, above of=c1, node distance=.25cm] (11) {};
		\node[main node, below of=c1, node distance=.25cm] (12) {};
		\node[main node, draw=none, above of=c1, node distance=.5cm] (011) {};
		\node[main node, draw=none, below of=c1, node distance=.5cm] (012) {};
		\node[main node, right of=c1] (22) {};
		\node[main node, above of=22, node distance=.5cm] (21) {};
		\node[main node, below of=22, node distance=.5cm] (23) {};
		\node[main node, right of=22,draw=none] (32) {$\dots$};
		\node[main node, above of=32,draw=none, node distance=.5cm] (31) {};
		\node[main node, below of=32,draw=none, node distance=.5cm] (33) {};
		\node[main node, above of=32,draw=none, node distance=.25cm] (031) {};
		\node[main node, below of=32,draw=none, node distance=.25cm] (032) {};
		\node[main node, right of=32,draw=none] (c4) {};
		\node[main node, above of=c4, node distance=.25cm] (41) {};
		\node[main node, below of=c4, node distance=.25cm] (42) {};
		\node[main node, draw=none, above of=c4, node distance=.5cm] (041) {};
		\node[main node, draw=none, below of=c4, node distance=.5cm] (042) {};
		\node[line width=1.5pt, color=gray, draw, rectangle, rounded corners, fit=(011)(012)](l1) {};
		\node[line width=1.5pt, color=gray, draw, rectangle, rounded corners, fit=(21)(23)](l2) {};
		\node[line width=1.5pt, color=gray, draw, rectangle, rounded corners, fit=(041)(042)](l3) {};
		
		\path[every node/.style={font=\sffamily\small}, line width =1.5pt]
		(11) edge [loop left, in=190, out=170, looseness=10] node {} (11)
		(12) edge [loop left, in=190, out=170, looseness=10] node {} (12)
		(11) edge [] node {} (21)
		(12) edge [] node {} (22)
		(12) edge [] node {} (23)
		(21) edge [bend left=5] node {} (31)
		(22) edge [] node {} (32)
		(23) edge [] node {} (032)
		(031) edge [] node {} (41)
		(32) edge [] node {} (42)
		
		(11) edge [loop left, color=red, in=195, out=165, looseness=11] node {} (11)
		(12) edge [loop left, color=red, in=195, out=165, looseness=11] node {} (12)
		(11) edge [color=red, bend left=5] node {} (21)
		(11) edge [color=red, bend left=5] node {} (22)
		(12) edge [color=red, bend left=5] node {} (23)
		(21) edge [color=red] node {} (031)
		(22) edge [color=red] node {} (031)
		(23) edge [color=red, bend left=5] node {} (33)
		(31) edge [color=red] node {} (41)
		(33) edge [color=red] node {} (42)
		
		(11) edge [loop left, color=blue, in=200, out=160, looseness=12] node {} (11)
		(12) edge [loop left, color=blue, in=200, out=160, looseness=12] node {} (12)
		(11) edge [color=blue, bend right=5] node {} (21)
		(11) edge [color=blue, bend right=5] node {} (22)
		(12) edge [color=blue, bend right=5] node {} (23)
		(21) edge [color=blue, bend right=5] node {} (31)
		(22) edge [color=blue] node {} (032)
		(23) edge [color=blue, bend right=5] node {} (33)
		(32) edge [color=blue] node {} (41)
		(032) edge [color=blue] node {} (42)
		;
		\end{tikzpicture}
	}
\end{center}
(see \textcite{Soares.2018}). Therein the cells can be partitioned into layers such that the feedforward structure respects these layers. In particular, if we collapse each layer to one cell, we are left with a feedforward chain. In this paper we investigate the general case that incorporates only the illustrative idea that a feedforward network should not contain any feedback (except for self-loops) (see \Cref{fig:genff}). Our definition includes feedforward chains and layered feedforward networks as special cases. Under the assumption of homogeneity and asymmetry of inputs we prove the aforementioned bifurcation result (see Summary of main results.).

Note that recently research has also extended to networks that do not exhibit a strict feedforward structure. In \textcite{Aguiar.2019b} the authors investigate the effect of feedback on the synchrony patterns of weighted feedforward networks with additive input structure. Furthermore, in \textcite{Gandhi.2020} the feedforward structure of transitive components is exploited to thoroughly investigate $1$- and $2$-parameter steady state bifurcations in fully inhomogeneous networks. A similar investigation is made in \textcite{Aguiar.2019c} for a specific $1$-parameter steady-state bifurcation scenario in homogeneous networks with asymmetric inputs. In this class of networks (or respectively in this bifurcation problem), however, amplification is generically not possible.

\paragraph{Structure of the article.}
This article is structured as follows. In \Cref{sec:pr} we briefly summarize the foundations of homogeneous coupled cell systems with asymmetric inputs. Feedforward structure for this class of networks and some immediate consequences are presented in \Cref{sec:td}. Finally, in \Cref{sec:bi}, all the results are used to compute the generic steady state bifurcations for feedforward networks. As before, we observe the amplification effect for our class of feedforward networks. However, due to the more complicated interaction structures, the picture becomes more complex than in feedforward chains or in layered feedforward networks. In particular, some expected amplifying branches may not exist. These results are illustrated in an example in \Cref{sec:ex}.

%% file: 5cellfeedforwardwithoddbifurcationsextnet.tex
\centering
\begin{tikzpicture}[->,
	>=stealth',
	shorten >=1pt,
	auto,
	node distance=1cm,
	main node/.style={line width=1.5pt, circle, scale = 3, draw, font=\sffamily\tiny, inner sep=1pt}]
	\node[main node] (1) {$1$};
	\node[main node, above left of=1] (2) {$2$};
	\node[main node, below left of=1] (3) {$3$};
	\node[main node, below left of=2] (4) {$4$};
	\node[main node, left of=4] (5) {$5$};
	\path[every node/.style={font=\sffamily\small}, line width =1.5pt]
	(2) edge [color = {red}] node {} (1)
	(5) edge [color = {red}] node {} (2)
	(4) edge [color = {red}] node {} (3)
	(5) edge [bend left = 5, color = {red}] node {} (4)
	(5) edge [in=190, out=170, looseness = 8, color = {red}] node {} (5)
	(3) edge [color = {blue}] node {} (1)
	(4) edge [color = {blue}] node {} (2)
	(5) edge [color = {blue}] node {} (3)
	(5) edge [bend left = -5, color = {blue}] node {} (4)
	(5) edge [in=195, out=165, looseness = 8, color = {blue}] node {} (5)
	(4) edge [color = {grey}] node {} (1)
	(5) edge [bend left =-5, color = {grey}] node {} (2)
	(5) edge [bend left =-5, color = {grey}] node {} (3)
	(5) edge [bend left =-15, color = {grey}] node {} (4)
	(5) edge [in=205, out=155, looseness = 8, color = {grey}] node {} (5)
	(5) edge [in = -100, out = -80, looseness = 1.2, color = {magenta}] node {} (1)
	(5) edge [bend left =5, color = {magenta}] node {} (2)
	(5) edge [bend left =5, color = {magenta}] node {} (3)
	(5) edge [bend left =15, color = {magenta}] node {} (4)
	(5) edge [in=200, out=160, looseness = 8, color = {magenta}] node {} (5)
	;
\end{tikzpicture}

%% file: II_genfeedforward_pr.tex
\section{Preliminaries: Homogeneous Networks with Asymmetric Inputs}
\label{sec:pr}
We consider systems of ordinary differential equations with the underlying structure of a homogeneous coupled cell network, as described in \textcite{Nijholt.2017d, Nijholt.2017c, Nijholt.2016, Rink.2013, Rink.2014, Rink.2015}.
\begin{defi}[homogeneous coupled cell system with asymmetric inputs, Def. 2.1 in \cite{Rink.2014}]
	Let the set of \emph{nodes} (or \emph{cells}) of a network be labeled by $C=\lbrace p_1,\dotsc,p_N \rbrace$ and denote the network interactions in the form of distinct \emph{input maps} $\Sigma = \lbrace \sigma_1, \dotsc , \sigma_n \rbrace$ where each $\sigma_i \colon C \to C$ characterizes one specific \emph{input type}. To each cell we attach the same \emph{internal state space} $V$ which is a finite dimensional real vector space. The \emph{total phase space} is $\bigoplus_{p \in C} V \cong V^N$ with coordinates chosen according to the cells of the network: $x=(x_{p_1},\dotsc,x_{p_N})^T$. The evolution of the state $x_i \in V$ of cell $p_i$ is governed by a function $f\colon V^n \to V$ via its inputs. The network dynamics is governed by the ordinary differential equations
	\begin{equation}
		\label{eq:netvf}
		\dot{x} = \gamma_f(x) = \begin{pmatrix}
			f (x_{\sigma_1(p_1)}, \dotsc, x_{\sigma_n(p_1)}) \\
			f (x_{\sigma_1(p_2)}, \dotsc, x_{\sigma_n(p_2)}) \\
			\vdots \\
			f (x_{\sigma_1(p_N)}, \dotsc, x_{\sigma_n(p_N)})
		\end{pmatrix}. 
	\end{equation}
	These network vector fields are also referred to as \emph{admissible maps} or \emph{admissible vector fields}. Each cell receives precisely one input of each type, hence the term \emph{asymmetric inputs}.
	\hspace*{\fill}$\triangle$
\end{defi}
\noindent
We make one additional assumption on the set of input maps $\Sigma$. We want it to include the identity map $\sigma_1 = \Id \colon C \to C$, which is natural, as it only means that the evolution of each cell's state depends on its own state. 
In order to investigate the inputs that a specific cell $p \in C$ receives, we define $\Sigma(p)=\{\sigma(p) \mid \sigma\in\Sigma\}$ and $\SSigma(p)=\Sigma(p)\setminus\{p\}$ to denote the sets of cells that $p$ receives an input from with and without self-loops respectively. Since $\Id\in\Sigma$, we have $\SSigma(p)\subsetneq\Sigma(p)$ for all $p \in C$.

In the remainder of this section, we recall some useful facts and definitions of homogeneous coupled cell systems. We start with the following characterization of linear admissible maps. Let us denote the algebra of linear maps on the internal phase space of the coupled cell system by $\gl{V}$.
\begin{prop}
	\label{prop:linadm}
	For each $\sigma \in \Sigma$ define a linear map $B_\sigma \colon V^N \to V^N$ by $(B_{\sigma}(x))_p = x_{\sigma(p)}$. Then any linear admissible map $L$ is of the form
	\begin{equation}
		\label{eq:linadm}
		(Lx)_p = \sum_{\sigma \in \Sigma} b_\sigma (B_\sigma (x))_p
	\end{equation}
	where $b_\sigma \in \gl{V}$ are linear maps on $V$ independent of $p$.
\end{prop}
\begin{proof}
	A linear admissible map $L$ is uniquely defined by linear internal dynamics given by a map $\mathfrak{l} \colon V^n \to V$, i.e. $L = \gamma_{\mathfrak{l}}$. As $\mathfrak{l}$ is linear and its arguments are labeled by the input maps $\sigma\in\Sigma$, we find $b_\sigma \in \gl{V}$ such that
	\[ \mathfrak{l}(Y) = \sum_{\sigma \in \Sigma} b_\sigma Y_\sigma, \]
	where $Y=(Y_\sigma)_{\sigma\in\Sigma} \in V^n$. Then for $x\in V^N$ the $p$-th entry of $Lx$ depends on the entries of cells that $p$ receives an arrow from, i.e. there is $\sigma \in \Sigma$ such that $\sigma(p)=q$. We obtain
	\[ (Lx)_p = \mathfrak{l}(x_{\sigma_1(p)}, \dotsc, x_{\sigma_n(p)}) = \sum_{\sigma \in \Sigma} b_\sigma x_{\sigma(p)} = \sum_{\sigma \in \Sigma} b_\sigma (B_\sigma (x))_p. \]
\end{proof}
\noindent
In particular, we will use the following straightforward special case for one-dimensional internal dynamics:
\begin{cor}
	\label{cor:linadm}
	If $V=\RR$, any linear admissible map $L$ is a real linear combination of the $B_\sigma$, i.e.
	\[ L = \sum_{\sigma \in \Sigma} b_\sigma B_\sigma \]
	for some $b_\sigma \in \RR$.
\end{cor}

Finally we state two definitions. The first makes the concept of a path precise in our setting. The second defines subnetworks of a network as subsets of the set of cells that are not influenced by any cell outside of the subset.
\begin{defi}
	For two cells $p, q \in C$ let $\omega = \{p_1, \dotsc, p_k\}$ be a \emph{(loop-free) path from $p$ to $q$} if $p_1 = p, p_k = q$, the $p_i$ are pairwise non-equal, and there exist $\sigma_1, \dotsc, \sigma_{k-1} \in \Sigma$ such that \mbox{$\sigma_1(p_k) = p_{k-1}, \dotsc, \sigma_{k-1}(p_2) = p_1$}. Define \mbox{$\Omega_{p,q} = \{ \omega \text{ path from } p \text{ to } q \}$}.
	\hspace*{\fill}$\triangle$
\end{defi}
\begin{remk}
	Note that the convention of denoting a path as a set implies that it does not contain any self-loops, i.e. it is loop-free by definition.
	\hspace*{\fill}$\triangle$
\end{remk}
\begin{defi}
	A \emph{cycle of length $k$} is a path $\omega = \{p_1, \dotsc, p_k\}$ with the additional property that there exists $\sigma\in\Sigma$ such that $\sigma(p_1)=p_k$. A cycle of length $1$ is also called a \emph{self-loop}.
\end{defi}
\begin{defi}
	A subset $B \subset C$ of cells of a homogeneous coupled cell system is called a \emph{subnetwork} if there are no arrows in the network starting outside of $B$, that target a cell inside of $B$. In other words $\sigma(b) \in B$ for all $b \in B$ and $\sigma \in \Sigma$.
	\hspace*{\fill}$\triangle$
\end{defi}

%% file: II_genfeedforward_td.tex
\section{Feedforward networks}
\label{sec:td}
In this section, we provide a general notion of feedforward networks. We then observe, that it can equivalently be defined in terms of a partial order on the cells. Finally, we explore some direct consequences of these definitions for later use.

\begin{defi}
	A homogeneous coupled cell network is called a \emph{feedforward network} if it has no cycles of length $2$ or more (consisting of arrows not necessarily of the same type). Put differently, a network is a feedforward network if the only cycles are self-loops.
	\hspace*{\fill}$\triangle$
\end{defi}
\begin{remk}
	The name \emph{feedforward network} is natural in the sense that all the \emph{feedback} a cell can receive is via a self-loop.
	\hspace*{\fill}$\triangle$
\end{remk}
Next, we define the preorder
\[ p \nole q \quad \iff \quad \text{there is a path from } q \text{ to } p. \]
As every cell is coupled to itself via the internal dynamics the preorder is obviously reflexive. On the other hand, if there is a path from cell $q$ to cell $p$ and a path from cell $r$ to cell $q$ the joint path goes from $r$ to $p$, which makes the preorder transitive as well.

We encode the situation that $p \nole q$ but $p\ne q$ by $p \nol q$. In the definition we do not exclude the possibility that there exists a path from $q$ to $p$ and one from $p$ to $q$ for two different cells $p$ and $q$. In that case we have $p \nole q$ and $q \nole p$ even though $p \ne q$. Thus, in general the preorder is no partial order.

\begin{prop}
	\label{prop:equivalence}
	A homogeneous coupled cell network is a feedforward network if and only if $\nole$ is a partial order.
\end{prop}
\begin{proof}
	Assume the network is not a feedforward network. Then there is a cycle $\left\{p_1, \dotsc, p_k\right\}$ with $k \ge 2$ in the network, meaning there are input maps $\sigma_1, \dotsc , \sigma_k$ such that $\sigma_1 (p_1) = p_k$ and \mbox{$\sigma_i(p_i) = p_{i-1}$} for $2 \le i \le k$. 
	In particular, there is a path from $p_1$ to $p_k$, which implies $p_k \nole p_1$, and a path from $p_k$ to $p_1$, which implies $p_1 \nole p_k$. As $p_1 \ne p_k$, this shows the preorder is no partial order.
	
	Now assume that the preorder is not a partial order. Then there are cells $p \ne q$ with $p \nole q$ and $q \nole p$. This implies the existence of 
	a path $\{p, p_1, \dotsc, p_{k-2}, q\}$ of length $k \ge 2$ from $q$ to $p$ as well as a path $\{q, p'_1, \dotsc, p'_{l-2}, p\}$ of length $l\ge2$ from $p$ to $q$. The concatenated path $\{p, p_1, \dotsc, p_{k-2}, q, p'_1, \dotsc, p'_{l-2}\}$ is a cycle of length $k+l-2 \ge 2$. Thus, the network is not a feedforward network.
\end{proof}

In the remainder of this section, we collect some consequences of the definition of feedforward networks. Note that, since the network contains only finitely many cells, there are well-defined maximal elements with respect to the partial order $\nole$. By definition, these are cells that do not receive any inputs from other cells.
\begin{lemma}
	\label{lem:maxfix}
	A cell $p \in C$ is maximal with respect to $\nole$ if and only if all its inputs are from itself.
\end{lemma}
\begin{proof}
	This follows almost directly from the definition of the partial order $\nole$. We prove the statement by contraposition. Assume, $\sigma(p) \ne p$ for some $\sigma \in \Sigma$. Then there is an arrow from $\sigma(p)$ to $p$. This yields $\sigma(p) \nog p$ so that $p$ is not maximal. On the other hand, if $q \nog p$ then there is a path $\{p, p_1, \dotsc, p_k, q\}$ from $q$ to $p$. In particular, there is an input map $\sigma\in\Sigma$ such that $\sigma(p)=p_1$, i.e., $p$ receives an input from $p_1\ne p$.
\end{proof}

The following results refine the notions of paths and subnetworks for feedforward networks.
\begin{cor}
	\label{cor:maxpath}
	Stating the definition of a maximal cell in terms of arrows in the network immediately proves that every cell $p$ is either maximal itself or there is a path from a maximal cell $\overline{p}$ to $p$. Furthermore, it is obvious that $q \nog p$ implies that $q$ is on some path from a maximal cell to $p$, i.e., there is a maximal cell $\overline{p}$ such that $q \in \omega \in \Omega_{\overline{p}, p}$ (the set of paths from $\overline{p}$ to $p$).
\end{cor}
\begin{lemma}
	A subnetwork of a feedforward network is again a feedforward network.
\end{lemma}
\begin{proof}
	Let $C$ be the set of cells of a feedforward network. In particular, there are no cycles of length two or greater of cells in $C$. Hence, there are also no cycles of length two or greater of cells in any subset $B\subset C$.
\end{proof}
\begin{lemma}
	\label{cor:subnetmax}
	A subset of cells $B \subset C$ of a feedforward network defines a subnetwork if and only if $p \in B$ implies $q \in B$ for all $q \nog p$. This, in turn, yields that every subnetwork contains at least one maximal cell.
\end{lemma}
\begin{proof}
	The result follows from the definitions of a subnetwork and that of the partial order $\nole$.
\end{proof}
\begin{cor}
	Let $\emptyset \ne B\subset C$ be a subnetwork and $p \in C\setminus B$ such that $q \in B$ for all $q \nog p$. Then also $B \cup \{p\}$ is a subnetwork.
\end{cor}
\begin{proof}
	The result follows immediately from the fact that $\sigma(B)\subset B$ for all $\sigma\in\Sigma$, since $B$ is a subnetwork, and $\sigma(p)\noge p$ for all $\sigma\in\Sigma$.
\end{proof}
\begin{lemma}
	Let $\emptyset\ne B \subsetneq C$ be a non-trivial subnetwork of a feedforward network and assume $B$ contains all maximal cells. Then there exists at least one $p \in C\setminus B$ such that $q \in B$ for all $q \nog p$. That is, there are cells that `surround' the subnetwork.
\end{lemma}
\begin{proof}
	Since $B$ contains all maximal cells, there is a path from a cell in $B$ to any $p \in C\setminus B$, i.e., there is a cell $\overline{p} \in B$ such that $\overline{p} \nog p$. Assume there is no $p \in C\setminus B$ such that $q\in B$ for all $q \nog p$. Then for all $p \notin B$ there must be $p' \in C \setminus B$ such that $p' \nog p$. As $C \setminus B$ is finite and contains no maximal cells, this implies that there are $p, p' \in C \setminus B$ with $p' \nog p$ and $p \nog p'$. This contradicts the assumption that $C$ is a feedforward network.
\end{proof}

The definition of a feedforward network yields additional structure in linear admissible maps for corresponding coupled cell systems. We choose a labeling of the nodes $C = \{p_1, \dotsc, p_N\}$ such that it holds that
\begin{equation}
	\label{eq:ordering}
	p_i \nole p_j \quad \implies \quad i\le j.
\end{equation}
This ordering is not unique as some elements may not be related by the partial order. Recall that we denote the algebra of linear maps on the internal phase space $V$ by $\gl{V}$.
\begin{lemma}
	\label{lem:triangular}
	Choosing the ordering of cells according to \eqref{eq:ordering} for a feedforward network yields that any linear admissible map can be represented by an upper triangular matrix with entries in $\gl{V}$. In particular, if $V=\RR$ we identify $\gl{V}\cong\RR$ to see that the linear admissible maps can be represented by real upper triangular matrices.
\end{lemma}
\begin{proof}
	The lemma follows directly from \Cref{prop:linadm}. Fix an input map $\sigma\in\Sigma$. For an input $\sigma(p)$ of cell $p$ it holds that $\sigma(p) \noge p$ and therefore its index is greater than or equal to that of $p$. This shows that the linear map $B_\sigma$ is upper triangular with non-zero entries $\idv\in\gl{V}$. Then \eqref{eq:linadm} implies that any linear admissible map $L$ is upper triangular with entries in $\gl{V}$ as well.
\end{proof} 
\begin{remk}
	We could have used upper triangularity of linear admissible maps as a definition for feedforward networks as well. Assume the cells $\lbrace p_1, \dotsc, p_N \rbrace$ are labeled such that all linear admissible maps are upper triangular. Suppose $p_i \nole p_j$ and $p_j \nole p_i$. Then there exists a path $\{p_j, p_{l_1}, \dotsc, p_{l_k}, p_i\}$ with $0\le k \le N$ and $1 \le l_r \le N$ from $p_j$ to $p_i$. In particular, there exist input maps $\sigma_1, \dotsc, \sigma_{k+1}\in\Sigma$ such that $\sigma_1(p_i)=p_{l_k},\dotsc, \sigma_{k+1}(p_{l_1})=p_j$.
	These maps describe the input structure that is reflected by the linear admissible maps. Hence, there exists such a map with non-zero values in the $(i,l_k)$-th, in the $(l_k,l_{k-1})$-th, and so on until the $(l_1,j)$-th entry. Due to the upper triangular structure, this implies $i \le l_k, l_k \le l_{k-1}, \dotsc, l_1 \le j$ and, in particular, $i \le j$. Applying the same argument to $p_j \nole p_i$ implies $j \le i$. As a result, we obtain $i=j$ and therefore $p_i = p_j$. Hence $\nole$ is a partial order. Summarizing, the upper triangular structure of linear admissible maps is equivalent to $\nole$ being a partial order on $C$.
	\hspace*{\fill}$\triangle$
\end{remk}

\Cref{lem:triangular} has immediate consequences for generic $1$-parameter Hopf bifurcations in feedforward networks. These can only occur, when the internal dynamics is at least $2$-dimensional. A thorough investigation of steady state bifurcations in feedforward networks can be found in \Cref{sec:bi}.
\begin{theorem}
	\label{thm:ffhopf}
	In a $1$-parameter bifurcation in a feedforward network with one-dimensional internal dynamics there cannot be a pair of conjugate imaginary eigenvalues at the synchronous bifurcation point. On the other hand, if the internal dynamics is at least $2$-dimensional, a $1$-parameter bifurcation in which a pair of complex eigenvalues crosses the imaginary axis is possible.
\end{theorem}
\begin{proof}
	In order for a bifurcation to occur, a $1$-parameter family of linear admissible maps has to have an eigenvalue/a pair of complex conjugate eigenvalues that crosses/cross the imaginary axis at the bifurcation point. In the case $V=\RR$ all linear admissible maps are real upper triangular matrices. Their eigenvalues are the diagonal elements which are real. Hence, only real eigenvalues can cross the imaginary axis. On the other hand, when the internal dynamics is in $V\cong \RR^d$ with $d\ge 2$, linear admissible maps are upper triangular with entries in $\gl{V}$. Thus the eigenvalues of a linear admissible map are the union of the eigenvalues of all diagonal elements which are arbitrary elements in $\gl{V}$ (some of which might be related). In particular, there are possible diagonal elements with complex eigenvalues.
\end{proof}
\begin{remk}
	In particular, the emergence (or collapse) of periodic solutions in a bifurcation -- as in classical and non-classical Hopf bifurcations -- requires a pair of complex conjugate eigenvalues to cross the imaginary axis. The previous theorem shows that this can only occur in feedforward networks, if the internal dynamics is at least $2$-dimensional.
	\hspace*{\fill}$\triangle$
\end{remk}

In the case $V = \RR$ the linear admissible maps are real upper triangular matrices (compare to \Cref{cor:linadm}). These have their eigenvalues on the diagonal. We introduce the following definition that turns out to be useful for determining diagonal entries of linear admissible maps.
\begin{defi}
	Given a homogeneous coupled cell network (not necessarily of feedforward type), we define an equivalence relation $\multimapboth$ on the nodes as follows: $p \multimapboth q$ if
	\begin{equation}
		\label{eq:deflooptype}
		\{\sigma \in \Sigma \mid \sigma(p) = p\} = \{\sigma \in \Sigma \mid \sigma(q) = q\}.
	\end{equation}
	We denote the sets involved by $\LL_p = \{\sigma \in \Sigma \mid \sigma(p) = p\}$. If $p \multimapboth q$ then we say that $p$ and $q$ have the same \textit{loop-type}. In a network, two nodes have the same loop-type if and only if they have the same self-loops (of the same type given by $\LL_p$).
	\hspace*{\fill}$\triangle$
\end{defi}
\begin{cor}
	\label{cor:maxloopeq}
	From \Cref{lem:maxfix} we obtain $\LL_p = \Sigma$ for all maximal $p \in C$ in a feedforward network. In particular, for $p$ maximal, $p \loopeq q$, if and only if $q$ is maximal.
\end{cor}
\begin{remk}
	The loop-type of a node $p$ can easily be read off from the admissible vector fields of a network; it is given by those entries of the response function through which the variable $x_p$ depends on $x_p$.
	\hspace*{\fill}$\triangle$
\end{remk}
The following result shows the importance of the loop-type relation: there is a one-to-one correspondence between loop-types and eigenvalues of a linear admissible map.
\begin{theorem}
	\label{thm:genev}
	In a feedforward network, the number of different loop-types of the nodes (that is, the number of equivalence classes under $\multimapboth$) equals the maximal number of different eigenvalues a linear admissible map can have.
\end{theorem}
\begin{proof}
	As a linear admissible map for a feedforward network is upper triangular (\Cref{lem:triangular}), the maximal number of eigenvalues is just the maximal number of distinct values the diagonal entries can attain. We define the standard basis of the total phase space using the Kronecker delta $\delta_{p,q}$ -- which equals $1$ if $p=q$ and $0$ otherwise -- as $\{Y^q\}_{q\in C}$ given by $(Y^q_p)_{p \in C} = (\delta_{p,q})_{p\in C} \in \RR^N$. Note that this basis respects the labeling of cells chosen as in \eqref{eq:ordering} so that the linear admissible maps with respect to this basis are upper triangular. Recall furthermore from \Cref{cor:linadm} that the algebra of linear admissible maps is spanned by the adjacency matrices $B_{\sigma}$ for $\sigma \in \Sigma$, defined by $(B_{\sigma}(x))_p = x_{\sigma(p)}$ for $p \in C$. That is, any linear admissible map is of the form
	\[\sum_{\sigma \in \Sigma}b_{\sigma}B_{\sigma}, \]
	with $b_\sigma\in\RR$ for all $\sigma\in\Sigma$. For a node $p \in C$, the $(p,p)$-entry of this matrix is given by 
	\begin{equation*}
		\left[\sum_{\sigma \in \Sigma}b_{\sigma}B_{\sigma}(Y^p)\right]_p = \sum_{\sigma \in \Sigma}b_{\sigma}Y^p_{\sigma(p)} = \sum_{\sigma \in \Sigma}b_{\sigma}\delta_{p, \sigma(p)} = \sum_{\sigma(p) = p} b_{\sigma}\, .
	\end{equation*}
	Hence, for two nodes $p, q \in C$ the $(p,p)$-entry and the $(q,q)$-entry are always the same, if and only if $p \multimapboth q$. This shows that the number of different eigenvalues is at most the number of loop-types. If $p$ and $q$ do not have the same loop-type, then for a dense open set of values $(b_{\sigma})_{\sigma \in \Sigma}$, the $p$-th and $q$-th diagonal entries of $\sum_{\sigma \in \Sigma}b_{\sigma}B_{\sigma}$ are distinct. Intersecting these sets for all pairs of nodes with a different loop-type, we find a dense open set of values $(b_{\sigma})_{\sigma \in \Sigma}$ for which $\sum_{\sigma \in \Sigma}b_{\sigma}B_{\sigma}$ has as many different eigenvalues as there are loop-types. This proves the theorem.
\end{proof}

\begin{remk}
	\label{rem:layered}
	As mentioned in the introduction, an often considered generalization of feedforward chains is that of so-called \emph{layered feedforward networks}. In such a network the cells are partitioned in layers and the feedforward structure is only with respect to these layers, i.e., if we collapse each layer to a single node we obtain a feedforward chain. Most notably we would like to mention \textcite{Soares.2018}. Therein, layered feedforward networks, their quotients and lifts as well as the lifting bifurcation problem are investigated thoroughly. It is shown that these networks exhibit similar steady state bifurcations as the ones presented in \Cref{sec:bi}. As a matter of fact, completing the set of input maps for the networks considered in \cite{Soares.2018}, we obtain networks that are included in our framework of feedforward networks. However, no self-loops (except for the internal dynamics governed by $\Id\in\Sigma$ and for the maximal cells) are possible. As a result, the branching patterns in the more general case considered in \Cref{sec:bi} are more complex than in \cite{Soares.2018}. Note that in our more general definition of feedforward networks, we may also group nodes together in layers, where we allow for self-loops and arrows that skip layers. In fact, it can be shown that we may give a definition in terms of layers that is equivalent to the definitions in \Cref{prop:equivalence}. However, we decided not to include the precise definition here, as we deemed the `graphical' and `order theoretic' definitions to be more natural.
	\hspace*{\fill}$\triangle$
\end{remk}

%% file: II_genfeedforward_bi.tex
\section{Amplified Steady State Bifurcations in Feedforward Networks}
\label{sec:bi}
In this section, we classify generic steady state bifurcations in one-parameter families of coupled cell network vector fields as in \eqref{eq:netvf} where the underlying structure is that of a feedforward network. It turns out, that the most useful definition for the computations in this section is that the set of cells $C=\lbrace p_1,\dotsc,p_N \rbrace$ is partially ordered with respect to $\nole$. Assuming the labeling of nodes is according to \eqref{eq:ordering}, the parameter-dependent dynamics on the total phase space $V^N \cong \bigoplus_{p \in C} V$ is governed by
\begin{equation}
\label{eq:netvfparam}
\dot{x} = \gamma_f(x,\lambda) = \begin{pmatrix}
f (x_{\sigma_1(p_1)}, \dotsc, x_{\sigma_n(p_1)}, \lambda) \\
f (x_{\sigma_1(p_2)}, \dotsc, x_{\sigma_n(p_2)}, \lambda) \\
\vdots \\
f (x_{\sigma_1(p_N)}, \dotsc, x_{\sigma_n(p_N)}, \lambda)
\end{pmatrix},
\end{equation}
where $\lambda \in \RR$ and $\gamma_f(x, \lambda)_p$ depends only on $q \in C$ with $q \noge p$.

We aim at investigating generic bifurcations from a fully synchronous steady state. Without loss of generality, we may assume this to be the origin and the bifurcation to occur for $\lambda=0$. Hence, we assume 
\[ \gamma_f(0,0) = 0, \]
which implies $f(0,0) = 0$.

Due to the implicit function theorem, a bifurcation of steady states can only occur if the linearization $D_x \gamma_f (0,0)$ is non-invertible. As the inputs of $f$ are labeled by the input maps $\sigma \in \Sigma$, we may define
\[ a_\sigma = \partial_\sigma f(0,0), \]
which is an arbitrary linear map on $V$, i.e. $a_\sigma \in \gl{V}$. Furthermore, recall the definition 
\[ \LL_p = \{ \sigma \in \Sigma \mid \sigma (p) = p \} \]
as in \Cref{eq:deflooptype}. Then $p$ and $q$ are of the same loop-type, $p \loopeq q$, if and only if $\LL_p = \LL_q$. We compute the linearization to be the block-triangular matrix
\begin{equation}
	\label{eq:lin}
	D_x \gamma_f (0,0) = \begin{pmatrix}
		\sum_{\sigma \in \LL_{p_1}} a_\sigma	& \bullet								& \dots		& \bullet	\\
		0										& \sum_{\sigma \in \LL_{p_2}} a_\sigma	& \ddots	& \vdots 	\\
		\vdots									& \ddots								& \ddots	& \bullet	\\
		0										& \dots									& 0			& \sum_{\sigma \in \LL_{p_N}} a_\sigma
	\end{pmatrix}.
\end{equation}
Then, the linearization $D_x \gamma_f(0,0)$ is non-invertible, if and only if there is a node $\overline{p} \in C$ such that \mbox{$0 \in \spec{\sum_{\sigma \in \LL_{\overline p}} a_\sigma}$}. Note that $\sum_{\sigma \in \LL_p} a_\sigma = \sum_{\sigma \in \LL_{q}} a_\sigma$, if $p \multimapboth q$. Thus, $0 \in \spec{\sum_{\sigma \in \LL_p} a_\sigma}$ for all $p \multimapboth \overline{p}$. We call these nodes \emph{critical}. Furthermore, as $a_\sigma \in \gl{V}$ arbitrary, generically $\sum_{\sigma \in \LL_{p}} a_\sigma$ does not have an eigenvalue $0$ if $p$ is not of the same loop-type as $\overline{p}$ (compare to \Cref{thm:genev}). Hence, we say and assume the loop-type of $\overline{p}$ is \emph{critical} and all other loop-types are \emph{non-critical}.

The steady state bifurcation problem is to find solutions to
\begin{equation}
\label{eq:Barb}
\gamma_f (x, \lambda) = 0
\end{equation}
locally around $(0,0) \in \bigoplus_{p \in C} V \times \RR$ for a generic $f$ satisfying the basic bifurcation assumptions (B), i.e., 
\begin{enumerate}[label=(B.\roman*)]
	\item $f(0,0) = 0$; \label{enum:B1}
	\item there exists $\overline{p} \in C$ such that $0 \in \spec{\sum_{\sigma \in \LL_{\overline{p}}} a_\sigma}$. \label{enum:B2}
\end{enumerate}

An immediate observation following the definition of criticality is 
\begin{lemma}
	\label{lem:maxcrit}
	Maximal cells are either all critical or all non-critical.
\end{lemma}
\begin{proof}
	This follows from the fact that every maximal cell $p$ receives all inputs from itself, i.e. $\LL_{p} = \Sigma$ -- compare to \Cref{lem:maxfix} and \Cref{cor:maxloopeq}.
\end{proof}
\noindent
Hence, when the maximal cells are critical, all non-maximal cells may be assumed not to be and vice versa when a non-maximal cell is critical all maximal cells may be assumed not to be.

We explicitly compute the generic steady state bifurcation behavior in individual cells in the case of one-dimensional internal dynamics, that is $V = \RR$. The results and their proofs are rather notation heavy. Furthermore, they consist of overlapping inductive definitions and statements. However, the branching solutions can be summarized informally as follows:
\begin{enumerate}[label=(\roman*)]
	\item The steady state solutions in maximal cells grow asymptotically as $\sim \lambda$, if these maximal cells are non-critical, and as $\sim \pm \sqrt{|\lambda|}$, if they are critical.
	\item The steady state solutions of non-critical cells are, to leading order, linear in their inputs, i.e., $x_p \sim \max_{\sigma\in\Sigma} |x_{\sigma(p)}|$. Hence, to leading order they have the same asymptotics as the leading order of their inputs.
	\item The steady state solutions of critical cells grow, to leading order, as the square root of their lowest order inputs, i.e., $x_p \sim \pm \sqrt{\max_{\sigma\in\Sigma}|x_{\sigma(p)}|}$. We refer to this phenomenon as \emph{amplification}.
\end{enumerate}
These results follow from an inductive investigation of the bifurcations in individual cells starting from the maximal cells and following the order $\nole$. Additionally, we have to carefully distinguish the cases when a bifurcation occurs for positive or negative values of $\lambda$. Furthermore, taking square roots of inputs is only possible if the signs of inputs are suitable. This results in restrictions on system parameters -- i.e. Taylor coefficients or partial derivatives of the governing function $f$.

The key ingredient in the computations in this section is the fact that for an arbitrary cell $p \in C$ the function $\gamma_f(x, \lambda)_p = f (x_{\sigma_1(p)}, \dotsc, x_{\sigma_n(p)}, \lambda)$ depends only on those $q \in C$ with $q \noge p$. Together with the bifurcation assumption \ref{enum:B1} this allows us to Taylor expand the governing function to obtain the following expanded bifurcation equation
\begin{equation}
	\label{eq:taylor}
	\begin{split}
	0 = f (x_{\sigma_1(p)}, \dotsc, x_{\sigma_n(p)}, \lambda)	&= \sum_{\sigma, \tau \in \LL_{p}} f_{\sigma\tau} x_p^2 + \sum_{\sigma \in \LL_{p}} a_\sigma x_p + \sum_{\sigma \in \LL_{p}} f_{\sigma \lambda} \lambda x_p \\
																&\phantom{=} + \sum_{p \nol q = \tau(p) \colon \atop \tau \notin \LL_p} a_\tau x_q + \ell \lambda + \sum_{p \nol q = \tau(p) \colon \atop \tau \notin \LL_p} f_{\tau\lambda} \lambda x_q \\
																&\phantom{=} + 2 \sum_{\substack{\sigma \in \LL_p \\ \tau \notin \LL_p \colon \\ p \nol q = \tau(p)}} f_{\sigma\tau} x_p x_q + \sum_{\substack{\sigma, \tau \notin \LL_p \\ p \nol q,s = \sigma(p),\tau(p)}} f_{\sigma\tau} x_q x_s + f_{\lambda\lambda} \lambda^2 \\
																&\phantom{=} + \OO \left( \| ((x_q \mid q \noge p), \lambda) \|^3 \right).
	\end{split}
\end{equation}
The constants $f_{\sigma\tau}, f_{\tau\lambda}, \ell$, and $f_{\lambda\lambda}$ are defined to be partial derivatives of $f$ in directions labeled by the input functions $\sigma \in \Sigma$, similar to the $a_\sigma$ before, or by the parameter $\lambda$:
\begin{align*}
	f_{\sigma\tau}		&= \frac{1}{2} \partial_{\sigma\tau} f(0,0) \quad \text{(hence } f_{\sigma\tau} = f_{\tau\sigma} \text{)}, \\
	f_{\sigma\lambda}	&= \partial_{\tau\lambda} f(0,0) = \partial_{\lambda\tau} f(0,0),\\
	\ell				&= \partial_\lambda f(0,0), \\
	f_{\lambda\lambda}	&= \frac{1}{2} \partial_{\lambda\lambda} f(0,0).
\end{align*}
Note that in the case $V=\RR$ all these constants are real numbers and especially \ref{enum:B2} becomes $\sum_{\sigma \in \LL_{\overline{p}}} a_\sigma = 0$ for a critical cell $\overline{p}$. We will use \eqref{eq:taylor} to inductively solve the bifurcation equation \eqref{eq:Barb} with respect to $\nole$. It turns out that this results in significantly different solutions depending on the two possible cases presented in \Cref{lem:maxcrit} -- either the critical cells are maximal or not. Even though the computations follow the general idea outlined above in both cases, it is convenient to separate the investigations.

\begin{remk}
	Note that, due to the partial order, the bifurcation problem (B) can be `restricted' to subnetworks with the `same' branching pattern restricted to the subnetwork. By definition, the dynamics of a subnetwork $B \subset C$ is governed by the ordinary differential equations
	\[ \dot x_p = f(x_{\sigma_1(p)}, \dotsc, x_{\sigma_n(p)}, \lambda), \]
	where $f$ is as in \eqref{eq:netvfparam} and $p \in B$. In particular, if the network undergoes a steady state bifurcation according to the bifurcation assumption (B), then the bifurcation assumption (B) (and the Taylor expansion \eqref{eq:taylor}) holds for the subnetwork as well. In particular, if $B$ contains a critical cell, the subnetwork undergoes a steady state bifurcation with the branching pattern of the full network restricted to cells in $B$. The linearization of the restricted vector field arises from \eqref{eq:lin}, by removing rows corresponding to cells in $C\setminus B$. Hence, a cell $p \in B$ is critical for the reduced bifurcation problem, if and only if it is critical for the full bifurcation problem.
	\hspace*{\fill}$\triangle$
\end{remk}


\subsection{The critical cells are maximal}
\label{subsec:maxcrit}
We start with the simpler case by assuming the maximal cells to be critical. From \Cref{lem:maxcrit} we know that this is equivalent to $\sum_{\sigma\in\Sigma} a_\sigma = 0$ and to all non-maximal cells not being critical generically. This greatly simplifies the computations as there are fewer cases to take care of. Furthermore, the bifurcation condition \ref{enum:B2} is equivalent to
\begin{equation}
	\label{eq:maxcrit}
	\sum_{\sigma \in \LL_{p}}a_\sigma = - \sum_{\tau \notin \LL_{p}} a_\tau
\end{equation}
for every cell $p\in C$, which is a useful identity. We start by computing the bifurcation behavior of \eqref{eq:taylor} for maximal cells, where we detect a saddle node bifurcation. Then we proceed inductively with respect to $\nole$ to all non-maximal cells, whose state variables mimic the bifurcation behavior of their inputs.
\begin{lemma}
	\label{lem:c1-st1}
	Let $p \in C$ be maximal. The state variable $x_p$ generically bifurcates as in one of the following two cases:
	\begin{enumerate}[label=(\roman*)]
		\item \[ x_p(\lambda) = D_p^\pm \cdot \sqrt{\lambda} + \OO(|\lambda|) \]
		for $\lambda > 0$ small, if  
		\[ \frac{\ell}{\sum_{\sigma, \tau \in \Sigma} f_{\sigma\tau}} <0. \]
		Here
		\[ D_p^\pm = \pm \sqrt{-\frac{\ell}{\sum_{\sigma, \tau \in \Sigma} f_{\sigma\tau}}} \ne 0. \]
		\item \[ x_p(\lambda) = D_p^\pm \cdot \sqrt{-\lambda} + \OO(|\lambda|) \]
		for $\lambda < 0$ small, if 
		\[ \frac{\ell}{\sum_{\sigma, \tau \in \Sigma} f_{\sigma\tau}} >0. \]
		Here
		\[ D_p^\pm = \pm \sqrt{\frac{\ell}{\sum_{\sigma, \tau \in \Sigma} f_{\sigma\tau}}} \ne 0. \]
	\end{enumerate}
\end{lemma}
\begin{proof}
	As a maximal cell $p$ only receives inputs from itself (see \Cref{lem:maxfix}), the equations to be solved \eqref{eq:taylor} only depend on the state variable of that specific cell. Hence, the bifurcation behavior is exactly the same for all maximal cells. As $\sum_{\sigma \in \Sigma} a_\sigma = 0$, \eqref{eq:taylor} becomes
	\begin{equation}
	\label{eq:c1-st1}
	0 = \sum_{\sigma, \tau \in \Sigma} f_{\sigma\tau} x_p^2 + \ell \lambda + \OO(|\lambda|^2 + |x_p|^3 + |x_p||\lambda|).
	\end{equation}
	In both cases we employ the standard method to detect saddle node bifurcations (see for example \textcite{Murdock.2003}). Assume $\ell / \sum_{\sigma, \tau \in \Sigma} f_{\sigma\tau} <0$. We introduce a new variable $x_p=\mu y$ where $\mu = \sqrt{\lambda}$ for small $\lambda>0$, hence $\mu>0$. \Cref{eq:c1-st1} transforms into
	\[ 0 = \sum_{\sigma, \tau \in \Sigma} f_{\sigma\tau} \mu^2y^2 + \ell \mu^2 + \OO(|\mu|^4 + |\mu|^3|y|). \]
	As $\mu >0$, we may divide by $\mu^2$ and obtain
	\[ 0 = \sum_{\sigma, \tau \in \Sigma} f_{\sigma\tau} y^2 + \ell + \OO(|\mu|^2 + |\mu||y|) = g(y, \mu). \]
	For $\mu = 0$ this is equation is solved by
	\[ \overline{y}^\pm = \pm \sqrt{- \frac{\ell}{\sum_{\sigma, \tau \in \Sigma} f_{\sigma\tau}}}. \]
	Furthermore, $\frac{\partial}{\partial y} g(\overline{y}^\pm, 0) = 2 \sum_{\sigma, \tau \in \Sigma} f_{\sigma\tau} \overline{y}^\pm$, which, generically, does not equal $0$. Hence, by the implicit function theorem, we obtain two branches of solutions
	\[ \overline{Y}^\pm (\mu) = \overline{y}^\pm + \OO(|\mu|). \]
	Transforming back into the original variables we obtain two branches
	\[ x_p (\lambda) = \pm \sqrt{- \frac{\ell}{\sum_{\sigma, \tau \in \Sigma} f_{\sigma\tau}}} \cdot \sqrt{\lambda} + \OO(|\lambda|) \]
	for small $\lambda>0$, which completes the proof for the first case.
	
	The case $\ell / \sum_{\sigma, \tau \in \Sigma} f_{\sigma\tau} <0$ is analogous. We introduce $x_p = \mu y$ with $\mu = \sqrt{-\lambda}$ for $\lambda<0$ with small absolute value. Then the proof is the same.
\end{proof}
\begin{remk}
	\begin{enumerate}[label=(\roman*)]
		\item The critical maximal cells simultaneously undergo a saddle node bifurcation, that is two steady state branches exist for either only positive or only negative values of $\lambda$. The sign of $\lambda$ for which branching solutions occur is the same for all maximal cells. We refer to the former as the \emph{supercritical} and to the latter as the \emph{subcritical} case. Note that a generic $f$ satisfying the bifurcation conditions (B) always fulfills one of the two assumptions from the previous lemma, as generically $\ell, \sum_{\sigma, \tau \in \Sigma} f_{\sigma\tau} \ne 0$. However, we see that the generic bifurcation behavior is different in different regions of system parameter space, i.e. partial derivatives of $f$.
		\item Note, furthermore, that in both cases the equations for maximal cells are completely uncoupled. Hence, for a specific branch of solutions not all maximal cells need to evolve according to the same branch. In particular the choice of sign in $D_p^\pm$ may differ in different maximal cells. As a result, globally, when restricting only to maximal cells, we obtain $2^m$ branches of solutions where $m$ is the number of maximal cells.
		\hspace*{\fill}$\triangle$
	\end{enumerate}
\end{remk}

For the non-maximal cells, we proceed inductively. We assume to know a specific branching pattern for all cells above a given cell $p$ and compute the solutions for that cell. Hence, we need to distinguish between the super- and subcritical cases.
\begin{lemma}[supercritical case]
	\label{lem:c1-st2.1}
	Let $p$ be non-maximal. Assume for all $q\nog p$
	\[ x_q(\lambda) = d_q \cdot \sqrt{\lambda} + \OO(|\lambda|) \]
	for small $\lambda>0$ and some $d_q \in \RR \setminus \{0\}$. Then
	\[ x_p(\lambda) = D_p \cdot \sqrt{\lambda} + \OO(|\lambda|) \]
	for small $\lambda >0$, where
	\[ D_p = - \frac{\sum_{\tau \notin \LL_p} a_\tau d_{\tau(p)}}{\sum_{\sigma \in \LL_{p}} a_\sigma} = \frac{\sum_{\tau \notin \LL_p} a_\tau d_{\tau(p)}}{\sum_{\tau \notin \LL_{p}} a_\tau} \ne 0. \]
\end{lemma}
\begin{proof}
	As the maximal cells are critical and $p$ is non-maximal, $p$ is non-critical. Using the assumption on all $q \nog p$, \eqref{eq:taylor} becomes
	\[ 0 = \sum_{\sigma \in \LL_{p}} a_\sigma x_p + \sum_{\substack{p \nol q = \tau(p) \colon \\ \tau \notin \LL_p}} a_\tau d_q \sqrt{\lambda} + \OO\left( |x_p|^2 + |\lambda| + |x_p| \cdot \sqrt{|\lambda|} \right). \]
	Since $p$ is non-critical, $\sum_{\sigma \in \LL_{p}} a_\sigma \ne 0$. Thus, by the implicit function theorem, we obtain that this equation is uniquely solved by
	\[ x_p (\lambda) = - \frac{\sum_{\tau \not\in \LL_p} a_\tau d_{\tau(p)}}{\sum_{\sigma \in \LL_{p}} a_\sigma} \cdot \sqrt{\lambda} + \OO\left(|\lambda|\right). \]
	The second representation of the coefficient follows from \eqref{eq:maxcrit} which completes the proof.
\end{proof}

\begin{lemma}[subcritical case]
	\label{lem:c1-st2.2}
	Let $p$ be non-maximal. Assume for all $q\nog p$
	\[ x_q(\lambda) = d_q \cdot \sqrt{-\lambda} + \OO(|\lambda|) \]
	for small $\lambda<0$ and some $d_q \in \RR \setminus \{0\}$. Then
	\[ x_p(\lambda) = D_p \cdot \sqrt{- \lambda} + \OO(|\lambda|) \]
	for small $\lambda <0$, where
	\[ D_p = - \frac{\sum_{\tau \not\in \LL_p} a_\tau d_{\tau(p)}}{\sum_{\sigma \in \LL_{p}} a_\sigma} = \frac{\sum_{\tau \notin \LL_p} a_\tau d_{\tau(p)}}{\sum_{\tau \notin \LL_{p}} a_\tau} \ne 0. \]
\end{lemma}
\begin{proof}
	The proof is completely analogous to the previous one.
\end{proof}

The maximal cells simultaneously determine whether a bifurcation occurs super- or subcritically. The non-maximal cells have no further influence, as we have seen in the previous two lemmas. We may, therefore, perform an inductive proof with respect to $\nole$, summarizing \Cref{lem:c1-st1,lem:c1-st2.1,lem:c1-st2.2}, to obtain
\begin{theorem}
	\label{thm:maxcrit}
	Under the bifurcation assumption \emph{(B)} and assuming the maximal cells in $C$ to be critical, the state variables $x_p$ for all $p\in C$ bifurcate according to one of the following two asymptotics.
	\begin{enumerate}[label=(\roman*)]
		\item (supercritical) 
		\[ x_p(\lambda) = D_p \cdot \sqrt{\lambda} + \OO(|\lambda|) \quad \text{for small} \quad \lambda>0, \quad \text{if} \quad \frac{\ell}{\sum_{\sigma, \tau \in \Sigma} f_{\sigma\tau}} <0; \]
		\item (subcritical)
		\[ x_p(\lambda) = D_p \cdot \sqrt{- \lambda} + \OO(|\lambda|) \quad \text{for small} \quad \lambda<0, \quad \text{if} \quad \frac{\ell}{\sum_{\sigma, \tau \in \Sigma} f_{\sigma\tau}} >0. \]
	\end{enumerate}
	Therein the coefficients $D_p$ 
	are defined recursively. For $p$ maximal they are
	\begin{enumerate}[label=(\roman*)]
		\item 
		(supercritical)
		\[ D_p \in \left\{ +\sqrt{- \frac{\ell}{\sum_{\sigma,\tau\in\Sigma} f_{\sigma\tau}}}, -\sqrt{- \frac{\ell}{\sum_{\sigma,\tau\in\Sigma} f_{\sigma\tau}}} \right\}; \]
		\item 
		(subcritical)
		\[ D_p \in \left\{ +\sqrt{\frac{\ell}{\sum_{\sigma,\tau\in\Sigma} f_{\sigma\tau}}}, -\sqrt{\frac{\ell}{\sum_{\sigma,\tau\in\Sigma} f_{\sigma\tau}}} \right\}. \]
	\end{enumerate}
	The remaining (non-maximal) ones are defined via
	\[ D_p = - \frac{\sum_{\tau \notin \LL_p} a_\tau D_{\tau(p)}}{\sum_{\sigma \in \LL_{p}} a_\sigma} = \frac{\sum_{\tau \notin \LL_p} a_\tau D_{\tau(p)}}{\sum_{\tau \notin \LL_{p}} a_\tau}. \]
	There are $2^m$, with $m = \# \{ p \in C \mid p \text{ maximal} \}$, different branches of steady states in both cases, that are determined by the choices of $D_p$ for $p$ maximal.
\end{theorem}
\begin{remk}
	The branching solutions in \Cref{thm:maxcrit} are the same as the ones described in Proposition 5.1 in \textcite{Soares.2018} for layered feedforward networks (compare to \Cref{rem:layered}) and in Proposition~5.7 of \textcite{Aguiar.2019c} investigating feedforward structure of transitive components. Here we extend these results by the explicit computation of the leading coefficients.
	\hspace*{\fill}$\triangle$
\end{remk}
\begin{remk}
	Note that by restriction to the invariant fully synchronous subspace, we obtain a fully synchronous saddle node bifurcation similar to the proof of \Cref{lem:c1-st1}. This will be made more precise in the case of non-critical maximal cells (see \Cref{lem:abff_maxnoncritfullysynch} in \Cref{subsec:maxnoncrit}) where it is of great importance. Hence two of the branching solutions provided by \Cref{thm:maxcrit} necessarily describe this branch. It can easily be seen from the recursive formulas using \eqref{eq:maxcrit} that these are exactly the ones where the coefficients of maximal cells all have the same sign. This is also to be expected, as in all other cases not even the maximal cells are synchronous.
	\hspace*{\fill}$\triangle$
\end{remk}
\begin{remk}
	The results of \Cref{thm:maxcrit} are only fully accurate if none of the $D_p$ vanish. This can be seen to be the case generically for the maximal cells. The system parameter $\ell$ does not vanish generically. Hence, $D_p \ne 0$ for $p$ maximal. In particular, the coefficients of the fully synchronous branching solutions do not vanish generically. For the remaining cells, this has to be checked in the inductive computation. The structure of the network could force some $D_p$ to vanish identically, a case that we do not expect to occur and have never encountered in any examples. Furthermore, the branching statement holds true generically in the sense that the maximal cells undergo a saddle note bifurcation and there is a unique solution for \eqref{eq:taylor} for $p$ non-maximal due to the implicit function theorem. In the case that $D_p$ vanishes, the leading square root order increases. 
	
	Finally, additional structure in the network can give us the means to prove that the leading coefficients do not vanish generically for all cells. For example, this holds true if there is an input map $\kappa\in\Sigma$ such that $\kappa(p)$ is maximal for all $p$. This is the case if the network is a `semigroup network' in the language of \textcite{Rink.2014}.
	
	\hspace*{\fill}$\triangle$
\end{remk}

\subsection{The critical cells are non-maximal}
\label{subsec:maxnoncrit}
Next, we assume that the maximal cells are non-critical. In particular, $\sum_{\sigma \in \Sigma}a_\sigma \ne 0$ under the condition of genericity. The general strategy for finding branching solutions remains the same as in the previous part. We solve \eqref{eq:taylor} for a given cell $p$ assuming knowledge of its inputs. However, the considerations, especially for non-maximal cells, become more involved, as we have to distinguish whether a cell is critical in each step. Inductively, this provides all possible solutions for all cells. As a result of the multitude of different cases, the explicit bifurcation patterns -- governed by simultaneous solutions of \eqref{eq:taylor} for all cells -- become a lot more complex than when maximal cells are critical. Once again, we have to distinguish branching solutions that exist for positive and negative values of the bifurcation parameter $\lambda$. We refer to these cases as \emph{super-} and \emph{subcritical}, as before. The result for the subcritical case, however, can be obtained as a corollary from the supercritical case, as we see in \Cref{thm:abff_maxnoncritsub}.

\subsubsection{Root subnetworks}
As we will see, branches of steady states for the bifurcation problem (B) are determined by subnetworks that are `surrounded' by critical cells. Similar to the maximal cells in \Cref{subsec:maxcrit}, these serve as the starting point for the inductive investigation of the entire network. All cells within the subnetwork remain synchronous. Non-trivial solutions branch off in cells that are not in those subnetworks. More precisely, the non-trivial solutions will be higher-order saddle node branches whose amplitude increases the lower a cell is in the network.

\begin{defi}[Root subnetwork]
	Let $C$ be the set of cells of a feedforward network and assume the bifurcation scenario (B). We call a non-trivial subnetwork $\emptyset \ne B \subsetneq C$ a \emph{root subnetwork}, if it contains all maximal cells and, if for every cell $p \in C \setminus B$ such that $q \in B$ for all $q \nog p$, it holds that $p$ is critical.
\end{defi}
\begin{ex}
	In the classical $N$ cell feedforward chain a subnetwork contains the first $k$ cells for any $1\le k\le N$. Furthermore, there are exactly two loop-types given by the maximal cell on the one hand and all other cells on the other. Hence, for non-maximal critical cells it can readily be seen that each subnetwork is a root subnetwork. As was shown in \cite{Rink.2013}, all generic branches of steady states are of the form that there is $1 \le k \le N$ such that the first $k$ cells in the chain remain synchronous, while the states of the remaining cells branch off in higher order saddle node branches whose amplitude increases with the distance to the synchronous cells.
\end{ex}

Solutions in the root subnetworks behave as the following fully synchronous branch.
\begin{lemma}
	\label{lem:abff_maxnoncritfullysynch}
	Recall that the bifurcation assumption \emph{\ref{enum:B1}} implies that there exists a fully synchronous steady state at the bifurcation point: $\gamma_f(0,0)=0$. If maximal cells are non-critical, this fully synchronous steady state persists under (small) parameter variations. It grows as
	\begin{equation}
		\label{eq:abff_maxnoncritfullysynchX}
		x_p (\lambda) = X(\lambda) = D \lambda + R \lambda^2 + \OO \left( |\lambda|^3 \right)
	\end{equation}
	for $|\lambda|$ small and all $p \in C$. Therein
	\begin{align}
		D &= - \frac{\ell}{\sum_{\sigma \in \Sigma} a_\sigma},	\label{eq:abff_maxnoncritfullysynchD} \\
		R &= - \frac{\sum_{\sigma, \tau \in \Sigma} f_{\sigma \tau} \ell^2 - \sum_{\sigma\in\Sigma} a_\sigma \sum_{\sigma \in \Sigma} f_{\sigma\lambda} \ell + \left( \sum_{\sigma \in \Sigma} a_\sigma \right)^2 f_{\lambda\lambda}}{\left( \sum_{\sigma \in \Sigma} a_\sigma \right)^3},	\label{eq:abff_maxnoncritfullysynchR}
	\end{align}
	which generically do not vanish.
\end{lemma}
\begin{proof}
	Consider the fully synchronous subspace $\{ x_{p_1} = \dotso = x_{p_N} \} \subset \bigoplus_{p \in C} \RR$. This can readily be seen to be invariant under the flow induced by network vector fields of the form \eqref{eq:netvfparam}. Choosing a coordinate $y$ for this subspace, the bifurcation problem becomes the same for all cells $p \in C$. The Taylor expanded equation \eqref{eq:taylor} is
	\[ 0 = \sum_{\sigma \in \Sigma} a_\sigma y + \ell \lambda + \OO(|y|^2 + |y||\lambda| + |\lambda|^2). \]
	The implicit function theorem yields a unique branch of solutions
	\[ x_p (\lambda) = y(\lambda) = -\frac{\ell}{\sum_{\sigma \in \Sigma} a_\sigma} \lambda + \OO \left( |\lambda|^2 \right) \]
	for $f(y, \dotsc, y, \lambda) = 0$ with $|\lambda|$ small. Performing second order implicit differentiation -- for which we omit the details --, we compute
	\[ y'' (\lambda) = - 2 \cdot \frac{\sum_{\sigma, \tau \in \Sigma} f_{\sigma \tau} \ell^2 - \sum_{\sigma\in\Sigma} a_\sigma \sum_{\sigma \in \Sigma} f_{\sigma\lambda} \ell + \left( \sum_{\sigma \in \Sigma} a_\sigma \right)^2 f_{\lambda\lambda}}{\left( \sum_{\sigma \in \Sigma} a_\sigma \right)^3}. \]
	As both derivatives generically do not vanish, this completes the proof.
\end{proof}
\begin{cor}
	\label{cor:abff_subnetfullsynch}
	Let $B \subset C$ be a subnetwork and assume bifurcation assumption \emph{(B)} holds with non-maximal critical cells. For all cells $p \in B$ there exists a synchronous branch of steady states
	\begin{equation}
		\label{eq:abff_subnetfullsynch}
		x_p (\lambda) = X(\lambda) = D \lambda + R \lambda^2 + \OO \left( |\lambda|^3 \right)
	\end{equation}
	for $|\lambda|$ small solving the bifurcation equation \eqref{eq:taylor}. Therein
	\begin{align}
		D &= - \frac{\ell}{\sum_{\sigma \in \Sigma} a_\sigma},	\label{eq:abff_subnetfullsynchD} \\
		R &= - \frac{\sum_{\sigma, \tau \in \Sigma} f_{\sigma \tau} \ell^2 - \sum_{\sigma\in\Sigma} a_\sigma \sum_{\sigma \in \Sigma} f_{\sigma\lambda} \ell + \left( \sum_{\sigma \in \Sigma} a_\sigma \right)^2 f_{\lambda\lambda}}{\left( \sum_{\sigma \in \Sigma} a_\sigma \right)^3},	\label{eq:abff_subnetfullsynchR}
	\end{align}
	which generically do not vanish.
\end{cor}
\begin{proof}
	Restrict the system of bifurcation equations \eqref{eq:Barb} to the cells of the subnetwork $B$, i.e. discard all the equations for variables $x_p$ with $p \notin B$, and apply \Cref{lem:abff_maxnoncritfullysynch}.
\end{proof}
\begin{remk}
	Note that the previous corollary does not make any claims concerning uniqueness of the solution branch -- except inside the synchrony space. It also does not cover the impact of cells $p \in C\setminus B$, which might prevent the existence of this branch for the entire network, as we will see later.
	\hspace*{\fill}$\triangle$
\end{remk}
\begin{lemma}
	\label{lem:abff_maxnoncritnonmaxcritT}
	Let $p \in C$ be non-maximal and critical. Assume that all cells $q \nog p$ are in the fully synchronous state $x_q(\lambda)=X(\lambda)$. Define
	\begin{align*}
		A	&= \sum_{\sigma, \tau \in \LL_p} f_{\sigma\tau}, \\
		B	&= \sum_{\sigma \in \LL_{p}} f_{\sigma \lambda} + 2 \sum_{\sigma \in \LL_p, \tau \notin \LL_p } f_{\sigma \tau} D, \\
		C	&= \sum_{\tau \notin \LL_p} a_\tau R + \sum_{\tau \notin \LL_p} f_{\tau \lambda} D + \sum_{\sigma, \tau \notin \LL_p} f_{\sigma \tau} D^2 + f_{\lambda \lambda},
	\end{align*}
	with $D,R$ as in \eqref{eq:abff_maxnoncritfullysynchD} and \eqref{eq:abff_maxnoncritfullysynchR}. The solutions to \eqref{eq:taylor} generically undergo a transcritical bifurcation
	\[ x_p (\lambda) = \Dsup_p^\pm \lambda + \OO\left( |\lambda|^2 \right), \]
	for $|\lambda|$ small, where
	\begin{align*}
		\Dsup_p^+ &= -\frac{\ell}{\sum_{\sigma \in \Sigma} a_\sigma} \\
		\Dsup_p^- &= \dfrac{\ell}{\sum_{\sigma \in \Sigma} a_\sigma} \cdot \left( 1 + 2 \dfrac{\sum_{\sigma \in \LL_p, \tau \notin \LL_p} f_{\sigma\tau}}{\sum_{\sigma, \tau \in \LL_p} f_{\sigma\tau}} \right) - \dfrac{\sum_{\sigma \in \LL_p} f_{\sigma\lambda}}{\sum_{\sigma, \tau \in \LL_p} f_{\sigma\tau}}.
	\end{align*}
	We will see in the proof of \Cref{lem:rootgen} below that the branch with the coefficient $\Dsup_p^+$ is the fully synchronous branch $X(\lambda)$.
\end{lemma}
\begin{proof}
	The bifurcation assumption (B) implies
	\[ \sum_{\tau \not\in \LL_p} a_\tau D = \sum_{\tau \in \Sigma} a_\tau D = - \sum_{\tau \in \Sigma} a_\tau\frac{\ell}{\sum_{\sigma \in \Sigma} a_\sigma} = - \ell, \]
	as $\sum_{\sigma \in \LL_p} a_\sigma =0$. In particular, the linear terms in $\lambda$ in \eqref{eq:taylor} vanish and the equation becomes
	\begin{equation*}
		\begin{split}
			0	&= \sum_{\sigma, \tau \in \LL_p} f_{\sigma\tau} x_p^2 + \sum_{\sigma \in \LL_{p}} f_{\sigma \lambda} \lambda x_p +  \sum_{\tau \not\in \LL_p} a_\tau R \lambda^2 + \sum_{\tau \not\in \LL_p} f_{\tau \lambda} D \lambda^2 \\
			&\phantom{=} + 2 \sum_{\sigma \in \LL_p, \tau \notin \LL_p} f_{\sigma \tau} D \lambda x_p + \sum_{\sigma, \tau \notin \LL_p} f_{\sigma \tau} D^2 \lambda^2 + f_{\lambda \lambda} \lambda^2 + \OO \left( |x_p|^3 + |x_p|^2|\lambda| + |x_p||\lambda|^2 + |\lambda|^3 \right) \\
			&= A x_p^2 + B \lambda x_p + C \lambda^2 + \OO \left( |x_p|^3 + |x_p|^2|\lambda| + |x_p||\lambda|^2 + |\lambda|^3 \right).
		\end{split}
	\end{equation*}
	Similar to the previous proofs, we employ the standard method to detect transcritical bifurcations (see for example \textcite{Murdock.2003}) by introducing a new variable $x_p = \lambda y$. The equation becomes
	\begin{align*}
		0	&= A \lambda^2 y^2 + B \lambda^2 y + C \lambda^2 + \OO \left( |y||\lambda|^3 + |\lambda|^3 \right) \\
		&= \lambda^2 \left( A y^2 + B y + C + \OO \left( |y||\lambda| + |\lambda| \right) \right) = \lambda^2 g(y,\lambda).
	\end{align*}
	Note that the coefficients $A,B$ and $C$ are composed of the second order partial derivatives of $f$. Hence, generically $A, B^2- 4AC \ne 0$. By \Cref{lem:abff_maxnoncritDcritsuper} in the appendix, generically $B^2 - 4AC >0$ and we obtain two solutions:
	\[ \overline{y}^\pm = \frac{-B \pm \sqrt{B^2 - 4AC}}{2A}. \]
	Furthermore, $\frac{\partial}{\partial y} g(\overline{y}^\pm, 0) = 2A \overline{y}^\pm + B$ which, generically, does not vanish. Hence, by the implicit function theorem, we obtain two branches of solutions
	\[ \overline{Y}^\pm (\lambda) = \overline{y}^\pm + \OO(|\lambda|) \]
	for small $\lambda>0$. Transforming back into the original coordinates, we obtain
	\[ x_p (\lambda) = \frac{-B \pm \sqrt{B^2 - 4AC}}{2A} \cdot \lambda + \OO\left(|\lambda|^2\right). \]
	The computation of the leading coefficients can be found in \Cref{lem:abff_maxnoncritDcritsuper} in the appendix.
\end{proof}

\begin{lemma}
	\label{lem:rootgen}
	For every generically existing branch of bifurcating steady states of \eqref{eq:Barb} with non-maximal critical cells, there is a unique root subnetwork $B \subset C$ such that all $x_p$ with $p\in B$ remain synchronous, i.e.,
	\[ x_p(\lambda) = X(\lambda) \quad \text{for all} \quad p \in B \]
	for $|\lambda|$ small and $X(\lambda)$ as in \Cref{lem:abff_maxnoncritfullysynch}. Furthermore, for each $p \in C\setminus B$ such that $q \in B$ for all $q \nog p$ the steady state solution of \eqref{eq:Barb} grows as
	\[ x_p(\lambda) = \Dsup_p^- \lambda + \OO(|\lambda|^2) \]
	for $|\lambda|$ small with
	\[
	\Dsup_p^- = \dfrac{\ell}{\sum_{\sigma \in \Sigma} a_\sigma} \cdot \left( 1 + 2 \dfrac{\sum_{\sigma \in \LL_p, \tau \notin \LL_p} f_{\sigma\tau}}{\sum_{\sigma, \tau \in \LL_p} f_{\sigma\tau}} \right) - \dfrac{\sum_{\sigma \in \LL_p} f_{\sigma\lambda}}{\sum_{\sigma, \tau \in \LL_p} f_{\sigma\tau}}.
	\]
\end{lemma}
\begin{proof}
	We prove the statement by investigating the bifurcation equations \eqref{eq:taylor} and constructing the subnetwork $B$ inductively. We begin with the maximal cells. As $\sigma(p)=p$ for each maximal $p$ and for each $\sigma\in\Sigma$, we see that
	\[ B_0 = \{ p \mid p \text{ is maximal}\} \]
	is a subnetwork. The bifurcation equation \eqref{eq:taylor} for maximal cells is the same as the one restricted to the fully synchronous subspace in the proof of \Cref{lem:abff_maxnoncritfullysynch} and the implicit function theorem guarantees the existence of a unique branch of solutions. As $B_0$ is a subnetwork, \Cref{cor:abff_subnetfullsynch} implies that this is necessarily the synchronous branch of solutions $x_p(\lambda)=X(\lambda)$ (note that this also follows from the fact that \eqref{eq:taylor} for maximal cells is the same as the one restricted to the fully synchronous subspace in the proof of \Cref{lem:abff_maxnoncritfullysynch}). As the maximal cells do not depend on any other cells, their state variables always branch according to this solution, independent of the branching behavior of the entire network. In particular, we may assume the existence of a subnetwork $B$ of synchronous cells for any branch from now on, since $B_0 \subset B$.
	
	It remains to show the properties of $B$ and of cells $p \in C\setminus B$ with $q \in B$ for all $q \nog p$. To that end let us assume the existence of a subnetwork $\tilde{B} \subset C$ such that $x_p(\lambda) = X(\lambda)$ for all $p \in \tilde{B}$. Then we consider a cell $p \in C \setminus \tilde{B}$ such that $q \in \tilde{B}$ for all $q \nog p$. We consider two cases.
	
	First, assume $p$ is non-critical. Similar to before, the derivative of \eqref{eq:taylor} with respect to $x_p$ equals $\sum_{\sigma \in \LL_p} a_\sigma $ which does not vanish generically according to bifurcation assumption (B). By the implicit function theorem the equation has a unique solution $x_p(\lambda)$ for $|\lambda|$ small. On the other hand, also $\tilde{B} \cup \{p\}$ is a subnetwork which contains $\tilde{B}$ as a subnetwork itself. Furthermore, equation \eqref{eq:taylor} for $p$ is the same in the full network as well as in $\tilde{B} \cup \{p\}$. Hence, according to \Cref{cor:abff_subnetfullsynch} the solution for cell $p$ is necessarily the same as the synchronous branch, $x_p(\lambda)=X(\lambda)$, and we may set $\tilde{B} \cup \{p\} \subset B$. In particular, if $p$ is non-critical and $q \in B$ for all $q \nog p$, then also $p \in B$.
	
	Second, we assume that $p$ is critical. The investigation is similar. We are in the situation of \Cref{lem:abff_maxnoncritnonmaxcritT}, which tells us that generically the bifurcation equation \eqref{eq:taylor} for cell $p$ has precisely two solution branches
	\[ x_p^\pm(\lambda) = D_p^\pm \lambda + \OO\left( |\lambda|^2 \right) \]
	with coefficients
	\begin{align*}
		D_p^+ &= -\frac{\ell}{\sum_{\sigma \in \Sigma} a_\sigma} \\
		D_p^- &= \dfrac{\ell}{\sum_{\sigma \in \Sigma} a_\sigma} \cdot \left( 1 + 2 \dfrac{\sum_{\sigma \in \LL_p, \tau \notin \LL_p} f_{\sigma\tau}}{\sum_{\sigma, \tau \in \LL_p} f_{\sigma\tau}} \right) - \dfrac{\sum_{\sigma \in \LL_p} f_{\sigma\lambda}}{\sum_{\sigma, \tau \in \LL_p} f_{\sigma\tau}}.
	\end{align*}
	On the other hand, also $\tilde{B} \cup \{p\}$ is a subnetwork which contains $\tilde{B}$ as a subnetwork itself. Furthermore, equation \eqref{eq:taylor} for $p$ is the same in the full network as well as in $\tilde{B} \cup \{p\}$. Hence, according to \Cref{cor:abff_subnetfullsynch} the synchronous solution for cell $p$ is also generic in the full network. Hence, necessarily one of the branching solutions $x_p^\pm (\lambda)$ is the same as the synchronous branch. As was shown in \Cref{lem:abff_maxnoncritnonmaxcritT}, the leading coefficient generically only matches the leading coefficient of the synchronous branch for one choice of sign: $x_p^+(\lambda)=X(\lambda)$. In that case we once again set $\tilde{B} \cup \{p\} \subset B$. For the other choice we generically have $x_p^-(\lambda) \not\equiv X(\lambda)$ so that $p \notin B$. This completes the proof.
\end{proof}
\begin{remk}
	Note that the previous lemma does not make any claims concerning existence of a solution branch for each root subnetwork $B\subset C$. In particular, it does not cover the impact of cells $p \in C\setminus B$, which might prevent the existence of this branch for the entire network, as we will see later. Furthermore, the result does not exclude the possibility that a cell $p \in C \setminus B$ branches according to the synchronous branch of solutions. However, we will see below that this is not to be expected generically.
	\hspace*{\fill}$\triangle$
\end{remk}
\begin{remk}
	The argument identifying the synchronous branch in the previous proof is a special case of quiver symmetry. In particular, the subnetworks $B$ and $B \cup \{p\}$ induce two non-classical symmetries each. On one hand, there is the inclusion of the total phase space of the subnetworks into the total phase space of the original network. On the other hand, each subnetwork gives rise to a quotient network by identifying all cells within the subnetwork. The dynamics on the original network respects these maps, which in turn guarantees genericity of the synchronous branch in the equations corresponding to the subnetwork. Note that these symmetries are parts of larger structures, namely the \emph{subnetwork quiver} and the \emph{quotient quiver}. For more details, see \cite{Nijholt.2020}.
	\hspace*{\fill}$\triangle$
\end{remk}

\subsubsection[Solving Equation (3.4)]{Solving \Cref{eq:taylor}}
In this subsection we set up the technical tools for the inductive proof of the bifurcation result for the entire network. In particular, we solve \eqref{eq:taylor} in individual cells under specific assumptions on the branching solutions for cells above with respect to $\nole$.
%
As we have seen in \Cref{lem:rootgen}, any generically existing branch for the entire network is determined by a root subnetwork $B\subset C$ in which the cells evolve according to the fully synchronous branch. Hence, we fix such a root subnetwork. Furthermore, we focus on $\lambda \ge 0$ and assume the following input scenarios for a fixed non-maximal cell $p \in C\setminus B$ as an inductive hypothesis:
\begin{itemize}
	\item[(\textbf{H})] For all $q \nog p$ and small $\lambda>0$ the solution to \eqref{eq:taylor} has the asymptotics
	\[ x_q (\lambda) = \dsup_q \cdot \lambda^{2^{-\xi_q}} + \OO\left(|\lambda|^{2^{-(\xi_q-1)}}\right), \]
	where $\dsup_q \in \RR \setminus \{ 0 \}$ and the $\xi_q$ are integers with $0 \le \xi_q$ that define the square root order of the branching solution of cell $q$.
\end{itemize}
Under the assumption (\textbf{H}) we define the quantity
\begin{equation}
\label{eq:abff_xi}
\xisup_p = \max_{q \in\SSigma(p)} \xi_q,
\end{equation}
to be the highest (and thus `leading') square root order of inputs into cell $p$. Then, $\xisup_p = 0$, if and only if all inputs into cell $p$ evolve linearly in $\lambda$ up to leading order. To further simplify, we denote the subset of cells $q \in\SSigma(p)$ which are of highest square root order in $\lambda$ by $\Qsup_p$. That is
\begin{equation}
	\label{eq:abff_Qsup}
	\Qsup_p = \left\{ q \in\SSigma(p)\ \left|\ x_q (\lambda) = \dsup_q \cdot \lambda^{2^{-\xisup_p}} + \OO\left(|\lambda|^{2^{-(\xisup_p-1)}}\right) \text{ and } \dsup_q\ne 0 \right. \right\}.
\end{equation}
In the case $\xisup_p=0$ all state variables $x_q$ for $q \in\SSigma(p)$ evolve linearly in $\lambda$ up to leading order. Hence, $\Qsup_p = \SSigma(p)$. Note that $\xisup_p$ and $\Qsup_p$ are only defined for $p$ non-maximal.

Depending on $\xisub_p$, we formulate the following non-degenericity conditions:
\begin{itemize}
	\item[(L)]
	\[ \xisup_p=0 \quad \text{and} \quad \sum_{\tau \not\in \LL_p} a_\tau \dsup_{\tau(p)} + \ell \ne 0, \]
	\item[(SN)] \hfill 
	\[ \xisup_p>0 \quad \text{and} \quad \sum_{\tau \colon \tau(p)\in \Qsup_p} a_\tau \dsup_{\tau(p)} \ne 0. \]
\end{itemize}
The conditions (L) and (SN) guarantee that the leading order terms in $\lambda$ in \eqref{eq:taylor} do not vanish.
In what follows, we observe that the two cases lead to a continuation of the trivial solution at the bifurcation point, which is linear up to leading order (L) and a higher order saddle node bifurcation (SN) respectively.

Under the given inductive assumptions, we prove statements providing branching solutions to \eqref{eq:taylor} for non-maximal cells. The technical proofs are postponed until the appendix \Cref{sec:ap1}. We start with the case that $p$ is non-critical. 
\begin{lemma}
	\label{lem:abff_maxnoncritnonmaxnoncrit}
	Let $p \in C\setminus B$ be non-critical. Under assumptions \emph{(\textbf{H})} and \emph{(L)} or \emph{(SN)}, depending on the value of $\xisup_p$, \eqref{eq:taylor} has the unique solution
	\[ x_p (\lambda) = \dsup_p \cdot \lambda^{2^{-\xisup_p}} + \OO\left(|\lambda|^{2^{-(\xisup_p-1)}}\right) \]
	for $\lambda > 0$ small, where $\dsup_p \ne 0$ is given by
	\[ \dsup_p = \begin{cases}
	- \dfrac{\sum_{\tau \notin \LL_p} a_\tau \dsup_{\tau(p)} + \ell}{\sum_{\sigma \in \LL_{p}} a_\sigma},	&\text{if} \quad \xisup_p = 0; \\[15pt]
	- \dfrac{\sum_{\tau \colon \tau(p) \in \Qsup_p} a_\tau \dsup_{\tau(p)}}{\sum_{\sigma \in \LL_{p}} a_\sigma},		&\text{if} \quad \xisup_p >0.
	\end{cases} \]
\end{lemma}

The corresponding result for critical non-maximal cells is proved in multiple lemmas. We distinguish between the cases (L) and (SN).
\begin{lemma}
	\label{lem:abff_maxnoncritnonmaxcritL}
	Let $p \in C\setminus B$ be critical. Assume \emph{(\textbf{H})} and \emph{(L)} to hold true. The solutions to \eqref{eq:taylor} generically bifurcate as in one of the following two cases:
	\begin{enumerate}[label=(\roman*)]
		\item If
		\[ \frac{\sum_{\tau \notin \LL_p} a_\tau \dsup_{\tau(p)} + \ell}{\sum_{\sigma,\tau \in \LL_{p}} f_{\sigma\tau}} >0, \]
		there are no branching solutions.
		\item If
		\[ \frac{\sum_{\tau \notin \LL_p} a_\tau \dsup_{\tau(p)} + \ell}{\sum_{\sigma,\tau \in \LL_{p}} f_{\sigma\tau}} <0, \]
		the solutions undergo a saddle node bifurcation
		\[ x_p (\lambda) = \dsup_p^\pm \cdot \sqrt{\lambda} + \OO\left(|\lambda|\right) \]
		for $\lambda > 0$ small, where
		\[ \dsup_p^\pm = \pm \sqrt{- \frac{\sum_{\tau \notin \LL_p} a_\tau \dsup_{\tau(p)} + \ell}{\sum_{\sigma,\tau \in \LL_{p}} f_{\sigma\tau}}} \ne 0. \]
	\end{enumerate}
\end{lemma}


Now, we turn to the branching solutions for critical cells $p$ with the additional non-degenericity condition (SN).
\begin{lemma}
	\label{lem:abff_maxnoncritnonmaxcritSN}
	Let $p \in C \setminus B$ be critical. Assume the induction hypothesis \emph{(\textbf{H})} and the additional non-degenericity condition \emph{(SN)} to hold true. The solutions to \eqref{eq:taylor} generically bifurcate as in one of the following two cases:
	\begin{enumerate}[label=(\roman*)]
		\item If
		\[ \frac{\sum_{\tau \colon \tau(p) \in \Qsup_p} a_\tau \dsup_{\tau(p)}}{\sum_{\sigma, \tau \in \LL_p} f_{\sigma\tau}} > 0, \]
		there are no branching solutions.
		\item If
		\[ \frac{\sum_{\tau \colon \tau(p) \in \Qsup_p} a_\tau \dsup_{\tau(p)}}{\sum_{\sigma, \tau \in \LL_p} f_{\sigma\tau}} < 0, \]
		the solutions undergo a higher order saddle node bifurcation as
		\[ x_p (\lambda) = \dsup_p^\pm \cdot \lambda^{2^{-(\xisup_p+1)}} + \OO\left( |\lambda|^{2^{-\xisup_p}} \right), \]
		for $\lambda>0$ small, where
		\[ \dsup_p^\pm = \pm \sqrt{- \frac{\sum_{\tau \colon \tau(p) \in \Qsup_p} a_\tau \dsup_{\tau(p)}}{\sum_{\sigma, \tau \in \LL_p} f_{\sigma\tau}}} \ne 0. \]
	\end{enumerate}	
\end{lemma}

\subsubsection{Branches of steady states for the entire network}
The results in the previous section form the technical background for the inductive investigation of branching solutions to the bifurcation problem (B) with non-maximal critical cells. Subtleties arise while investigating which cases from \Cref{lem:abff_maxnoncritnonmaxnoncrit,lem:abff_maxnoncritnonmaxcritL,lem:abff_maxnoncritnonmaxcritSN} can generically occur, when \eqref{eq:taylor} is solved for all $p \in C$ simultaneously. 

For a given root subnetwork $\Bsup \subset C$, the branching behavior of cells $p \notin \Bsup$ is determined by the number of critical cells `in between' $\Bsup$ and $p$. In particular, we will need the quantity
\begin{equation}
	\label{eq:musup}
	\musup_p = \max_{\overline{p} \in \Bsup} \max_{\omega \in \Omega_{\overline{p}, p}} \# \left\{q \in \omega \mid q \text{ critical}, q \notin \Bsup \right\} - 1,
\end{equation}
which is the maximal number of critical cells $q \notin \Bsup$ along paths from any cell in $\Bsup$ to $p$ (recall that $\Omega_{\overline{p}, p}$ denotes the set of all paths from $\overline{p}$ to $p$ without any loops). For convenience we set $\mu_p=0$ for $p \in B$. It can readily be seen via induction on the partial order $\nole$ that $\mu_p$ can alternatively be characterized iteratively:
\begin{equation}
	\label{eq:muiteratively}
	\musup_p = \begin{cases}
		0								\qquad &\text{for} \quad p \in B, \\
		0								\qquad &\text{for} \quad p \notin B \text{ with } q \in B \text{ for all } q \nog p, \\
		\max_{q \nog p} \musup_q		\qquad &\text{for} \quad p \text{ non-critical}, \\
		\max_{q \nog p} \musup_q + 1	\qquad &\text{for} \quad p \text{ critical}.
	\end{cases}
\end{equation}
The second line is necessary because for any cell $p \in C\setminus B$ with $q \in B$ for all $q \nog p$ we have that $p$ is critical and all inputs come from inside $\Bsup$ so that $\musup_p=0$.

The main results in \Cref{thm:abff_maxnoncritsuper,thm:abff_maxnoncritsub,thm:abff_maxnoncritgen} below describe branches for the entire network. They combine \Cref{lem:rootgen} with an inductive investigation of cells outside of root subnetworks employing \Cref{lem:abff_maxnoncritnonmaxnoncrit,lem:abff_maxnoncritnonmaxcritL,lem:abff_maxnoncritnonmaxcritSN}. The number of critical cells along paths from the root subnetwork $\mu_p$ will be the induction parameter.

\begin{table}[h!]
	\begin{displaymath}
		\begin{array}{l|l}
			\Dsup_p																									& \text{Case} \\ \hline \\
			-\dfrac{\ell}{\sum_{\sigma \in \Sigma} a_\sigma}	& p \in \Bsup \\[15pt]
			- \dfrac{\sum_{\tau \notin \LL_p} a_\tau \Dsup_{\tau(p)} + \ell}{\sum_{\sigma \in \LL_{p}} a_\sigma}	& p \text{ non-critical}, \musup_p = 0 \\[15pt]
			- \dfrac{\sum_{\tau \colon \tau(p) \in \Qn_p} a_\tau \Dsup_{\tau(p)}}{\sum_{\sigma \in \LL_{p}} a_\sigma}	& p \text{ non-critical}, \musup_p > 0 \\[15pt]
			\dfrac{\ell}{\sum_{\sigma \in \Sigma} a_\sigma} \cdot \left( 1 + 2 \dfrac{\sum_{\sigma \in \LL_p, \tau \notin \LL_p} f_{\sigma\tau}}{\sum_{\sigma, \tau \in \LL_p} f_{\sigma\tau}} \right) - \dfrac{\sum_{\sigma \in \LL_p} f_{\sigma\lambda}}{\sum_{\sigma, \tau \in \LL_p} f_{\sigma\tau}}		& p \text{ critical}, \musup_p = 0 \\[15pt]
			\pm \sqrt{- \dfrac{\sum_{\tau \notin \LL_p} a_\tau \Dsup_{\tau(p)} + \ell}{\sum_{\sigma,\tau \in \LL_{p}} f_{\sigma\tau}}}	& p \text{ critical}, \musup_p=1 \\[15pt]
			\pm \sqrt{- \dfrac{\sum_{\tau \colon \tau(p) \in \Qn_p} a_\tau \Dsup_{\tau(p)}}{\sum_{\sigma, \tau \in \LL_p} f_{\sigma\tau}}}	& p \text{ critical}, \musup_p>1
		\end{array}
	\end{displaymath}
	\caption{Leading coefficients for the root subnetwork $\Bsup$ in supercritically branching steady states.}
	\label{tab:Dsup}
	
	\begin{displaymath}
		\begin{array}{l|l}
			\Dsub_p																									& \text{Case} \\ \hline \\
			\dfrac{\ell}{\sum_{\sigma \in \Sigma} a_\sigma}	& p \in \Bsub \\[15pt]
			- \dfrac{\sum_{\tau \notin \LL_p} a_\tau \Dsub_{\tau(p)} - \ell}{\sum_{\sigma \in \LL_{p}} a_\sigma}	& p \text{ non-critical}, \musub_p = 0 \\[15pt]
			- \dfrac{\sum_{\tau \colon \tau(p) \in \Qn_p} a_\tau \Dsub_{\tau(p)}}{\sum_{\sigma \in \LL_{p}} a_\sigma}			& p \text{ non-critical}, \musub_p > 0 \\[15pt]
			- \dfrac{\ell}{\sum_{\sigma \in \Sigma} a_\sigma} \cdot \left( 1 + 2 \dfrac{\sum_{\sigma \in \LL_p, \tau \notin \LL_p} f_{\sigma\tau}}{\sum_{\sigma, \tau \in \LL_p} f_{\sigma\tau}} \right) + \dfrac{\sum_{\sigma \in \LL_p} f_{\sigma\lambda}}{\sum_{\sigma, \tau \in \LL_p} f_{\sigma\tau}}		& p \text{ critical}, \musub_p = 0 \\[15pt]
			\pm \sqrt{ - \dfrac{\sum_{\tau \not\in \LL_p} a_\tau \Dsub_{\tau(p)} - \ell}{\sum_{\sigma, \tau \in \LL_p} f_{\sigma\tau}}}	& p \text{ critical}, \musub_p=1 \\[15pt]
			\pm \sqrt{- \dfrac{\sum_{\tau \colon \tau(p) \in \Qn_p} a_\tau \Dsub_{\tau(p)}}{\sum_{\sigma, \tau \in \LL_p} f_{\sigma\tau}}}	& p \text{ critical}, \musub_p>1
		\end{array}
	\end{displaymath}
	\caption{Leading coefficients for the root subnetwork $\Bsub$ in subcritically branching steady states.}
	\label{tab:Dsub}
\end{table}
\begin{theorem}[Supercritical branches]
	\label{thm:abff_maxnoncritsuper}
	Consider a feedforward network with cells $C$ and input maps $\Sigma$. Assume bifurcation assumption \emph{(B)} with non-maximal critical cells. For every root subnetwork $\Bsup \subset C$ for which we may define real non-vanishing coefficients as in \Cref{tab:Dsup} such that the inequalities
	\begin{equation}
		\label{eq:ineq_super}
		\begin{cases}
			\dfrac{\sum_{\tau \notin \LL_p} a_\tau \Dsup_{\tau(p)} + \ell}{\sum_{\sigma,\tau \in \LL_{p}} f_{\sigma\tau}} < 0, \quad &\text{for } p \text{ critical with } \musup_p = 1;\footnotemark \\
			\dfrac{\sum_{\tau \colon \tau(p) \in \Qn_p} a_\tau \Dsup_{\tau(p)}}{\sum_{\sigma, \tau \in \LL_p} f_{\sigma\tau}} <0, \quad
			&\text{for } p \text{ critical with } \musup_p > 1;
		\end{cases}
	\end{equation}
	\footnotetext{Note that this assumption implies that the coefficients in the first and fourth row do not coincide due to \Cref{lem:abff_maxnoncritDcritsuper}.}
	hold, there is a branch of steady states such that
	\begin{equation}
		\label{eq:maxnoncritsuper}
		x_p(\lambda) = \begin{cases}
			X(\lambda) = \Dsup_p \lambda + \OO \left( |\lambda| \right) \qquad &\text{for} \quad p \in \Bsup;\\
			\Dsup_p \cdot \lambda^{2^{-\musup_p}} + \OO \left( |\lambda|^{2^{-(\musup_p-1)}} \right) \qquad &\text{for} \quad p \notin \Bsup,
		\end{cases}
	\end{equation}
	for $\lambda >0$ small, where $\musup_p$ is the maximal number of critical cells $q \notin \Bsup$ along paths from any cell in $\Bsup$ to $p$ (see \eqref{eq:musup}) and
	\[ \Qn_p = \left\{ q \in\SSigma(p) \mid \musup_q = \max_{s \nog p} \musup_s \right\} \]
	(compare to \eqref{eq:abff_Qsup}). In particular, cells in $\Bsup$ are synchronous with $X(\lambda)$ as in \eqref{eq:abff_maxnoncritfullysynchX} while all cells not in $\Bsup$ are not synchronous to those in $\Bsup$.
\end{theorem}

\begin{remk}
	The assumptions on the coefficients $D_p$ in the theorem are made to guarantee the genericity conditions (L) and (SN), which are nonvanishing conditions on the nominator of $\Dsup_p$. By specifying the sign, we restrict to supercritical branches in this theorem. In \Cref{thm:abff_maxnoncritsub} below, we consider the corresponding situation for subcritically branching solutions. 
\end{remk}

\begin{proof}[Proof of \Cref{thm:abff_maxnoncritsuper}]
	We have seen in \Cref{cor:abff_subnetfullsynch} that the bifurcation equation \eqref{eq:taylor} for cells $p$ in an arbitrary subnetwork $B \subset C$ can be solved by the fully synchronous branch $x_p(\lambda) = X(\lambda)$ independent of the system parameters. Thus, for a root subnetwork $\Bsup$ it suffices to investigate cells $p \notin \Bsup$ which we will do with nested iterative arguments. The main induction is with respect to $\musup_p$ as indicated by the subheadings. The base case $\mu_p=0$ requires another inductive investigation with respect to the partial order $\nole$.
	
	\paragraph{Base case $\musup_p = 0$:}
	Cells $p\in\Bsup$, for which $\musup_p=0$ satisfy the statement of the theorem. Consider a cell $p \notin \Bsup$ with $\musup_p=0$. We have to distinguish two cases. If $p$ is critical, the definition of $\mu_p$ shows that necessarily $q \in B$ for all $q \nog p$. In particular, this is the situation of \Cref{lem:abff_maxnoncritnonmaxcritT} and we generically 
	obtain two branches of solutions for the bifurcation equation \eqref{eq:taylor} for $p$:
	\[ x_p(\lambda) = D_p^\pm \lambda + \OO\left(|\lambda|^2\right) \]
	for $|\lambda|$ small with coefficients
	\begin{align*}
		D_p^+ &= -\frac{\ell}{\sum_{\sigma \in \Sigma} a_\sigma} \\
		D_p^- &= \dfrac{\ell}{\sum_{\sigma \in \Sigma} a_\sigma} \cdot \left( 1 + 2 \dfrac{\sum_{\sigma \in \LL_p, \tau \notin \LL_p} f_{\sigma\tau}}{\sum_{\sigma, \tau \in \LL_p} f_{\sigma\tau}} \right) - \dfrac{\sum_{\sigma \in \LL_p} f_{\sigma\lambda}}{\sum_{\sigma, \tau \in \LL_p} f_{\sigma\tau}}
	\end{align*}
	as in \Cref{lem:abff_maxnoncritnonmaxcritT}. The assumption that these coefficients do not vanish implicitly implies that we are indeed in the generic situation. We consider only the second case
	\[ x_p(\lambda) = D_p^- \lambda + \OO\left(|\lambda|^2\right), \]
	i.e., $\Dsup_p=D_p^-$ (compare to the proof of \Cref{lem:rootgen}).
	
	We investigate the case that $p$ is non-critical inductively with respect to $\nole$. The argument above indicates that there must be a critical cell $q \notin \Bsup$ with $q \nog p$ such that $s \in \Bsup$ for all $s \nog q$ and no other critical cells along any path from $q$ to $p$. There may, however, be additional non-critical cells along paths from $q$ to $p$. First, assume that for all $q\in\SSigma(p)$ either $q \in \Bsup$ or $q$ is critical with $s\in\Bsup$ for all $s \nog q$ and $x_q(\lambda) = D_q \lambda + \OO\left(|\lambda|^2\right)$ as above. We define
	\[ \Dsup_p = - \dfrac{\sum_{\tau \notin \LL_p} a_\tau \Dsup_{\tau(p)} + \ell}{\sum_{\sigma \in \LL_{p}} a_\sigma}. \]
	By assumption $\Dsup_p \ne 0$, which implies
	\[ \sum_{\tau \notin \LL_p} a_\tau \Dsup_{\tau(p)} + \ell \ne 0. \]
	In combination with the fact that $x_q(\lambda)$ grows linearly in $\lambda$ for all $q \nog p$ this is condition (L) so that we are in the situation of \Cref{lem:abff_maxnoncritnonmaxnoncrit}. Hence, there is precisely one solution to the bifurcation equation \eqref{eq:taylor} for cell $p$
	\[ x_p(\lambda) = \Dsup_p \lambda + \OO\left(|\lambda|^2\right) \]
	for $|\lambda|$ small.
	
	Now consider an arbitrary non-critical cell $p$ with $\musup_p=0$. As an inductive hypothesis, we assume that for all non-critical cells $q\in\SSigma(p)$ with $\musup_q=0$ the bifurcation equation \eqref{eq:taylor} is uniquely solved by
	\[ x_q(\lambda) = \Dsup_q \lambda + \OO\left(|\lambda|^2\right) \]
	for $|\lambda|$ small with $\Dsup_q$ as in \Cref{tab:Dsup}. Note that with this assumption we have characterized the branch for all cells $q \in\SSigma(p)$, as these are either in $\Bsup$, critical with $s \in \Bsup$ for all $s \nog q$, or non-critical with $\musup_q=0$. As before, we define
	\[ \Dsup_p = - \dfrac{\sum_{\tau \notin \LL_p} a_\tau \Dsup_{\tau(p)} + \ell}{\sum_{\sigma \in \LL_{p}} a_\sigma}. \]
	By assumption $\Dsup_p \ne 0$, which implies (L). By \Cref{lem:abff_maxnoncritnonmaxnoncrit} there is precisely one solution to the bifurcation equation \eqref{eq:taylor} for cell $p$
	\[ x_p(\lambda) = \Dsup_p \lambda + \OO\left(|\lambda|^2\right) \]
	for $|\lambda|$ small. By induction with respect to $\nole$ this characterization of generically branching solutions to the bifurcation equation \eqref{eq:taylor} holds for all non-critical cells $p$ with $\musup_p=0$.
	
	\paragraph{Induction step $\musup_p > 0$:}
	We investigate cells $p \notin \Bsup$ with $\musup_p > 0$ inductively with respect to $\musup_p$ focusing on branching solutions for $\lambda>0$. Fix a cell $p \notin \Bsup$ with $\musup_p > 0$ and, as an inductive hypothesis, assume
	\begin{equation}
		\label{eq:abff_as-super-nonmaxbranch}
		x_q (\lambda) = \Dsup_q \lambda^{2^{-\musup_q}} + \OO \left( |\lambda|^{2^{-(\musup_q-1)}} \right)
	\end{equation}
	for $\lambda>0$ small with non-vanishing real coefficients $\Dsup_q$ as in \Cref{tab:Dsup} for all $q \in\SSigma(p)$ with $q \notin \Bsup$. For $q \in\SSigma(p)$ define
	\[ \Mu_q = \max_{s \nog q} \musup_s \]
	such that 
	\[ \Qn_q = \left\{ s \in\SSigma(q) \mid \musup_s = \Mu_q \right\}. \]
	Note that $\musup_q = \Mu_q$, if $q$ is non-critical, and $\musup_q = \Mu_q +1$, if $q$ is critical.
	We have to distinguish multiple cases.
	
	\subparagraph{Case $p$ non-critical, $\Mu_p=0$:}
	This situation has already been investigated in the base case, as $\Mu_p=\musup_p$ for $p$ non-critical.
	
	\subparagraph{Case $p$ non-critical, $\Mu_p>0$:}
	We define
	\[ \Dsup_p = - \frac{\sum_{\tau \colon \tau(p) \in \Qn_p} a_\tau \Dsup_{\tau(p)}}{\sum_{\sigma \in \LL_{p}} a_\sigma}. \]
	By assumption $\Dsup_p \ne 0$, which implies (SN). By \Cref{lem:abff_maxnoncritnonmaxnoncrit} there is precisely one solution to the bifurcation equation \eqref{eq:taylor} for cell $p$
	\[ x_p(\lambda) = \Dsup_p \lambda + \OO\left(|\lambda|^2\right) \]
	for $\lambda > 0$ small. Since $p$ is non-critical, we additionally have $\Mu_p=\musup_p$, and the coefficient matches the third row in \Cref{tab:Dsup}.
	
	\subparagraph{Case $p$ critical, $\Mu_p=0$:}
	The situation is as in \Cref{lem:abff_maxnoncritnonmaxcritL} (assumptions (\textbf{H}) and (L)). Since
	\[ \frac{\sum_{\tau \notin \LL_p} a_\tau \Dsup_{\tau(p)} + \ell}{\sum_{\sigma,\tau \in \LL_{p}} f_{\sigma\tau}} < 0 \]
	by assumption, we obtain two branching solutions to the bifurcation equation \eqref{eq:taylor} in cell $p$
	\[ x_p (\lambda) = \Dsup_p^\pm \lambda^{2^{-1}} + \OO (|\lambda|)  \]
	with
	\[ \Dsup_p^\pm	= \pm \sqrt{- \dfrac{\sum_{\tau \notin \LL_p} a_\tau \Dsup_{\tau(p)} + \ell}{\sum_{\sigma,\tau \in \LL_{p}} f_{\sigma\tau}}}, \]
	which we assumed to be real and nonzero. Since $p$ is critical, $\Mu_p=0$ implies $\musup_p=1$. Hence, the coefficient matches the fifth row in \Cref{tab:Dsup}.
	
	\subparagraph{Case $p$ critical, $\Mu_p>0$:}
	The situation is as in \Cref{lem:abff_maxnoncritnonmaxcritSN} (assumptions (\textbf{H}) and (SN)). Since
	\[ \frac{\sum_{\tau \colon \tau(p) \in \Qn_p} a_\tau \Dsup_{\tau(p)}}{\sum_{\sigma, \tau \in \LL_p} f_{\sigma\tau}} < 0 \]
	by assumption, we obtain two branching solutions to the bifurcation equation \eqref{eq:taylor} in cell $p$
	\[ x_p (\lambda) = \Dsup_p^\pm \cdot \lambda^{2^{-(\Mu_p+1)}} + \OO\left( |\lambda|^{2^{-\Mu_p}} \right), \]
	for $\lambda>0$ small, where
	\[ \Dsup_p^\pm = \pm \sqrt{- \frac{\sum_{\tau \colon \tau(p) \in \Qn_p} a_\tau \dsup_{\tau(p)}}{\sum_{\sigma, \tau \in \LL_p} f_{\sigma\tau}}} \]
	which we assumed to be real and nonzero. Since $p$ is critical, $\musup_p=\Mu_p+1$ and the coefficient matches the sixth row in \Cref{tab:Dsup}.
\end{proof}
\noindent
As a corollary we obtain an analogous result for subcritical branches.
\begin{theorem}[Subcritical branches]
	\label{thm:abff_maxnoncritsub}
	Consider a feedforward network with cells $C$ and input maps $\Sigma$. Assume bifurcation assumption \emph{(B)} with non-maximal critical cells. For every root subnetwork $\Bsub \subset C$ for which we may define real non-vanishing coefficients as in \Cref{tab:Dsub} such that the inequalities
	\begin{equation}
		\label{eq:ineq_sub}
		\begin{cases}
			\dfrac{\sum_{\tau \notin \LL_p} a_\tau \Dsub_{\tau(p)} - \ell}{\sum_{\sigma,\tau \in \LL_{p}} f_{\sigma\tau}} < 0, \quad &\text{for } p \text{ critical with } \musub_p = 1;\footnotemark \\
			\dfrac{\sum_{\tau \colon \tau(p) \in \Qn_p} a_\tau \Dsub_{\tau(p)}}{\sum_{\sigma, \tau \in \LL_p} f_{\sigma\tau}} <0, \quad
			&\text{for } p \text{ critical with } \musub_p > 1;
		\end{cases}
	\end{equation}
	\footnotetext{Note that this assumption implies that the coefficients in the first and fourth row do not coincide due to \Cref{lem:abff_maxnoncritDcritsuper}.}
	hold, there is a branch of steady states such that
	\begin{equation}
		\label{eq:maxnoncritsub}
		x_p(\lambda) = \begin{cases}
			X(\lambda) = \Dsub_p \lambda + \OO \left( |\lambda| \right) \qquad &\text{for} \quad p \in \Bsub;\\
			\Dsub_p \cdot (-\lambda)^{2^{-\musub_p}} + \OO \left( |\lambda|^{2^{-(\musub_p-1)}} \right), \qquad &\text{for} \quad p \notin \Bsub,
		\end{cases}
	\end{equation}
	for $\lambda <0$ small, where $\musub_p$ is the maximal number of critical cells $q \notin \Bsub$ along paths from any cell in $\Bsub$ to $p$ (see \eqref{eq:musup}) and
	\[ \Qn_p = \left\{ q \in\SSigma(p) \mid \musub_q = \max_{s \nog p} \musub_s \right\} \]
	(compare to \eqref{eq:abff_Qsup}). In particular, cells in $\Bsub$ are synchronous with $X(\lambda)$ as in \eqref{eq:abff_maxnoncritfullysynchX} while all cells not in $\Bsub$ are not synchronous to those in $\Bsub$.
\end{theorem}
\begin{proof}
	Substituting the parameter $\mu = - \lambda$ and the system parameters $\iota =- \ell, \theta_{\sigma\lambda} = - f_{\sigma\lambda}$, the governing function Taylor expands as
	\[ f (x_{\sigma_1(p)}, \dotsc, x_{\sigma_n(p)}, \lambda) = \sum_{\sigma \in \Sigma} a_\sigma x_{\sigma(p)} + \iota \mu + \sum_{\sigma, \tau \in \Sigma} f_{\sigma\tau} x_{\sigma(p)}x_{\tau(p)} + \sum_{\sigma \in \Sigma} \theta_{\sigma\lambda} \mu x_{\sigma(p)} + f_{\lambda\lambda} \mu^2 + \text{h.o.t.}, \]
	as in \eqref{eq:taylor}. The result follows immediately from \Cref{thm:abff_maxnoncritsuper} for $\mu > 0$.
\end{proof}

Finally, we complete our considerations by observing that all occurring branches of steady states are as described in \Cref{thm:abff_maxnoncritsuper,thm:abff_maxnoncritsub}, as long as the non-degenericity conditions (L) and (SN) can be satisfied by the leading coefficients.

\begin{theorem}
	\label{thm:abff_maxnoncritgen}
	Consider a feedforward network with cells $C$ and input maps $\Sigma$. Assume bifurcation assumption \emph{(B)} with non-maximal critical cells and generic system parameters. 
	Every occurring branch of steady states that is not the fully synchronous continuation (see \Cref{lem:abff_maxnoncritfullysynch}) is either as in \Cref{thm:abff_maxnoncritsuper} or as in \Cref{thm:abff_maxnoncritsub}, if leading coefficients can be chosen as in \Cref{tab:Dsup,tab:Dsub} such that non-degenericity conditions (L) and (SN) are satisfied. In particular, there is a root subnetwork $\Bsup \subset C$ such that the cells outside of $\Bsup$ branch super- or subcritically with asymptotics determined by the number of critical cells along paths into these cells.
\end{theorem}
\begin{proof}
	The fully synchronous continuation of the bifurcation point exists for all values of system parameters. \Cref{lem:rootgen} shows that for all other generically existing branches there is a root subnetwork $\Bsup \subset C$ such that $x_p (\lambda) = X(\lambda)$ for all $p \in \Bsup$. Fix a root subnetwork $B\subset C$ and investigate the remaining cells. By assumption, (L) and (SN) hold, which implies that the conditions of \Cref{lem:abff_maxnoncritnonmaxcritL,lem:abff_maxnoncritnonmaxcritT,lem:abff_maxnoncritnonmaxcritSN} are satisfied and branching solutions to the bifurcation equation \eqref{eq:taylor} are as in these lemmas (for subcritically branching solutions this requires a substitution as in the proof of \Cref{thm:abff_maxnoncritsub}).
	
	Note that for a critical cell $p \notin B$ with $\mu_p=1$ inequality \eqref{eq:ineq_super} is satisfied if and only if \eqref{eq:ineq_sub} is not and vice versa. In particular, these cells indicate either a super- or a subcritical branch of steady states as in \Cref{thm:abff_maxnoncritsuper,thm:abff_maxnoncritsub}. If two critical cells $p,p' \notin B$ with $\mu_p=\mu_{p'}=1$ indicate a branch of steady states for opposing signs of $\lambda$ there is no branch generated by $B$.
	
	If the direction of branching is determined uniquely, for a critical cell $p \notin B$ with $\mu_p>1$, there is a choice in signs for the coefficients $D_q$ with $q\in \Qn_p$. Due to (SN), exactly half of these signs satisfy the inequality \eqref{eq:ineq_super} or \eqref{eq:ineq_sub} respectively. If two critical cells $p,p'\notin B$ with $\mu_p=\mu_{p'}>1$ require opposing signs for a cell $q \in \Qn_p \cap \Qn_{p'}$ the branch fails to exist and $B$ does not generate a branching solution. Otherwise, the solutions branch as in \Cref{thm:abff_maxnoncritsuper,thm:abff_maxnoncritsub}.
\end{proof}

\begin{remk}
	Throughout this section we have used the term `generic' on multiple occasions. It is used to indicate that the corresponding statement holds true for an open and dense subset of the system parameters, i.e., the low-order partial derivatives of the governing function $f$. In particular, where used it guarantees that we do not divide by $0$ or that leading order terms in the equations we investigate do not vanish.
	
	However, proving genericity of (L) and (SN) is not possible due to their inductive nature. We cannot thoroughly exclude the possibility that the network structure forces one of the corresponding weighted sums of leading coefficients to vanish identically, an issue that we have never encountered. It would lead to vanishing leading order coefficients in \Cref{tab:Dsup,tab:Dsub} and therefore to lower square root orders. When applying our results to a specific network, one computes the leading coefficients algorithmically anyway (see \Cref{rem:c2-algorithm} below). In this process the conditions (L) and (SN) are checked.
	
	Finally, the following reasoning underlines, why we do not expect (L) or (SN) not to hold generically. For any cell $p \notin B$ that is either non-critical or critical with $\musup_p =1$ there exists a cell $q \nog p$ with $q$ critical, $q \notin B$, and $s \in B$ for all $s \nog q$. Hence, there are cells $q \nog p$ for which $\Dsup_q$ is as in the fourth row of \Cref{tab:Dsup} (\Cref{tab:Dsub}) as well as cells $q'\nog p$ for which $\Dsup_{q'}$ is as in the first row of \Cref{tab:Dsup} (\Cref{tab:Dsub}). In particular, some of the terms in the weighted sum
	\[ \sum_{\tau \not\in \LL_p} a_\tau \Dsup_{\tau(p)} \]
	depend only on the $a_\sigma$ and $\ell$, while others also depend on the $f_{\sigma\tau}$. We do not expect this weighted sum to equal $-\ell$ identically, which would be necessary for (L) to be violated.
	
	A similar argument for (SN) can be made. For $p \notin B$ non-critical with $\musup_p>0$ or for $p \notin B$ critical with $\musup_p>1$ there exists a cell $q \nog p$ with $\musup_q \ge 1$. For such a cell the coefficient $\Dsup_q$ consists of square roots of terms depending on the $a_\sigma$. By definition this holds true for all $q \in\Qn_p$. In particular,
	\[ \sum_{\tau \colon \tau(p)\in \Qn_p} a_\tau \Dsup_{\tau(p)} \]
	is a weighted sum of these square root terms with weights given by the $a_\tau$. Again, we do not expect this weighted sum to vanish identically.
\end{remk}
\begin{remk}
	The branching solutions in \Cref{thm:abff_maxnoncritsuper,thm:abff_maxnoncritsub,thm:abff_maxnoncritgen} contain those that are described in Section 6 in \textcite{Soares.2018} for layered feedforward networks as a special case. Therein, generically all non-maximal cells are critical if the maximal cells are non-critical. This generalization is due to the fact that the class of feedforward networks satisfying the equivalent definitions in \Cref{prop:equivalence} contains layered feedforward networks (compare to \Cref{rem:layered}).
	\hspace*{\fill}$\triangle$
\end{remk}
\begin{remk}
	\label{rem:c2-algorithm}
	\Cref{thm:abff_maxnoncritsuper,thm:abff_maxnoncritsub,thm:abff_maxnoncritgen} provide a constructive method to determine all possible branching solutions for the bifurcation problem (B). The first step is to determine all root subnetworks $B \subset C$. A solution branch is computed as $x_p$ staying in the fully synchronous state for all $p \in B$ and the states of the remaining cells being determined by the number of critical cells in between $p$ and $B$ where the leading coefficients $\Dsup_p$ are chosen according to \Cref{tab:Dsup,tab:Dsub}. In each cell, we have to check whether \eqref{eq:ineq_super} or \eqref{eq:ineq_sub} can be satisfied, i.e. whether the root subnetwork generates a branch. Critical cells outside of $B$ determine existence and direction of branches of steady states as in the proof of \Cref{thm:abff_maxnoncritgen}. The asymptotic order can be determined inductively via \eqref{eq:muiteratively}.
	\hspace*{\fill}$\triangle$
\end{remk}
\begin{remk}
	\label{rem:c2-systemparam}
	Note that the conditions determining the existence of branching solutions in \Cref{thm:abff_maxnoncritsuper,thm:abff_maxnoncritsub} depend only on the system parameters -- the partial derivatives of $f$. This implies the existence of different branches in different regions of system parameter space which may also vary according to the direction of branching -- super- or subcritical.
	\hspace*{\fill}$\triangle$
\end{remk}

%% file: II_genfeedforward_ex.tex
\section{An example}
\label{sec:ex}
We illustrate the analytic results from \Cref{thm:maxcrit,thm:abff_maxnoncritsuper,thm:abff_maxnoncritsub,thm:abff_maxnoncritgen} in the network in \Cref{fig:numerics1} that was numerically investigated in the introduction. In the case of non-maximal critical cells we employ the algorithm presented in \Cref{rem:c2-algorithm} and highlight the peculiarities mentioned in \Cref{rem:c2-systemparam}. Consider the feedforward network given by the graph in \Cref{fig:numerics1}. 
Each arrow color corresponds to one input map $\sigma \colon C \to C$. Note that we have not drawn an arrow for $\sigma_1 = \Id$ corresponding to the internal dynamics which we implicitly assume to be there. This network is clearly a feedforward network as it does not contain any cycles besides self-loops (note that it is not a layered feedforward network as in \cite{Soares.2018}). Its only maximal cell is cell $5$. Furthermore, it possesses two different loop-types $\LL_5 = \Sigma$ and $\LL_1 = \LL_2 = \LL_3 = \LL_4 =\{\Id\}$, see \eqref{eq:deflooptype}. As a matter of fact, this network is a fundamental network in the language of \textcite{Rink.2014} but this is not important for the upcoming investigations. Assuming a one-dimensional internal phase space $x_i \in V = \RR$ and additional dependence on a real parameter $\lambda \in \RR$, the corresponding dynamics is governed by
\input{5cellfeedforwardwithoddbifurcationsvfextnet}as in \eqref{eq:numerics1}. We want to investigate bifurcations of steady states as in the bifurcation scenario (B) described in the beginning of \Cref{sec:bi}. That is, we assume
\[ \gamma_f (0,0) = 0. \]
The linearization at this steady state is
\begin{equation*}
	D_x \gamma_f (0,0) = \begin{pmatrix*}[r]
		a_{\Id}	& \textcolor{red}{a_{\sigma_2}}	& \textcolor{blue}{a_{\sigma_3}}	& \textcolor{grey}{a_{\sigma_4}}	& \textcolor{magenta}{a_{\sigma_5}}	\\
		0		& a_{\Id}	& 0	&	\textcolor{blue}{a_{\sigma_3}}	& \textcolor{red}{a_{\sigma_2}} + \textcolor{grey}{a_{\sigma_4}} + \textcolor{magenta}{a_{\sigma_5}}	\\
		0		& 0	& a_{\Id}	& \textcolor{red}{a_{\sigma_2}}	& \textcolor{blue}{a_{\sigma_3}} + \textcolor{grey}{a_{\sigma_4}} + \textcolor{magenta}{a_{\sigma_5}}	\\
		0		& 0	& 0	& a_{\Id}	& \textcolor{red}{a_{\sigma_2}} + \textcolor{blue}{a_{\sigma_3}} + \textcolor{grey}{a_{\sigma_4}} + \textcolor{magenta}{a_{\sigma_5}} \\
		0		& 0	& 0	& 0	& a_{\Id} + \textcolor{red}{a_{\sigma_2}} + \textcolor{blue}{a_{\sigma_3}} + \textcolor{grey}{a_{\sigma_4}} + \textcolor{magenta}{a_{\sigma_5}}
	\end{pmatrix*}.
\end{equation*}
Herein we define $a_\sigma = \partial_\sigma f(0,0)$ (see \eqref{eq:lin}). Furthermore, we follow our convention $\sigma_1 = \Id$. The other partial derivatives are abbreviated accordingly again:
\begin{alignat*}{4}
	& f_{\sigma\tau}		&&= \frac{1}{2} \partial_{\sigma\tau} f(0,0),	\qquad	&& f_{\sigma\lambda}	&&= \partial_{\sigma\lambda} f(0,0), \\
	& \ell					&&= \partial_\lambda f(0,0),	\qquad					&& f_{\lambda\lambda}	&&= \frac{1}{2} \partial_{\lambda\lambda} f(0,0).
\end{alignat*}

The eigenvalues of the linearization are in one-to-one correspondence with the loop-types of the network. This can also easily be read off of the matrix. As a matter of fact the linearization has two eigenvalues $\sum_{\sigma \in \LL_5}a_\sigma\! =\! a_{\Id} + a_{\sigma_2}\! + a_{\sigma_3}\! + a_{\sigma_4}\! + a_{\sigma_5}$, which is simple, and $\sum_{\sigma \in \LL_1}a_\sigma = a_{\Id}$ which has algebraic multiplicity $4$. For a steady state bifurcation to occur, the linearization has to have an eigenvalue $0$. Generically -- i.e. for a generic choice of system parameters --, in such a point only one of the two eigenvalues vanishes. Under these assumptions we investigate generic solutions to 
\[ \gamma_f (x,\lambda) = 0 \]
close to the bifurcation point.

\paragraph{Case \RN{1}}
Let us investigate the case $\sum_{\sigma \in \LL_5}a_\sigma\! =\! a_{\Id} + a_{\sigma_2}\! + a_{\sigma_3}\! + a_{\sigma_4}\! + a_{\sigma_5}\! = 0$ and $\sum_{\sigma \in \LL_1}a_\sigma = a_{\Id} \ne 0$ first. In particular, this means the maximal cell is critical while all other cells are not. As a result, all branches of steady state solutions are given in \Cref{thm:maxcrit}. As there is only one maximal cell, necessarily all branches are fully synchronous. We obtain two different saddle node branches depending on the system parameters. If
\[ \frac{\ell}{\sum_{\sigma, \tau \in \Sigma} f_{\sigma\tau}} < 0 \]
we compute
\[ x_i(\lambda) = \pm \sqrt{- \frac{\ell}{\sum_{\sigma, \tau \in \Sigma} f_{\sigma\tau}}} \cdot \sqrt{\lambda} + \OO (|\lambda|). \]
for $i=1,\dotsc,5$. Note that therein the choice of sign is the same for all cells simultaneously yielding exactly two fully synchronous branches. On the other hand, if
\[ \frac{\ell}{\sum_{\sigma, \tau \in \Sigma} f_{\sigma\tau}} > 0 \]
we obtain
\[ x_i(\lambda) = \pm \sqrt{\frac{\ell}{\sum_{\sigma, \tau \in \Sigma} f_{\sigma\tau}}} \cdot \sqrt{-\lambda} + \OO (|\lambda|). \]
for $i=1,\dotsc,5$ accordingly. These branches exists for $|\lambda|$ small. Generically, no other cases are possible so that no other branching solutions exist.

\paragraph{Case \RN{2}}
Next, we turn to the case $K\! = \sum_{\sigma \in \LL_5}a_\sigma\! =\! a_{\Id} + a_{\sigma_2}\! + a_{\sigma_3}\! + a_{\sigma_4}\! + a_{\sigma_5}\! \ne 0$ and $\sum_{\sigma \in \LL_1}a_\sigma = a_{\Id} = 0$. Equivalently, the maximal cell is not critical but all the other cells are. \Cref{thm:abff_maxnoncritsuper,thm:abff_maxnoncritsub} provide all generic branching solutions and we employ the algorithm in \Cref{rem:c2-algorithm} to characterize them. As a first step, we have to determine all possible root subnetworks, i.e. subnetworks $B \subset C$ such that $p \notin \Bsup$ but $q \in \Bsup$ for all $q \nog p$ implies $p$ critical. The possible choices are $\{ 5 \}, \{ 4,5 \}, \{ 3,4,5 \}, \{ 2,4,5 \}, \{ 2,3,4,5 \}, \{ 1,2,3,4,5 \}$. For each solution branch the cells in exactly one of these subnetworks are in the fully synchronous state $X(\lambda)$
as in \Cref{lem:abff_maxnoncritfullysynch} while all others are not. The solutions for the variables of the remaining cells are computed iteratively with respect to the partial order $\nole$ according to the rules in \Cref{thm:abff_maxnoncritsuper,thm:abff_maxnoncritsub}. If for a cell $p$ the necessary inequalities \eqref{eq:ineq_super} or \eqref{eq:ineq_sub} are not satisfied for all choices of coefficients for cells $q \nog p$, this implies that the solution branch does not exist. We describe the branches as briefly as possible starting with the most simple case.

\subparagraph{Case \RN{2}.{\romannumeral 1})}
Assume $\Bsup = \{ 1,2,3,4,5 \}$. All cells remain in the fully synchronous state, i.e. 
\[ x_i (\lambda) = X(\lambda) = - \frac{\ell}{K} \lambda + \OO \left( |\lambda|^2 \right) \]
for $i=1,\dotsc,5$. This branch exists for $|\lambda|$ small independent of the sign and without any further restrictions on the system parameters.

\subparagraph{Case \RN{2}.{\romannumeral 2})}
Next, assume $\Bsup = \{ 2,3,4,5 \}$. We obtain
\[ x_5 (\lambda) = x_4 (\lambda) = x_3 (\lambda) = x_2 (\lambda) = X(\lambda). \]
As cell $1$ is critical but not in the fully synchronous state, this leaves
\[ x_1 (\lambda) = \left( \frac{\ell}{K} \left(1+ 2 \frac{f_{\Id\sigma_2} + f_{\Id\sigma_3} + f_{\Id\sigma_4} + f_{\Id\sigma_5}}{f_{\Id\Id}}\right) - \frac{f_{\Id\lambda}}{f_{\Id\Id}} \right) \lambda + \OO \left(|\lambda|^2\right). \]
This branch exists without any further restrictions on the system parameters as well.

\subparagraph{Case \RN{2}.{\romannumeral 3})}
For $\Bsup = \{ 3,4,5 \}$ we obtain
\[ x_5 (\lambda) = x_4 (\lambda) = x_3 (\lambda) = X(\lambda) = - \frac{\ell}{K} \lambda + \OO \left( |\lambda|^2 \right) \]
and abbreviate $\Dsup_5 = \Dsup_4 = \Dsup_3 = - \ell / K$. As cell $2$ is critical but $2 \notin \Bsup$, we obtain
\begin{align*}
	x_2 (\lambda)	&= \left( \frac{\ell}{K} \left(1+ 2 \frac{f_{\Id\sigma_2} + f_{\Id\sigma_3} + f_{\Id\sigma_4} + f_{\Id\sigma_5}}{f_{\Id\Id}}\right) - \frac{f_{\Id\lambda}}{f_{\Id\Id}} \right) \lambda + \OO \left(|\lambda|^2\right) \\
	&= \Dsup_2 \lambda + \OO \left( |\lambda|^2 \right).
\end{align*}
Then cell $1$ receives an input from a cell not in $\Bsup$. Thus, we have to distinguish two cases according to \eqref{eq:ineq_super}. If
\[ (*) = \frac{a_{\sigma_2} \Dsup_2 + a_{\sigma_3} \Dsup_3 + a_{\sigma_4} \Dsup_4 + a_{\sigma_5} \Dsup_5 + \ell}{f_{\Id\Id}} <0, \]
we obtain
\[ x_1 (\lambda) = \pm \sqrt{- \frac{a_{\sigma_2} \Dsup_2 + a_{\sigma_3} \Dsup_3 + a_{\sigma_4} \Dsup_4 + a_{\sigma_5} \Dsup_5 + \ell}{f_{\Id\Id}}} \cdot \sqrt{\lambda} + \OO(|\lambda|) \]
for small $\lambda > 0$. If, on the other hand, $(*)>0$, there is no supercritical solution branch with $\Bsup = \{ 3,4,5 \}$. On the other hand, the solutions for cells $2,3,4,5$ can be written as $x_i (\lambda) = E_i (-\lambda) + \OO (|\lambda|^2)$, where $E_i = -\Dsup_i$. The condition \eqref{eq:ineq_sub} for the existence of a subcritically branching solution for cell $1$ is
\[ \frac{a_{\sigma_2} E_2 + a_{\sigma_3} E_3 + a_{\sigma_4} E_4 + a_{\sigma_5} E_5 - \ell}{f_{\Id\Id}} <0. \]
Note that the left hand side of this inequality is $-(*)$. Hence, if $(*) >0$, we obtain
\[ x_1 (\lambda) = \pm \sqrt{\frac{a_{\sigma_2} \Dsup_2 + a_{\sigma_3} \Dsup_3 + a_{\sigma_4} \Dsup_4 + a_{\sigma_5} \Dsup_5 + \ell}{f_{\Id\Id}}} \cdot \sqrt{- \lambda} + \OO(|\lambda|). \]
If $(*)<0$ there is no solution for cell $1$ for $\lambda<0$. Summarizing we see that depending on the sign of $(*)$, the branch exists for precisely one sign of $\lambda$ -- this includes the solutions for cells $2,3,4,5$.

\subparagraph{Case \RN{2}.{\romannumeral 4})}
The considerations for $\Bsup = \{ 2,4,5 \}$ are almost identical to those made for $\Bsup = \{ 3,4,5 \}$. Exchanging cells $2$ and $3$ as well as the input maps $\sigma_2$ and $\sigma_3$ provides the solution branches.

\subparagraph{Case \RN{2}.{\romannumeral 5})}
The case $\Bsup = \{ 4,5 \}$ is very similar as well. Cells $4$ and $5$ remain in the fully synchronous state $X(\lambda)$. More precisely for cells $i = 2, \dotsc, 5$ we obtain
\[ x_i (\lambda) = \Dsup_i \lambda + \OO \left( |\lambda|^2 \right), \quad x_i (\lambda) = E_i \cdot (-\lambda) + \OO \left( |\lambda|^2 \right) \]
with 
\begin{align*}
	\Dsup_5 = -E_5 = \Dsup_4 = - E_4 &= - \frac{\ell}{K}, \\
	\Dsup_3 = -E_3 = \Dsup_2 = - E_2 &= \frac{\ell}{K} \left(1+ 2 \frac{f_{\Id\sigma_2} + f_{\Id\sigma_3} + f_{\Id\sigma_4} + f_{\Id\sigma_5}}{f_{\Id\Id}}\right) - \frac{f_{\Id\lambda}}{f_{\Id\Id}}
\end{align*}
for $|\lambda|$ small respectively. Similar to before we obtain
\begin{align*}
	x_1 (\lambda) &= \pm \sqrt{- \frac{a_{\sigma_2} \Dsup_2 + a_{\sigma_3} \Dsup_3 + a_{\sigma_4} \Dsup_4 + a_{\sigma_5} \Dsup_5 + \ell}{f_{\Id\Id}}} \cdot \sqrt{\lambda} + \OO(|\lambda|) \quad \text{or} \\
	x_1 (\lambda) &= \pm \sqrt{\frac{a_{\sigma_2} \Dsup_2 + a_{\sigma_3} \Dsup_3 + a_{\sigma_4} \Dsup_4 + a_{\sigma_5} \Dsup_5 + \ell}{f_{\Id\Id}}} \cdot \sqrt{- \lambda} + \OO(|\lambda|)
\end{align*}
for $\lambda >0$ or $\lambda<0$ respectively, if
\[ \frac{a_{\sigma_2} \Dsup_2 + a_{\sigma_3} \Dsup_3 + a_{\sigma_4} \Dsup_4 + a_{\sigma_5} \Dsup_5 + \ell}{f_{\Id\Id}} <0 \quad \text{or} \quad >0. \]
Once again, the solution branches in this case exist for precisely one sign of $\lambda$.

\subparagraph{Case \RN{2}.{\romannumeral 6})}
Finally, we investigate the case $\Bsup = \{ 5 \}$. The mechanism that relates the two cases -- i.e. super- or subcritically branching solutions -- is the same as in the previous cases. Therefore we omit the computational details. Cell $5$ remains in the fully synchronous state
\[ x_5 (\lambda) = X(\lambda) = \Dsup_5 \lambda + \OO \left( |\lambda|^2 \right). \]
For cell $4$ we obtain
\begin{align*}
	x_4 (\lambda)	&= \Dsup_4 \lambda + \OO \left( |\lambda|^2 \right) \\
	&= \left( \frac{\ell}{K} \left(1+ 2 \frac{f_{\Id\sigma_2} + f_{\Id\sigma_3} + f_{\Id\sigma_4} + f_{\Id\sigma_5}}{f_{\Id\Id}}\right) - \frac{f_{\Id\lambda}}{f_{\Id\Id}} \right) \lambda + \OO \left(|\lambda|^2\right).
\end{align*}
Considering cell $3$, we obtain
\begin{align*}
	x_3 (\lambda) &= \pm \sqrt{- \frac{a_{\sigma_2} \Dsup_4 + (a_{\sigma_3} + a_{\sigma_4} + a_{\sigma_5}) \Dsup_5 + \ell}{f_{\Id\Id}}} \cdot \sqrt{\lambda} + \OO(|\lambda|) \quad \text{or} \\
	x_3 (\lambda) &= \pm \sqrt{\frac{a_{\sigma_2} \Dsup_4 + (a_{\sigma_3} + a_{\sigma_4} + a_{\sigma_5}) \Dsup_5 + \ell}{f_{\Id\Id}}} \cdot \sqrt{- \lambda} + \OO(|\lambda|),
\end{align*}
if
\[ (*) = \frac{a_{\sigma_2} \Dsup_4 + (a_{\sigma_3} + a_{\sigma_4} + a_{\sigma_5}) \Dsup_5 + \ell}{f_{\Id\Id}} <0 \quad \text{or} \quad >0 \]
respectively. In particular, $\Bsub = \{ 5 \}$ does not provide a solution branch for $\lambda>0$, if $(*)>0$, or for $\lambda<0$, if $(*)<0$. Similarly, we obtain
\begin{align*}
	x_2 (\lambda) &= \pm \sqrt{- \frac{a_{\sigma_3} \Dsup_4 + (a_{\sigma_2} + a_{\sigma_4} + a_{\sigma_5}) \Dsup_5 + \ell}{f_{\Id\Id}}} \cdot \sqrt{\lambda} + \OO(|\lambda|) \quad \text{or} \\
	x_2 (\lambda) &= \pm \sqrt{\frac{a_{\sigma_3} \Dsup_4 + (a_{\sigma_2} + a_{\sigma_4} + a_{\sigma_5}) \Dsup_5 + \ell}{f_{\Id\Id}}} \cdot \sqrt{- \lambda} + \OO(|\lambda|),
\end{align*}
if
\[ (**) = \frac{a_{\sigma_3} \Dsup_4 + (a_{\sigma_2} + a_{\sigma_4} + a_{\sigma_5}) \Dsup_5 + \ell}{f_{\Id\Id}} <0 \quad \text{or} \quad >0 \]
respectively. In particular, $\Bsub = \{ 5 \}$ does not provide a solution branch for $\lambda>0$, if $(**)>0$, or for $\lambda<0$, if $(**)<0$. Hence, if $(*)$ and $(**)$ have opposite signs, neither of the two branches exists. If both have the same sign, we abbreviate the coefficients as $\pm\Dsup_3, \pm\Dsup_2, \pm E_3, \pm E_2$. We only need to investigate cell $1$ in that case. Consider $(*), (**)<0$. If
\[ (***) = \frac{\pm a_{\sigma_2} \Dsup_2 \pm a_{\sigma_3} \Dsup_3}{f_{\Id\Id}} <0, \]
we obtain
\[ x_1 (\lambda) = \pm \sqrt{- \frac{\pm a_{\sigma_2} \Dsup_2 \pm a_{\sigma_3} \Dsup_3}{f_{\Id\Id}}} \sqrt{\sqrt{\lambda}} + \OO \left( \sqrt{|\lambda|} \right). \]
If $(***)>0$, the solution branch does not exist. Note that $(***)$ depends on the choice of signs for the coefficients in cells $2$ and $3$. Therefore, half of the possible choices yields a negative sign of $(***)$ while the other half yields a positive sign. This is due to the fact that $(***)$ and $-(***)$ are both possible choices, while $(***)\ne0$ generically. Similarly, for $(*), (**) >0$ we obtain
\[ x_1 (\lambda) = \pm \sqrt{- \frac{\pm a_{\sigma_2} \Dsub_2 \pm a_{\sigma_3} \Dsub_3}{f_{\Id\Id}}} \sqrt{\sqrt{-\lambda}} + \OO \left( \sqrt{|\lambda|} \right), \]
if
\[ \frac{\pm a_{\sigma_2} \Dsub_2 \pm a_{\sigma_3} \Dsub_3}{f_{\Id\Id}} <0 \]
for admissible choices of signs.

We have therefore computed all generic branches of steady states. We summarize the results in \Cref{tab:ex_branches}.
\begin{table}[h!]
	\resizebox{\textwidth}{!}{%
		\renewcommand{\arraystretch}{2}
		\begin{tabular}{c|c|c|c}
			\textbf{Root subnetwork}	& \textbf{Existence condition}	& \textbf{Asymptotics}	& \textbf{Type} \\ \hline \hline
			$\{1,2,3,4,5\}$	& --	& $\left( \lambda, \lambda, \lambda, \lambda, \lambda \right)$	& continuation \\ \hline
			$\{2,3,4,5\}$	& --	& $\left( \lambda, \lambda, \lambda, \lambda, \lambda \right)$	& transcritical \\ \hline
			$\{3,4,5\}$		& $\lambda \cdot \dfrac{a_{\sigma_2} \Dsup_2 + a_{\sigma_3} \Dsup_3 + a_{\sigma_4} \Dsup_4 + a_{\sigma_5} \Dsup_5 + \ell}{f_{\Id\Id}} <0 $	& $\left( \sqrt{|\lambda|}, \lambda, \lambda, \lambda, \lambda \right)$	& $4\times$ saddle node \\ \hline
			$\{2,4,5\}$ & $\lambda \cdot \dfrac{a_{\sigma_2} \Dsup_2 + a_{\sigma_3} \Dsup_3 + a_{\sigma_4} \Dsup_4 + a_{\sigma_5} \Dsup_5 + \ell}{f_{\Id\Id}} <0 $	& $\left( \sqrt{|\lambda|}, \lambda, \lambda, \lambda, \lambda \right)$	& $4\times$ saddle node \\ \hline
			$\{4,5\}$ & $\lambda \cdot \dfrac{a_{\sigma_2} \Dsup_2 + a_{\sigma_3} \Dsup_3 + a_{\sigma_4} \Dsup_4 + a_{\sigma_5} \Dsup_5 + \ell}{f_{\Id\Id}} <0 $	& $\left( \sqrt{|\lambda|}, \lambda, \lambda, \lambda, \lambda \right)$	& $4\times$ saddle node \\ \hline
			$\{5\}$	&
			$\begin{aligned}
				\lambda\cdot\dfrac{a_{\sigma_2} \Dsup_4 + (a_{\sigma_3} + a_{\sigma_4} + a_{\sigma_5}) \Dsup_5 + \ell}{f_{\Id\Id}} &<0 \\
				\lambda\cdot\dfrac{a_{\sigma_3} \Dsup_4 + (a_{\sigma_2} + a_{\sigma_4} + a_{\sigma_5}) \Dsup_5 + \ell}{f_{\Id\Id}} &<0 \\
			\end{aligned} $
			&  $\left( \sqrt[4]{|\lambda|}, \sqrt{|\lambda|}, \sqrt{|\lambda|}, \lambda, \lambda \right)$ & $4\times$ amplified saddle node
		\end{tabular}
	}
	\caption{
		Summary of the generic bifurcation branches in Case \RN{2}. As branches are computed per root subnetwork, these are indicated in the first column. The entries in the second column show conditions that have to be satisfied in order for the branches with asymptotics as in the third column to exist. These conditions depend on the leading coefficients in certain cells for a specific branch and on the sign of the bifurcation parameter $\lambda$. Note that the coefficients are specific to the root subnetwork, i.e., the $D_i$ vary along the rows of the table. Regions in parameter space for which multiple conditions are satisfied allow for all of the corresponding branches. Note that the only cells within the root subnetworks follow the continuation of the fully synchronous solution. This is not displayed in the asymptotics in the third column, as cells outside of the root subnetwork may have the same leading order, but branch transcritically.
	}
	\label{tab:ex_branches}
\end{table}
We see that there are numerous ways in which a solution branch for the root subnetwork $\Bsup = \{ 5 \}$ fails to exist. These ultimately depend on the system parameters. Hence, there are different solutions in different regions of system parameter space. We briefly introduce two cases to illustrate that already this simple network produces unexpected -- compared to the summary of the amplification effect -- bifurcation scenarios.

Consider the bifurcation scenario as before with $a_{\sigma_3} = - 2 a_{\sigma_2}$ as well as $a_{\sigma_4}=a_{\sigma_5}=0$ and investigate $\Bsup = \{5\}$. We compute $(**) = - 2 \cdot (*)$ proving that generically $(*)$ and $(**)$ have opposite signs. Therefore, there are no branching solutions with $\Bsup = \{5\}$, as cell $2$ forces the branch to exist for $\lambda >0$ and cell $3$ forces it to exist for $\lambda<0$ or the other way around. This implies the existence of an open region in parameter space for which this issue occurs. The reason lies in the structure of the network. The two cells $2$ and $3$ receive the same inputs. However, the input from cell $4$ comes via different arrow types. As these types reflect various types of interactions, this can lead to one cell only amplifying its inputs `before' the bifurcation point and the other one `after' the bifurcation point $\lambda=0$. 

On the other hand, whenever $(*)$ and $(**)$ have the same sign, there is also a suitable choice of coefficients in cells $2$ and $3$ such that $(***) <0$, as was mentioned before. Hence, there is also generically a branching solution for cell $1$ resulting in the generic existence of the solution branch with $B=\{5\}$ for the entire network. In this context, genericity means that the solution branch exists for an open but not dense set of system parameters.

\begin{figure}[h!]
	\begin{subfigure}[h]{\textwidth}
		\begin{center}
			\resizebox{.7\linewidth}{!}{	
				\input{5cellfeedforwardwithoddbifurcationsextnet-bif1}%
			}%
		\end{center}%
		\subcaption{Parameters: $a_{\sigma_1} = 0, a_{\sigma_2} = 1, a_{\sigma_3} = -2, a_{\sigma_4} = 1, a_{\sigma_5} = 1, \ell = 1, f_{\Id\Id} = -1, f_{\Id\lambda} = 1, f_{\Id\sigma_2} + f_{\Id\sigma_3} + f_{\Id\sigma_4} + f_{\Id\sigma_5} = \frac{1}{2}$. Total number of branches: 5.}
		\label{fig:ex-bif1}
	\end{subfigure} \\%
	\begin{subfigure}[h]{\textwidth}
		\vspace{10pt}
		\begin{center}
			\resizebox{.7\linewidth}{!}{	
				\input{5cellfeedforwardwithoddbifurcationsextnet-bif2}%
			}%
		\end{center}%
		\subcaption{Parameters: $a_{\sigma_1} = 0, a_{\sigma_2} = 1, a_{\sigma_3} = \frac{1}{2}, a_{\sigma_4} = -\frac{1}{2}, a_{\sigma_5} = 0, \ell = -1, f_{\Id\Id} = 1, f_{\Id\lambda} = 0, f_{\Id\sigma_2} = 0, f_{\Id\sigma_3} = 0, f_{\Id\sigma_4} = 0, f_{\Id\sigma_5} = 0$. Total number of branches: 10.}
		\label{fig:ex-bif2}
	\end{subfigure}
	\caption{Depiction of the qualitative steady state bifurcations of the network in \Cref{fig:numerics1} with different parameter values. The diagrams describe each cells behavior separately. However, the branching is not independent of the other cells as described in \Cref{sec:ex}.}
	\label{fig:ex-bif}
\end{figure}
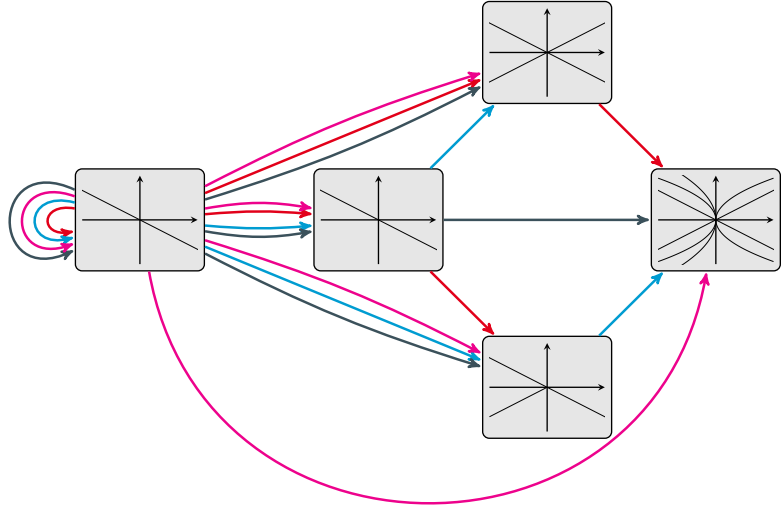
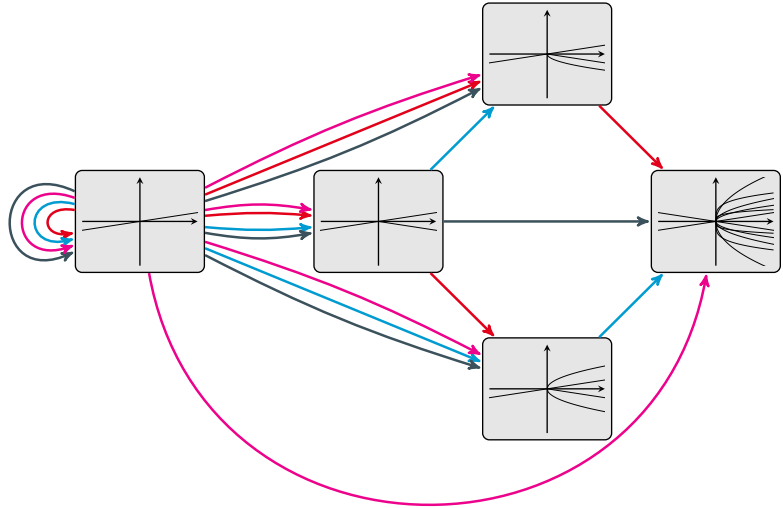

In \Cref{fig:ex-bif} we illustrate the steady state bifurcations for two different choices of parameter values. The qualitative bifurcation scenario is depicted for each cell separately. Note that for a non-maximal cell certain branches are only possible if cells above it are in a suitable state. This fact is not displayed in the figures. Both choices of parameters are generic but display different behavior. The amplification effect can be seen in both. However, in \Cref{fig:ex-bif1} the strongest amplification is $\sim \sqrt{\lambda}$ in cell $1$, whereas we also find a branch $\sim \sqrt[4]{\lambda}$ in cell $1$ in \Cref{fig:ex-bif2}.

%% file: 5cellfeedforwardwithoddbifurcationsvfextnet.tex
\begin{equation*}
\dot{x} = \gamma_{f} (x) =
\begin{pmatrix}
f(x_{1},\textcolor{red}{x_{2}},\textcolor{blue}{x_{3}},\textcolor{grey}{x_{4}},\textcolor{magenta}{x_{5}}, \lambda)\\
f(x_{2},\textcolor{red}{x_{5}},\textcolor{blue}{x_{4}},\textcolor{grey}{x_{5}},\textcolor{magenta}{x_{5}}, \lambda)\\
f(x_{3},\textcolor{red}{x_{4}},\textcolor{blue}{x_{5}},\textcolor{grey}{x_{5}},\textcolor{magenta}{x_{5}}, \lambda)\\
f(x_{4},\textcolor{red}{x_{5}},\textcolor{blue}{x_{5}},\textcolor{grey}{x_{5}},\textcolor{magenta}{x_{5}}, \lambda)\\
f(x_{5},\textcolor{red}{x_{5}},\textcolor{blue}{x_{5}},\textcolor{grey}{x_{5}},\textcolor{magenta}{x_{5}}, \lambda)\\
\end{pmatrix}
\end{equation*}

%% file: 5cellfeedforwardwithoddbifurcationsextnet-bif1.tex
\centering
\begin{tikzpicture}[->,
	>=stealth',
	shorten >=1pt,
	auto,
	node distance=5cm,
	plot node/.style={rectangle, rounded corners, thick, draw, fill=black!10}]
	\node[plot node] (1) {%
		\includegraphics{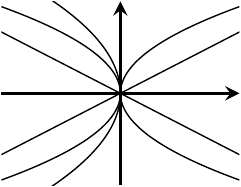}
	};
	\node[plot node, above left of=1] (2) {%
		\includegraphics{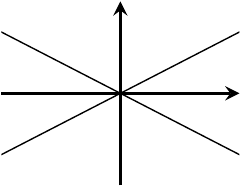}
	};
	\node[plot node, below left of=1] (3) {%
		\includegraphics{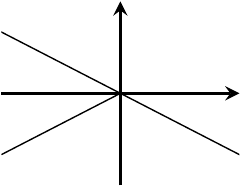}
	};
	\node[plot node, below left of=2] (4) {%
		\includegraphics{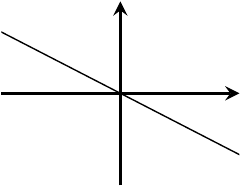}
	};
	\node[plot node, left of=4] (5) {%
		\includegraphics{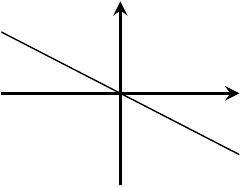}
	};
	\path[every node/.style={font=\sffamily\small}, line width =1.5pt]
	(2) edge [color = {red}] node {} (1)
	(5) edge [color = {red}] node {} (2)
	(4) edge [color = {red}] node {} (3)
	(5) edge [bend left = 5, color = {red}] node {} (4)
	(5) edge [in=190, out=170, looseness = 4, color = {red}] node {} (5)
	(3) edge [color = {blue}] node {} (1)
	(4) edge [color = {blue}] node {} (2)
	(5) edge [color = {blue}] node {} (3)
	(5) edge [bend left = -5, color = {blue}] node {} (4)
	(5) edge [in=195, out=165, looseness = 4, color = {blue}] node {} (5)
	(4) edge [color = {grey}] node {} (1)
	(5) edge [bend left =-5, color = {grey}] node {} (2)
	(5) edge [bend left =-5, color = {grey}] node {} (3)
	(5) edge [bend left =-10, color = {grey}] node {} (4)
	(5) edge [in=205, out=155, looseness = 4, color = {grey}] node {} (5)
	(5) edge [in = -100, out = -80, looseness = 1.45, color = {magenta}] node {} (1)
	(5) edge [bend left =5, color = {magenta}] node {} (2)
	(5) edge [bend left =5, color = {magenta}] node {} (3)
	(5) edge [bend left =10, color = {magenta}] node {} (4)
	(5) edge [in=200, out=160, looseness = 4, color = {magenta}] node {} (5)
;
\end{tikzpicture}

%% file: 5cellfeedforwardwithoddbifurcationsextnet-bif2.tex
\centering
\begin{tikzpicture}[->,
	>=stealth',
	shorten >=1pt,
	auto,
	node distance=5cm,
	plot node/.style={rectangle, rounded corners, thick, draw, fill=black!10}]
	\node[plot node] (1) {%
		\includegraphics{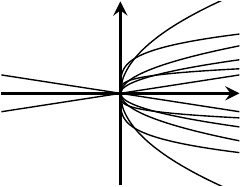}
	};
	\node[plot node, above left of=1] (2) {%
		\includegraphics{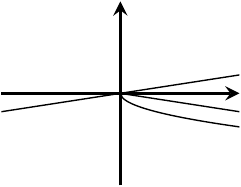}
	};
	\node[plot node, below left of=1] (3) {%
		\includegraphics{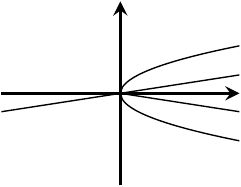}
	};
	\node[plot node, below left of=2] (4) {%
		\includegraphics{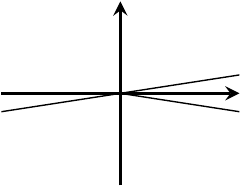}
	};
	\node[plot node, left of=4] (5) {%
		\includegraphics{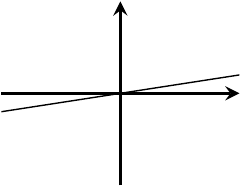}
	};
	\path[every node/.style={font=\sffamily\small}, line width =1.5pt]
	(2) edge [color = {red}] node {} (1)
	(5) edge [color = {red}] node {} (2)
	(4) edge [color = {red}] node {} (3)
	(5) edge [bend left = 5, color = {red}] node {} (4)
	(5) edge [in=190, out=170, looseness = 4, color = {red}] node {} (5)
	(3) edge [color = {blue}] node {} (1)
	(4) edge [color = {blue}] node {} (2)
	(5) edge [color = {blue}] node {} (3)
	(5) edge [bend left = -5, color = {blue}] node {} (4)
	(5) edge [in=195, out=165, looseness = 4, color = {blue}] node {} (5)
	(4) edge [color = {grey}] node {} (1)
	(5) edge [bend left =-5, color = {grey}] node {} (2)
	(5) edge [bend left =-5, color = {grey}] node {} (3)
	(5) edge [bend left =-10, color = {grey}] node {} (4)
	(5) edge [in=205, out=155, looseness = 4, color = {grey}] node {} (5)
	(5) edge [in = -100, out = -80, looseness = 1.45, color = {magenta}] node {} (1)
	(5) edge [bend left =5, color = {magenta}] node {} (2)
	(5) edge [bend left =5, color = {magenta}] node {} (3)
	(5) edge [bend left =10, color = {magenta}] node {} (4)
	(5) edge [in=200, out=160, looseness = 4, color = {magenta}] node {} (5)
;
\end{tikzpicture}

%% file: II_genfeedforward_ap.tex

\begin{appendix}
\renewcommand\thetheorem{\thesubsection.\arabic{theorem}}
\makeatletter
\@addtoreset{theorem}{subsection}
\makeatother
\renewcommand\theequation{\thesubsection.\arabic{equation}}
\makeatletter
\@addtoreset{equation}{subsection}
\makeatother
\section{Appendix}
\label{sec:ap1}
In this appendix we fill the gaps left in \Cref{sec:bi} by proving
\Cref{lem:abff_maxnoncritnonmaxnoncrit,lem:abff_maxnoncritnonmaxcritL,lem:abff_maxnoncritnonmaxcritSN} as well as the leading coefficients of branching solutions in \Cref{lem:abff_maxnoncritnonmaxcritT}. Recall that $\LL_p = \{ \sigma\in\Sigma \mid \sigma(p)=p \}$.

\begin{proof}[Proof of \Cref{lem:abff_maxnoncritnonmaxnoncrit}]
	The proofs for both cases are very similar and analogous to the proofs for \Cref{lem:c1-st2.1,lem:c1-st2.2}. Hence, we only sketch them here. We assume (\textbf{H}) and (L) first. \Cref{eq:taylor} becomes
	\[ 0 = \sum_{\sigma \in \LL_{p}} a_\sigma x_p + \sum_{\substack{p \nol q = \tau(p) \colon \\ \tau \notin \LL_p}} a_\tau \dsup_q \lambda + \ell \lambda + \OO \left( |x_p|^2 + |\lambda||x_p| + |\lambda|^2 \right). \]
	As $\sum_{\sigma \in \LL_{p}} a_\sigma \ne 0$ this is uniquely solved by 
	\[ x_p (\lambda) = -\frac{\sum_{\tau \notin \LL_p} a_\tau \dsup_{\tau(p)} + \ell}{\sum_{\sigma \in \LL_{p}} a_\sigma} \cdot \lambda + \OO\left(|\lambda|^2\right) \]
	for small $\lambda>0$, due to the implicit function theorem. Note that, because of assumption (L), the linear coefficient does not vanish.
	
	Next, assume (\textbf{H}) and (SN). \Cref{eq:taylor} becomes
	\[ 0 = \sum_{\sigma \in \LL_{p}} a_\sigma x_p + \sum_{\substack{\tau \notin \LL_p \colon \\ \tau(p) = q\in \Qsup_p}} a_\tau \dsup_q \lambda^{2^{-\xisup_p}}  + \OO \left( |x_p|^2 + |\lambda|^{2^{-\xisup_p}}|x_p| + |\lambda|^{2^{-(\xisup_p-1)}} \right). \]
	By the same argument as before, this is uniquely solved by
	\[ x_p (\lambda) = -\frac{\sum_{\tau \colon \tau(p) \in \Qsup_p} a_\tau \dsup_{\tau(p)}}{\sum_{\sigma \in \LL_{p}} a_\sigma} \cdot \lambda^{2^{-\xisub_p}} + \OO\left(|\lambda|^{2^{-(\xisup_p-1)}}\right) \]
	for small $\lambda>0$ with non-vanishing leading coefficient.
\end{proof}

\begin{proof}[Proof of \Cref{lem:abff_maxnoncritnonmaxcritL}]
	The proof is analogous to the one for \Cref{lem:c1-st1}, except for slightly different coefficients. Once again it uses the standard technique for detecting saddle node bifurcations as in \textcite{Murdock.2003}. As $\xisup_p=0$, \eqref{eq:taylor} becomes
	\[ 0 = \sum_{\sigma, \tau \in \LL_p} f_{\sigma\tau} x_p^2 + \sum_{\substack{\tau \not\in \LL_p \\ p \nol q = \tau(p)}} a_\tau \dsup_q \lambda + \ell \lambda + \OO\left( |x_p|^3 + |\lambda||x_p| + |\lambda|^2 \right). \]
	We introduce a new variable $x_p = \mu y$ where $\mu = \sqrt{\lambda}$ for small $\lambda>0$. The equation to be solved transforms into
	\[ 0 = \sum_{\sigma, \tau \in \LL_p} f_{\sigma\tau} \mu^2 y^2 + \sum_{\tau \not\in \LL_p} a_\tau d_{\tau(p)} \mu^2 + \ell \mu^2 + \OO\left( |\mu|^3|y| + |\mu|^4 \right). \]
	As $\mu>0$, we may divide by $\mu^2$ and obtain
	\[ 0 = \sum_{\sigma, \tau \in \LL_p} f_{\sigma\tau} y^2 + \sum_{\tau \not\in \LL_p} a_\tau \dsup_{\tau(p)} + \ell + \OO\left( |\mu||y| + |\mu|^2 \right) = g(y,\mu). \]
	If $\left( \sum_{\tau \notin \LL_p} a_\tau \dsup_{\tau(p)} + \ell \right) / \sum_{\sigma,\tau \in \LL_{p}} f_{\sigma\tau} >0$, the equation $g(y,0)=0$ has no real solutions. If, on the other hand, $\left( \sum_{\tau \notin \LL_p} a_\tau \dsup_{\tau(p)} + \ell \right) / \sum_{\sigma,\tau \in \LL_{p}} f_{\sigma\tau} <0$, there are two solutions to $g(y,0)=0$
	\[ \overline{y}^\pm = \pm \sqrt{- \frac{\sum_{\tau \not\in \LL_p} a_\tau \dsup_{\tau(p)} + \ell}{\sum_{\sigma, \tau \in \LL_p} f_{\sigma\tau}}}. \]
	Furthermore, $\frac{\partial}{\partial y} g(\overline{y}^\pm, 0) = 2 \sum_{\sigma, \tau \in \LL_p} f_{\sigma\tau} \overline{y}^\pm$, which does not vanish, due to assumption (L). Hence, by the implicit function theorem, we obtain two branches of solutions
	\[ \overline{Y}^\pm (\mu) = \overline{y}^\pm + \OO(|\mu|). \]
	Transforming back into the original variables, we obtain the two branches
	\[ x_p (\lambda) = \pm \sqrt{- \frac{\sum_{\tau \not\in \LL_p} a_\tau \dsup_{\tau(p)} + \ell}{\sum_{\sigma, \tau \in \LL_p} f_{\sigma\tau}}} \sqrt{\lambda} + \OO(|\lambda|) \]
	for small $\lambda>0$.
\end{proof}

\begin{proof}[Proof of \Cref{lem:abff_maxnoncritnonmaxcritSN}]
	Under the given assumptions \eqref{eq:taylor} becomes
	\[ 0 = \sum_{\sigma, \tau \in \LL_p} f_{\sigma\tau} x_p^2 + \sum_{\substack{\tau \notin \LL_p \\ \tau(p) = q \in \Qsup_p}} a_\tau \dsup_q \lambda^{2^{-\xisup_p}} + \OO\left( |x_p|^3 + |x_p||\lambda|^{2^{-\xisup_p}} + |\lambda|^{2^{-(\xisup_p-1)}} \right). \]
	Similar to previous proofs, we introduce new coordinates $x_p = \mu y$, where \mbox{$\mu = \sqrt{\lambda^{2^{-\xisup_p}}} = \lambda^{2^{-(\xisup_p+1)}}$} for small $\lambda>0$. The equation becomes
	\[ 0 = \sum_{\sigma, \tau \in \LL_p} f_{\sigma\tau} \mu^2 y^2 + \sum_{\substack{\tau \notin \LL_p \\ \tau(p) = q \in \Qsup_p}} a_\tau \dsup_q \mu^2 + \OO\left( |y||\mu|^3 + |\mu|^4 \right). \]
	As $\mu >0$, we may divide by $\mu^2$ to obtain
	\[ 0 = \sum_{\sigma, \tau \in \LL_p} f_{\sigma\tau} y^2 + \sum_{\substack{\tau \notin \LL_p \\ \tau(p) = q \in \Qsup_p}} a_\tau \dsup_q + \OO\left( |y||\mu| + |\mu|^2 \right) = g(y,\mu). \]
	If $\left( \sum_{\tau \colon \tau(p) \in \Qsup_p} a_\tau \dsup_{\tau(p)} \right) / \sum_{\sigma, \tau \in \LL_p} f_{\sigma\tau} > 0$ there are no solutions to $g(y,0)=0$ -- this proves the first case. If, on the other hand, $\left( \sum_{\tau \colon \tau(p) \in \Qsup_p} a_\tau \dsup_{\tau(p)} \right) / \sum_{\sigma, \tau \in \LL_p} f_{\sigma\tau} < 0$, there are two solutions to $g(y,0)=0$, as $\sum_{\sigma, \tau \in \LL_p} f_{\sigma\tau} \ne 0$ generically:
	\[ \overline{y}^\pm = \pm \sqrt{- \frac{\sum_{\tau \colon \tau(p) \in \Qsup_p} a_\tau \dsup_{\tau(p)}}{\sum_{\sigma, \tau \in \LL_p} f_{\sigma\tau}}}. \]
	Furthermore, $\frac{\partial}{\partial y} g(\overline{y}^\pm, 0) = 2 \sum_{\sigma, \tau \in \LL_p} f_{\sigma\tau} \overline{y}^\pm$, which, by the same argument, generically does not vanish. Hence, by the implicit function theorem, we obtain two branches of solutions
	\[ \overline{Y}^\pm (\mu) = \overline{y}^\pm + \OO(|\mu|). \]
	Transforming back into the original coordinates, we obtain
	\[ x_p (\lambda) = \pm \sqrt{- \frac{\sum_{\tau \colon \tau(p) \in \Qsup_p} a_\tau \dsup_{\tau(p)}}{\sum_{\sigma, \tau \in \LL_p} f_{\sigma\tau}}} \cdot \lambda^{2^{-(\xisup_p+1)}} + \OO\left( |\lambda|^{2^{-\xisup_p}} \right), \]
	for small $\lambda>0$, which completes the proof.
\end{proof}


\begin{lemma}
	\label{lem:abff_maxnoncritDcritsuper}
	Let $p \in C$ be non-maximal and critical. Assume \emph{(\textbf{H})} with
	\[ \dsup_q = -\frac{\ell}{\sum_{\sigma \in \Sigma} a_\sigma}, \qquad \Rsup_q = - \frac{\sum_{\sigma, \tau \in \Sigma} f_{\sigma \tau} \ell^2 - \sum_{\sigma\in\Sigma} a_\sigma \sum_{\sigma \in \Sigma} f_{\sigma\lambda} \ell + \left( \sum_{\sigma \in \Sigma} a_\sigma \right)^2 f_{\lambda\lambda}}{\left( \sum_{\sigma \in \Sigma} a_\sigma \right)^3} \]
	for all $q \nog p$ and define
	\begin{align*}
		A	&= \sum_{\sigma, \tau \in \LL_p} f_{\sigma\tau}, \\
		B	&= \sum_{\sigma \in \LL_{p}} f_{\sigma \lambda} + 2 \sum_{\sigma \in \LL_p, \tau \notin \LL_p } f_{\sigma \tau} \dsup_{\tau(p)}, \\
		C	&= \sum_{\tau \notin \LL_p} a_\tau \Rsup_{\tau(p)} + \sum_{\tau \notin \LL_p} f_{\tau \lambda} \dsup_{\tau(p)} + \sum_{\sigma, \tau \notin \LL_p} f_{\sigma \tau} \dsup_{\sigma(p)} \dsup_{\tau(p)} + f_{\lambda \lambda}, \\
		E	&= \sum_{\sigma \in \LL_p} f_{\sigma\lambda} - 2 \cdot \frac{\ell}{\sum_{\sigma\in\Sigma} a_\sigma} \sum_{\substack{\sigma \in \Sigma \\ \tau \in \LL_p}} f_{\sigma \tau}.
	\end{align*}
	Then generically
	\[ B^2 - 4AC = E^2 > 0 \]
	and
	\[ \frac{-B + E}{2A} = -\frac{\ell}{\sum_{\sigma \in \Sigma} a_\sigma}, \qquad \frac{-B - E}{2A} = M \]
	with 
	\[ M = \dfrac{\ell}{\sum_{\sigma \in \Sigma} a_\sigma} \cdot \left( 1 + 2 \dfrac{\sum_{\sigma \in \LL_p, \tau \notin \LL_p} f_{\sigma\tau}}{\sum_{\sigma, \tau \in \LL_p} f_{\sigma\tau}} \right) - \dfrac{\sum_{\sigma \in \LL_p} f_{\sigma\lambda}}{\sum_{\sigma, \tau \in \LL_p} f_{\sigma\tau}}. \]
	In particular, these equalities hold for the coefficients $d_p^\pm$ in \Cref{lem:abff_maxnoncritnonmaxcritT} if all $q \nog p$ are in the fully synchronous steady state.
\end{lemma}
\begin{proof}
	Let $p \in C$ be non-maximal. Assume
	\[ \dsup_q = - \frac{\ell}{\sum_{\sigma \in \Sigma} a_\sigma}	\quad \text{and} \quad \Rsup_q = - \frac{\sum_{\sigma, \tau \in \Sigma} f_{\sigma \tau} \ell^2 - \sum_{\sigma\in\Sigma} a_\sigma \sum_{\sigma \in \Sigma} f_{\sigma\lambda} \ell + \left( \sum_{\sigma \in \Sigma} a_\sigma \right)^2 f_{\lambda\lambda}}{\left( \sum_{\sigma \in \Sigma} a_\sigma \right)^3} \]
	for all $q \nog p$. A key observation is
	\[ \sum_{\tau \notin \LL_p} a_\tau = \sum_{\sigma \in \Sigma} a_\sigma, \]
	as $\sum_{\sigma \in \LL_{p}} a_\sigma = 0$. We denote this sum by $K$. Hence,
	\begin{align*}
		\sum_{\tau \notin \LL_p} a_\tau \Rsup_{\tau(p)}	&= - \frac{\sum_{\sigma, \tau \in \Sigma} f_{\sigma \tau} \ell^2 - K \sum_{\sigma \in \Sigma} f_{\sigma\lambda} \ell + K^2 f_{\lambda\lambda}}{K^2} \\
		& = - \frac{\ell^2}{K^2} \sum_{\sigma, \tau \in \Sigma} f_{\sigma\tau} + \frac{\ell}{K} \sum_{\sigma \in \Sigma} f_{\sigma\lambda} - f_{\lambda\lambda}.
	\end{align*}
	Note that
	\[ \sum_{\sigma, \tau \in \Sigma} f_{\sigma\tau} = \sum_{\sigma, \tau \in \LL_p} f_{\sigma\tau} + 2 \sum_{\substack{\sigma \in \LL_p \\ \tau \notin \LL_p}} f_{\sigma\tau} + \sum_{\sigma, \tau \notin \LL_p} f_{\sigma\tau}, \]
	where we have used $f_{\sigma\tau}=f_{\tau\sigma}$. Similar considerations occur frequently in the remainder of this proof. We use them without explicitly mentioning them. Furthermore, we compute
	\begin{align*}
		\sum_{\tau \notin \LL_p} f_{\tau \lambda} \dsup_{\tau(p)}							&= - \frac{\ell}{K} \sum_{\tau \notin \LL_p} f_{\tau \lambda}, \\
		\sum_{\sigma, \tau \notin \LL_p} f_{\sigma \tau} \dsup_{\sigma(p)} \dsup_{\tau(p)}	&= \frac{\ell^2}{K^2} \sum_{\sigma, \tau \notin \LL_p} f_{\sigma \tau} .
	\end{align*} 
	Thus, we obtain
	\begin{align*}
		C	&= \sum_{\tau \notin \LL_p} a_\tau \Rsup_{\tau(p)} + \sum_{\tau \notin \LL_p} f_{\tau \lambda} \dsup_{\tau(p)} + \sum_{\sigma, \tau \notin \LL_p} f_{\sigma \tau} \dsup_{\sigma(p)} \dsup_{\tau(p)} + f_{\lambda \lambda} \\
		&= - \frac{\ell^2}{K^2} \sum_{\sigma, \tau \in \Sigma} f_{\sigma\tau} + \frac{\ell}{K} \sum_{\sigma \in \Sigma} f_{\sigma\lambda} - f_{\lambda\lambda} - \frac{\ell}{K} \sum_{\tau \notin \LL_p} f_{\tau \lambda} + \frac{\ell^2}{K^2} \sum_{\sigma, \tau \notin \LL_p} f_{\sigma \tau} + f_{\lambda\lambda} \\
		&= - \frac{\ell^2}{K^2} \left( \sum_{\sigma, \tau \in \LL_p} f_{\sigma\tau} + 2 \sum_{\substack{\sigma \in \LL_p \\ \tau \notin \LL_p}} f_{\sigma\tau} \right)  + \frac{\ell}{K} \sum_{\sigma \in \LL_p} f_{\sigma\lambda}.
	\end{align*}
	Next, we compute
	\begin{align*}
		B^2	&= \left( \sum_{\sigma \in \LL_{p}} f_{\sigma \lambda} + 2 \sum_{\substack{\sigma \in \LL_p \\ \tau \notin \LL_p}} f_{\sigma \tau} \dsup_{\tau(p)} \right)^2 \\
		&= \left( \sum_{\sigma \in \LL_{p}} f_{\sigma \lambda} - 2 \cdot \frac{\ell}{K} \sum_{\substack{\sigma \in \Sigma \\ \tau \in \LL_p}} f_{\sigma \tau} + 2 \cdot \frac{\ell}{K} \sum_{\sigma, \tau \in \LL_p} f_{\sigma \tau} \right)^2 \\
		&= \left( \sum_{\sigma \in \LL_p} f_{\sigma\lambda} - 2 \cdot \frac{\ell}{K} \sum_{\substack{\sigma \in \Sigma \\ \tau \in \LL_p}} f_{\sigma \tau} \right)^2 + 4 \cdot L
	\end{align*}
	with
	\begin{align*}
		L	&= \frac{\ell}{K} \left( \sum_{\sigma \in \LL_p} f_{\sigma\lambda} - 2 \cdot \frac{\ell}{K} \sum_{\substack{\sigma \in \Sigma \\ \tau \in \LL_p}} f_{\sigma \tau} \right) \sum_{\sigma,\tau \in \LL_p} f_{\sigma \tau} + \frac{\ell^2}{K^2} \left( \sum_{\sigma,\tau \in \LL_p} f_{\sigma \tau} \right)^2 \\
		& = \sum_{\sigma,\tau \in \LL_p} f_{\sigma \tau} \cdot \left( \frac{\ell}{K} \sum_{\sigma \in \LL_p} f_{\sigma\lambda} - \frac{\ell^2}{K^2} \left( 2 \cdot \sum_{\substack{\sigma \in \Sigma \\ \tau \in \LL_p}} f_{\sigma \tau} - \sum_{\sigma,\tau \in \LL_p} f_{\sigma \tau} \right) \right) \\
		& = \sum_{\sigma,\tau \in \LL_p} f_{\sigma \tau} \cdot \left( \frac{\ell}{K} \sum_{\sigma \in \LL_p} f_{\sigma\lambda} - \frac{\ell^2}{K^2} \left( \sum_{\sigma, \tau \in \LL_p} f_{\sigma \tau} + 2 \cdot \sum_{\substack{\sigma \in \LL_p \\ \tau \notin \LL_p}} f_{\sigma \tau} \right) \right) = A \cdot C.
	\end{align*}
	
	Hence,
	\begin{align*}
		B^2 - 4 AC &= \left( \sum_{\sigma \in \LL_p} f_{\sigma\lambda} - 2 \cdot \frac{\ell}{K} \sum_{\substack{\sigma \in \Sigma \\ \tau \in \LL_p}} f_{\sigma \tau} \right)^2 + 4L - 4AC \\
		&= \left( \sum_{\sigma \in \LL_p} f_{\sigma\lambda} - 2 \cdot \frac{\ell}{K} \sum_{\substack{\sigma \in \Sigma \\ \tau \in \LL_p}} f_{\sigma \tau} \right)^2 \\
		&= E^2.
	\end{align*}
	Generically, this expression is positive, which allows us to compute
	\begin{align*}
		- B + E	&= 2 \cdot \frac{\ell}{K} \left( \sum_{\substack{\sigma \in \LL_p \\ \tau \notin \LL_p}} f_{\sigma \tau} - \sum_{\substack{\sigma \in \Sigma \\ \tau \in \LL_p}} f_{\sigma \tau} \right) \\
		&= - 2 \cdot \frac{\ell}{K} \sum_{\sigma, \tau \in \LL_p} f_{\sigma \tau} = - 2 A \cdot \frac{\ell}{K},
	\end{align*}
	proving
	\[ \frac{-B + E}{2A} = - \frac{\ell}{\sum_{\sigma \in \Sigma} a_\sigma}. \]
	
	On the other hand
	\begin{align*}
		- B - E	&= - \sum_{\sigma \in \LL_{p}} f_{\sigma \lambda} + 2 \cdot \frac{\ell}{K} \sum_{\substack{\sigma \in \LL_p \\ \tau \notin \LL_p}} f_{\sigma \tau} - \left( \sum_{\sigma \in \LL_p} f_{\sigma\lambda} - 2 \cdot \frac{\ell}{K} \sum_{\substack{\sigma \in \Sigma \\ \tau \in \LL_p}} f_{\sigma \tau} \right) \\
		&= 2 \cdot \left( \frac{\ell}{K} \left( \sum_{\sigma, \tau \in \LL_p} f_{\sigma \tau} + 2 \sum_{\substack{\sigma \in \LL_p \\ \tau \notin \LL_p}} f_{\sigma \tau} \right) - \sum_{\sigma \in \LL_{p}} f_{\sigma \lambda} \right),
	\end{align*}
	proving
	\[ \frac{- B - E}{2A} = \dfrac{\ell}{\sum_{\sigma \in \Sigma} a_\sigma} \cdot \left( 1 + 2 \dfrac{\sum_{\sigma \in \LL_p, \tau \notin \LL_p} f_{\sigma\tau}}{\sum_{\sigma, \tau \in \LL_p} f_{\sigma\tau}} \right) - \dfrac{\sum_{\sigma \in \LL_p} f_{\sigma\lambda}}{\sum_{\sigma, \tau \in \LL_p} f_{\sigma\tau}} = M. \]
\end{proof}
\end{appendix}